\documentclass[11pt]{article}

\usepackage{amssymb,amsbsy,amsthm,amsmath,graphicx,mathrsfs,epsfig,graphics,times}

\input xy

\xyoption{all}

\begin{document}

\title{\textbf{An alternative ending to ``Pleasant extensions retaining algebraic structure''}}
\author{Tim Austin}
\date{}

\maketitle


\newenvironment{nmath}{\begin{center}\begin{math}}{\end{math}\end{center}}

\newtheorem{thm}{Theorem}[section]
\newtheorem*{thm*}{Theorem}
\newtheorem{lem}[thm]{Lemma}
\newtheorem{prop}[thm]{Proposition}
\newtheorem{cor}[thm]{Corollary}
\newtheorem{conj}[thm]{Conjecture}
\newtheorem{dfn}[thm]{Definition}
\newtheorem{prob}[thm]{Problem}
\newtheorem{ques}[thm]{Question}
\theoremstyle{remark}
\newtheorem*{ex}{Example}


\newcommand{\s}{\sigma}
\renewcommand{\O}{\Omega}
\renewcommand{\S}{\Sigma}
\newcommand{\co}{\mathrm{co}}
\newcommand{\e}{\mathrm{e}}
\newcommand{\eps}{\varepsilon}
\renewcommand{\d}{\mathrm{d}}
\newcommand{\im}{\mathrm{i}}
\renewcommand{\k}{\kappa}
\renewcommand{\l}{\lambda}
\newcommand{\G}{\Gamma}
\newcommand{\g}{\gamma}
\renewcommand{\L}{\Lambda}
\renewcommand{\a}{\alpha}
\renewcommand{\b}{\beta}
\newcommand{\Sone}{\mathrm{S}^1}

\newcommand{\Aut}{\mathrm{Aut}}
\renewcommand{\Pr}{\mathrm{Pr}}
\newcommand{\Hom}{\mathrm{Hom}}
\newcommand{\id}{\mathrm{id}}
\newcommand{\BSp}{\mathsf{BSp}}
\newcommand{\WBSp}{\mathsf{WBSp}}
\newcommand{\Sys}{\underline{\mathsf{Sys}}}
\newcommand{\Rep}{\underline{\mathsf{Rep}}}
\newcommand{\Lat}{\mathrm{Lat}}
\newcommand{\FLat}{\mathrm{FLat}}
\newcommand{\SG}{\mathrm{SG}}
\newcommand{\CSG}{\mathrm{CSG}}
\newcommand{\Clos}{\mathrm{Clos}}
\newcommand{\pro}{\mathrm{pro}}
\newcommand{\nil}{\mathrm{nil}}
\newcommand{\dCL}{\mathrm{dCL}}
\newcommand{\rat}{\mathrm{rat}}
\newcommand{\Ab}{\mathrm{Ab}}
\newcommand{\CIS}{\mathrm{CIS}}

\newcommand{\bbN}{\mathbb{N}}
\newcommand{\bbR}{\mathbb{R}}
\newcommand{\bbZ}{\mathbb{Z}}
\newcommand{\bbQ}{\mathbb{Q}}
\newcommand{\bbT}{\mathbb{T}}
\newcommand{\bbD}{\mathbb{D}}
\newcommand{\bbC}{\mathbb{C}}
\newcommand{\bbP}{\mathbb{P}}

\newcommand{\A}{\mathcal{A}}
\newcommand{\B}{\mathcal{B}}
\newcommand{\C}{\mathcal{C}}
\newcommand{\E}{\mathcal{E}}
\newcommand{\F}{\mathcal{F}}
\newcommand{\I}{\mathcal{I}}
\newcommand{\calL}{\mathcal{L}}
\renewcommand{\P}{\mathcal{P}}
\newcommand{\U}{\mathcal{U}}
\newcommand{\W}{\mathcal{W}}
\newcommand{\Y}{\mathcal{Y}}
\newcommand{\Z}{\mathcal{Z}}

\newcommand{\frH}{\mathfrak{H}}
\newcommand{\frM}{\mathfrak{M}}

\newcommand{\bfX}{\mathbf{X}}
\newcommand{\bfY}{\mathbf{Y}}
\newcommand{\bfZ}{\mathbf{Z}}
\newcommand{\bfW}{\mathbf{W}}
\newcommand{\bfV}{\mathbf{V}}

\newcommand{\rmC}{\mathrm{C}}
\newcommand{\rmH}{\mathrm{H}}
\newcommand{\T}{\mathrm{T}}

\newcommand{\sfC}{\mathsf{C}}
\newcommand{\sfE}{\mathsf{E}}
\newcommand{\sfV}{\mathsf{V}}
\newcommand{\sfW}{\mathsf{W}}
\newcommand{\sfY}{\mathsf{Y}}
\newcommand{\sfZ}{\mathsf{Z}}

\newcommand{\uhr}{\!\!\upharpoonright}
\newcommand{\into}{\hookrightarrow}
\newcommand{\onto}{\twoheadrightarrow}

\newcommand{\bb}[1]{\mathbb{#1}}
\newcommand{\bs}[1]{\boldsymbol{#1}}
\newcommand{\fr}[1]{\mathfrak{#1}}
\renewcommand{\bf}[1]{\mathbf{#1}}
\renewcommand{\sf}[1]{\mathsf{#1}}
\renewcommand{\rm}[1]{\mathrm{#1}}
\renewcommand{\cal}[1]{\mathcal{#1}}

\renewcommand{\t}[1]{\tilde{#1}}

\newcommand{\mvee}{\hbox{$\bigvee$}}

\newcommand{\lfr}{\lfloor}
\newcommand{\rfr}{\rfloor}

\newcommand{\fin}{\nolinebreak\hspace{\stretch{1}}$\lhd$}
\newcommand{\tick}[1]{\nolinebreak\hspace{\stretch{1}}$\surd_{\mathrm{#1}}$}

\begin{abstract}
The culmination of the two recent papers~\cite{Aus--lindeppleasant1,Aus--lindeppleasant2} was a proof of the norm convergence in $L^2(\mu)$ of the quadratic nonconventional ergodic averages
\[\frac{1}{N}\sum_{n=1}^N(f_1\circ T_1^{n^2})(f_2\circ T_1^{n^2}T_2^n)\quad\quad f_1,f_2\in L^\infty(\mu)\]
associated to an arbitrary probability-preserving $\bbZ^2$-system $(X,\mu,T_1,T_2)$. This is a special
case of the Bergelson-Leibman conjecture on the norm convergence of
polynomial nonconventional ergodic averages~\cite{BerLei02}.

That proof relied on some new
machinery for extending probability-preserving $\bbZ^d$-systems to
obtain simplified asymptotic behaviour for various nonconventional averages such as the above.  The engine of this machinery is formed by some detailed structure theorems for the `characteristic factors' that are available for some such averages after ascending to a suitably-extended system. However, these new structure theorems underwent two distinct phases of development, separated by the discovery of some new technical results in Moore's cohomology theory for locally compact groups~\cite{Aus--cohomcty}.  That discovery enabled a significant improvement to the main structure theorem (Theorem 1.1 in~\cite{Aus--lindeppleasant1}), which in turn afforded a much shortened proof of convergence.  However, since the proof of convergence using the original structure theorem required some quite different ideas that are now absent from~\cite{Aus--lindeppleasant1,Aus--lindeppleasant2}, I have recorded it here in case it has some independent interest.
\end{abstract}

\parskip 0pt

\tableofcontents

\parskip 7pt

\section{Introduction}

This note records a proof of a new instance of the
Bergelson-Leibman Conjecture on norm convergence of polynomial
nonconventional ergodic averages:

\begin{thm}\label{thm:polyconv} If $T_1,T_2:\bbZ\curvearrowright
(X,\mu)$ are commuting invertible probability-preserving
transformations of a standard Borel probability space then the
averages
\[\frac{1}{N}\sum_{n=1}^N (f_1\circ T_1^{n^2})(f_2\circ T_1^{n^2}T_2^n)\]
converge in $L^2(\mu)$ as $N\to\infty$ for any $f_1,f_2 \in
L^\infty(\mu)$.
\end{thm}

The proof of the present paper has been superseded by an improved approach in~\cite{Aus--lindeppleasant1,Aus--lindeppleasant2}, enabled by a recent development in the cohomology of compact groups~(\cite{Aus--cohomcty}).  Originally, the proof recorded below was contained in a Part III to the sequence~\cite{Aus--lindeppleasant1,Aus--lindeppleasant2}, and I have maintained a presentation of it here in case it has any independent interest.

The proof of Theorem~\ref{thm:polyconv} follows a strategy that has emerged by stages in
work of Furstenberg~\cite{Fur77}, Conze and
Lesigne~\cite{ConLes84,ConLes88.1,ConLes88.2}, Furstenberg and
Weiss~\cite{FurWei96}, Host and Kra~\cite{HosKra01,HosKra05},
Ziegler~\cite{Zie07} and a number of others, and in the
papers~\cite{Aus--nonconv,Aus--newmultiSzem,Aus--lindeppleasant1,Aus--lindeppleasant2}
(see the introduction to~\cite{Aus--lindeppleasant1} for a more
complete history).  We seek an extension of an initially-given
system $(X,\mu,T_1,T_2)$, say
$\pi:(\t{X},\t{\mu},\t{T}_1,\t{T}_2)\to (X,\mu,T_1,T_2)$, such that
for the extended system the analogous nonconventional averages admit
a `simple' pair of factors
$\xi_i:(\t{X},\t{\mu},\t{T}_1,\t{T}_2)\to
(Y_i,\nu_i,S_{i,1},S_{i,2})$ that is \textbf{characteristic}, in
that
\[\frac{1}{N}\sum_{n=1}^N (f_1\circ \t{T}_1^{n^2})(f_2\circ \t{T}_1^{n^2}\t{T}_2^n) \sim \frac{1}{N}\sum_{n=1}^N (\sfE_{\t{\mu}}(f_1\,|\,\xi_1)\circ \t{T}_1^{n^2})(\sfE_{\t{\mu}}(f_2\,|\,\xi_2)\circ \t{T}_1^{n^2}\t{T}_2^n)\]
as $N\to\infty$ for all $f_1,f_2 \in L^\infty(\t{\mu})$, where we
write $f_N \sim g_N$ to denote that $\|f_N - g_N\|_2 \to 0$ as
$N\to\infty$.

These factors reduce our problem to proving convergence in case each
$f_i$ is $\xi_i$-measurable. Informally we refer to an extension that admits useful
characteristic factors for some averages as a \textbf{pleasant}
extension for those averages. The construction of a pleasant
extension in this paper will rely on some of the results from~\cite{Aus--lindeppleasant1,Aus--lindeppleasant2} (or, more precisely, from the incarnations of \cite{Aus--lindeppleasant1,Aus--lindeppleasant2} from before the above-mentioned re-write).  In~\cite{Aus--lindeppleasant1} we set up some general
machinery for constructing extensions of abstract
probability-preserving systems, showing in particular how to obtain
the useful property of satedness with respect to an idempotent class
of systems.  In~\cite{Aus--lindeppleasant2} we brought this
machinery to bear on the problem of obtaining pleasant extensions
for the linear nonconventional averages
\[\frac{1}{N}\sum_{n=1}^N(f_1\circ T^{n\bf{p}_1})(f_2\circ T^{n\bf{p}_2})(f_3\circ T^{n\bf{p}_3})\]
associated to a system $T:\bbZ^2\curvearrowright (X,\mu)$ and a
triple of distinct directions $\bf{p}_1$, $\bf{p}_2$, $\bf{p}_3 \in
\bbZ^2$ that lie in general position with the origin.  The
difficulty of that construction results from the insistence that the
pleasant extension should retain the algebraic relations among the
transformations $T^{\bf{p}_i}$ that must follow from the linear
dependence of the $\bf{p}_i$. In the previous version of~\cite{Aus--lindeppleasant2} we eventually obtained a description of some characteristic factors for these linear averages
that can be secured in an extended system as joins of various
isotropy factors and a two-step distal $\bbZ^2$-system with compact
Abelian fibres of a special kind called a `directional CL-system'.

\begin{thm}[Pleasant extensions for general-position triple linear averages]\label{thm:char-three-lines-in-2D}
For each $\bf{p}_1$, $\bf{p}_2$, $\bf{p}_3 \in \bbZ^2$ that are in
general position with the origin, let $(m_i,m_{ij},m_{ik})$ be a
relatively prime triple of nonzero integers such that $m_i\bf{p}_i +
m_{ij}(\bf{p}_i - \bf{p}_j) + m_{ik}(\bf{p}_i - \bf{p}_k) = \bs{0}$.
Then any system $T:\bbZ^2\curvearrowright (X,\mu)$ has an extension
$\pi:(\t{X},\t{\mu},\t{T}) \to (X,\mu,T)$ in which for every choice
of such $\bf{p}_1,\bf{p}_2,\bf{p}_3$ the averages
\[\frac{1}{N}\sum_{n=1}^N(f_1\circ \t{T}^{n\bf{p}_1})(f_2\circ \t{T}^{n\bf{p}_2})(f_3\circ \t{T}^{n\bf{p}_3}),\quad\quad f_1,f_2,f_3 \in L^\infty(\t{\mu}),\]
admit a characteristic triple of factors $\t{\xi}_i$, $i=1,2,3$, of
the form
\[\t{\xi}_i = \zeta_0^{\t{T}^{\bf{p}_i}}\vee \zeta_0^{\t{T}^{\bf{p}_i}=\t{T}^{\bf{p}_j}}\vee\zeta_0^{\t{T}^{\bf{p}_i}=\t{T}^{\bf{p}_k}}\vee \eta_i\]
where the target of $\eta_i$ is a $(\bf{p}_i,m_{ij}(\bf{p}_i -
\bf{p}_j),m_{ik}(\bf{p}_i - \bf{p}_k))$-directional CL-system (so
certainly a two-step Abelian system) when $\{i,j,k\} = \{1,2,3\}$.
\end{thm}

The definition of directional CL-systems will be given in
Subsection~\ref{subs:dCL-again} below.  The above theorem no longer appears in~\cite{Aus--lindeppleasant2}, because it was subsequently discovered that a cohomological argument using the new continuity results for Moore cohomology in~\cite{Aus--cohomcty} enabled an arbitrary directional CL-system to be factorized into further isotropy factors and a two-step pro-nilsystem.  This leads to an improved version of the above theorem in which $\eta_i$ may itself simply be taken to be a two-step $\bbZ^2$-pro-nilsystem, and this improvement in turn leads to a much-shortened proof of convergence.  The improved structural result now appears as Theorem 1.1 in~\cite{Aus--lindeppleasant2}, and the new proof of convergence is given in Section 5 of that paper.  However, the theorem above can still be quite quickly deduced from the arguments that appear in~\cite{Aus--lindeppleasant2}: Lemma 4.35 of~\cite{Aus--lindeppleasant2} provides solutions to the `directional CL-equations', and given this a fairly simple modification of the arguments from the current Subsection 4.6 of~\cite{Aus--lindeppleasant2} yields a proof of the above structure theorem in place of its newer improvement.

The purpose of the present note is to retain a record of the proof of Theorem~\ref{thm:polyconv} using Theorem~\ref{thm:char-three-lines-in-2D} above.  Since Theorem~\ref{thm:char-three-lines-in-2D} still lies within easy reach from the new contents of~\cite{Aus--lindeppleasant2}, I will simply assume it here.

A top-level outline of the older proof proceeds as follows. From a
careful study of the possible joinings among directional CL-systems, we will
be able to obtain a rather stronger characteristic-factor result for
our advertised nonconventional quadratic averages.  In part, the
extra strength of this result will derive from a reduction to
considering $(\G,\bf{n}_2,\bf{n}_3)$-directional CL-systems for a
pair of directions $\bf{n}_2$, $\bf{n}_3 \in \bbZ^2$ and a
finite-index subgroup $\G\leq \bbZ^2$, rather than
$(\bf{n}_1,\bf{n}_2,\bf{n}_3)$-directional CL-systems for a single
direction $\bf{n}_1$.

In terms of these systems our pleasant extensions for our polynomial
averages are as follows.

\begin{thm}[Pleasant extensions for quadratic averages]\label{thm:char-poly}
Any ergodic system of two commuting transformations
$T_1,T_2\curvearrowright (X,\mu)$ has an ergodic extension
$\pi:(\t{X},\t{\mu},\t{T}_1,\t{T}_2) \to (X,\mu,T_1,T_2)$ in which
the averages
\[\frac{1}{N}\sum_{n=1}^N(f_1\circ T_1^{n^2})(f_2\circ T_1^{n^2}T_2^n)\]
admit characteristic factors of the form
\[\xi_1 = \xi_2 := \bigvee_{m\geq 1}\zeta_0^{T^m_1}\vee \zeta_0^{T_2}\vee \bigvee_{h\geq 1}\eta_h,\]
where $\eta_h$ is a factor of $(\t{X},\t{\mu},\t{T}_1,\t{T}_2)$
whose target is a $(h\bbZ^2,(h,0),(0,h))$-directional CL-system for
the finite-index sublattice $h\bbZ^2 := \{(hm,hn):\ m,n\in\bbZ\}$.
\end{thm}

We will bring Theorem~\ref{thm:char-three-lines-in-2D} to bear on
proving Theorem~\ref{thm:char-poly} via the well-known van der
Corput estimate. Note that, unlike
Theorem~\ref{thm:char-three-lines-in-2D}, we will prove
Theorem~\ref{thm:char-poly} only for ergodic systems, and obtain
ergodic extensions as a result.  In fact the proof we give works
equally well without this additional requirement, but the version
formulated above will be more convenient for our proof of
convergence.

After proving Theorem~\ref{thm:char-poly}, we proceed towards the
proof of Theorem~\ref{thm:polyconv} through a careful analysis of
how functions measurable with respect to the factor $\xi_1 = \xi_2$
above behave upon composition with powers of $T_1$ and $T_2$.
Although our methods for controlling the images of functions upon
iterating an $(h\bbZ^2,(h,0),(0,h))$-directional CL-system are 
rather clumsy, we will find that the simplification afforded by
Theorem~\ref{thm:char-poly} is still enough to enable a more-or-less
direct proof of Theorem~\ref{thm:polyconv}.  In the present note this relies on adapting a strategy developed by Host
and Kra in~\cite{HosKra01} for the treatment of the triple linear
averages $\frac{1}{N}\sum_{n=1}^N(f_1\circ T^n)(f_2\circ
T^{2n})(f_3\circ T^{3n})$ for a single transformation $T$.

\textbf{Notational remark}\quad In this note we will make free use of notations and definitions introduced in~\cite{Aus--lindeppleasant1}
and~\cite{Aus--lindeppleasant2}.
\fin

\section{A cohomological proposition}\label{sec:cohom}

In the later stages of Section~\ref{sec:polys} below we will make
crucial use of a technical proposition allowing us to re-write
certain cocycles in a very explicit form. It will enable a final,
extremely concrete re-writing of the quadratic nonconventional
averages so that they are susceptible to a more direct analysis. We
prove the needed technical result in this section as
Proposition~\ref{prop:combining-indiv-cobdry-eqns}, preferring to
separate it from the main steps in the proof of
Theorem~\ref{thm:polyconv}. Surprisingly, this will rest on a
continuity result for certain measurable cohomology groups under
taking inverse limits of the base groups, which will apply after we
suitably re-cast the data we wish to simplify\footnote{This continuity result is a precursor from Moore's original papers of the more recent results of~\cite{Aus--cohomcty}.}. We will therefore
need to call on the measurable cohomology theory for compact Abelian
groups, as developed by Moore in his important sequence of
papers~\cite{Moo64(gr-cohomI-II),Moo76(gr-cohomIII),Moo76(gr-cohomIV)}.
We recall or prove those cohomological facts that we need in
Appendix~\ref{app:cohom}.

\textbf{Remark on notation}\quad We will write $\{\cdot\}:\Sone\to
[0,1)$ for the inverse to the bijection $\theta\mapsto
\rm{e}^{2\pi\rm{i}\theta}$ and $\lfr\cdot\rfr:\bbR\to \bbZ$ for the
usual `integer part' function, so these maps are related by the
equation $s-\lfr s\rfr = \{\rm{e}^{2\pi\rm{i}s}\}$, both sides of
which give the usual `fractional part' of $s\in \bbR$. \fin

\begin{prop}[Processing certain individual coboundary
equations]\label{prop:combining-indiv-cobdry-eqns} Suppose that $U$
is a compact metrizable Abelian group and $\psi:\bbZ^2\to U$ is a
homomorphism such that $\overline{\psi(\bbZ\bf{e}_1)}\cap
\overline{\psi(\bbZ\bf{e}_2)} = \{0\}$ and
$\overline{\psi(\bbZ\bf{e}_1)}\cdot \overline{\psi(\bbZ\bf{e}_2)}$
has finite index in $U$, and that $\s:\bbZ^2\times U\to\Sone$ is a
cocycle over the corresponding rotation action $R_\psi$ of $\bbZ^2$
on $(U,m_U)$. Suppose in addition that for each $i=1,2$ there are
Borel maps $b_i:U\to\Sone$ and $c_i:U\to \Sone$ so that $c_i$ is
$R_{\psi(\bf{e}_i)}$-invariant and
\[\s(\bf{e}_i,\,\cdot\,) = \Delta_{\psi(\bf{e}_i)}b_i\cdot c_i.\]

Then there are Borel maps $b_i':U\to \Sone$ such that each $c_i':=
c_i\cdot\Delta_{\psi(\bf{e}_i)}b'_i:U\to\Sone$ is a map of the form
\[c'_i(u) = \a_i(u)\cdot\exp\Big(2\pi\rm{i}\sum_{j=1}^{J_i} a_{i,j}(u)\{\chi_{i,j}(\psi(\bf{e}_i))\}\{\g_{i,j}(u)\}\Big)\]
for some function $\a_i:U\to \Sone$ that factorizes through a
\emph{finite} quotient group of $U$, functions $a_{i,j}:U\to\bbZ$
for $j=1,2,\ldots,J_i$ that also factorize through this finite
quotient group of $U$, and characters $\chi_{i,1}$, $\chi_{i,2}$,
\ldots, $\chi_{i,J_i}\in \widehat{U}$ and $\g_{i,1}$, $\g_{i,2}$,
\ldots, $\a_{i,J_i}\in \overline{\psi(\bbZ\bf{e}_i)}^\perp$.
Therefore we can write instead
\[\s(\bf{e}_i,\,\cdot\,) = \Delta_{\psi(\bf{e}_i)}(b_i\cdot b'_i)\cdot c'_i\]
with $c'_i$ a map of this special form.
\end{prop}

\textbf{Remarks}\quad\textbf{1}\quad Simply by playing around with
examples of functions $c_i$ that are already of the special form
appearing above, it is not hard to show that there are quite
nontrivial examples of $\bbZ^2$-systems admitting cocycles that
satisfy the conditions of this proposition.  For instance, let
$w_1,w_2\in \Sone$ be transcendental and algebraically independent
over $\bbQ$ when identified with classes in $\bbT:= \bbR/\bbZ$ and
also such that $0 < \{w_1\},\{w_2\} < 1/50$ and observe that if we
define $\theta' \in \Sone$ by $\{\theta'\} = \{w_1\}\{w_2\}$ then
the Borel map
\[c_1:(\Sone)^2\to \Sone:(t_2,z_2) \mapsto z_2\cdot \exp(-2\pi\rm{i}\{w_1\}\{t_2\}) = \exp(2\pi\rm{i}(\{z_2\} - \{w_1\}\{t_2\}))\]
satisfies
\[\Delta_{(w_2,\theta')}c_1(t_2,z_2) = -\lfr\{w_2\} + \{t_2\}\rfr\cdot w_1.\]

Now let $U := \Sone\times (\Sone)^2$ and $\psi:(m,n)\mapsto
(mw_1,n(w_2,\theta'))$ (this has dense image by the algebraic
independence of $\{w_1\}$, $\{w_2\}$ and $1$ over $\bbQ$), and
define
\[\s_1(\bf{e}_1,(t_1,t_2,z_2)) := c_1(t_2,z_2)\]
and
\[\s_2(\bf{e}_2,(t_1,t_2,z_2)) := -\lfr\{t_2\} + \{w_2\}\rfr\cdot t_1 = \Delta_{(0,w_2,\theta')}b(t_1,t_2,z_2)\cdot c_2(t_1)\]
where
\[b(t_1,t_2,z_2) = \exp(2\pi\rm{i}\{t_1\}\{t_2\})\quad\quad\hbox{and}\quad\quad c_2(t_1) = \exp(-2\pi\rm{i}\{t_1\}\{w_2\}).\]

We can now check immediately that
$\Delta_{w_1}\s(\bf{e}_2,\,\cdot\,) =
\Delta_{(w_2,\theta')}\s(\bf{e}_1,\,\cdot\,)$, so this does indeed
define a cocycle over the rotation action $R_\psi$ that admits
functions $b$, $c_1$ and $c_2$ as in the above proposition.
Furthermore, since $c_2(t_1)$ is $R_{(w_2,\theta')}$-invariant and
takes continuum-many different values, it cannot be an
$R_{(w_2,\theta')}$-quasi-coboundary (since for this to be true its
values would be restricted to the eigenvalue group of some rotation
on a compact metrizable Abelian group, and such an eigenvalue group
would be countable); thus, in a sense, this $c_2$ does not admit
further simplification in any obvious way, and similar remarks apply
to $c_1$.

The importance of the above proposition is that it tells us that all
such examples must be `finite-dimensional' up to cohomology, and as
the proof will show the reason behind this is very much a
cohomological one (in particular, it will rest on the continuity of
$\rmH^2(\cdot,\cdot)$ under inverse limits in the first argument,
recalled below). Although such a result seems quite surprising a
priori, we note that it does have a precedent in the study of
pro-nilsystems as characteristic factors, where it is shown that
towers of Abelian isometric extensions that are initially
characterized by the Conze-Lesigne equation and its higher-step
analogs can always be represented as inverse limits of
finite-dimensional examples (see, in
particular,~\cite{Rud93,HosKra05,Zie07}).

\quad\textbf{2}\quad It seems likely that a version of this result
is available without the simplifying assumption that
$\overline{\psi(\bbZ\bf{e}_1)}\cap \overline{\psi(\bbZ\bf{e}_2)} =
\{0\}$, but we make it here as this is the only case we will need
and this assumption does lead to a much lighter presentation. \fin

\textbf{Proof}\quad Let $w_i := \psi(\bf{e}_i)$ and $K_i :=
\overline{\psi(\bbZ\bf{e}_i)}$ for $i=1,2$. We will make use of the
cocycle condition
\[\Delta_{w_1}\s(\bf{e}_2,\,\cdot\,) = \Delta_{w_2}\s(\bf{e}_1,\,\cdot\,).\]

First, because there are only finitely many cosets of $K_1\cdot K_2$
in $U$ and these are preserved by both of the rotations $R_{w_1}$
and $R_{w_2}$, the desired conclusion clearly follows overall if we
prove it separately within each of these cosets, and so we now
simply assume that $U = K_1\cdot K_2$.  Given this, the condition
that $K_1 \cap K_2 = \{0\}$ means we may assume $U = K_1\times K_2$
and correspondingly denote points of $U$ as ordered pairs
$(u_1,u_2)$ in this product group.

Next, by adjusting the whole of $\s$ by $\Delta_\psi b_1$, we may
assume simply that $b_1 \equiv 1$. Given this, now substituting our
expressions for $\s(\bf{e}_i,\,\cdot\,)$ into the commutativity
condition we obtain
\[\Delta_{w_1}(\Delta_{w_2}b_2\cdot c_2) = \Delta_{w_2}c_1.\]

We will deduce our desired conclusion from this equation in several
small steps.

\quad\textbf{Step 1}\quad We first focus our attention on the map
$b_2$, with the goal of proving that it admits a factorization as
\[b_2(u_1,u_2) = \a(u_1,u_2)\cdot\rho_1(u_1,u_2)\cdot\rho_2(u_1,u_2)\cdot b_2'(u_1,u_2),\]
where $\a:K_1\times K_2\to \Sone$ factorizes through some finite
quotient group of $K_1\times K_2$, $\rho_1$ has the property that
that $\rho_1(u_1,\,\cdot\,)$ is a member of $\E(K_2)$ for
Haar-almost every $u_1 \in K_2$, $\rho_2$ has the symmetric property
and $b'_2$ is of the form
\[b_2'(u_1,u_2) = \exp\Big(2\pi\rm{i}\sum_{j=1}^J\{\g_j(u_1)\}\{\chi_j(u_2)\}\Big)\]
for some $\g_1$, $\g_2$, \ldots, $\g_J \in \widehat{K_1}$ and
$\chi_1$, $\chi_2$, \ldots, $\chi_J \in \widehat{K_2}$.  This will
occupy the first five steps (the bulk of the proof).

Our first step amounts to a simple re-interpretation of the various
data in hand. Observe that the right-hand side of the commutativity
equation above is $R_{w_1}$-invariant, while the left-hand side is
an $R_{w_1}$-coboundary.  This implies that
\begin{itemize}
\item $\Delta_{w_2}c_1$ takes values in $\widehat{K_1}(w_1)$, and
\item for almost every $u_2 \in K_2$ the map $\Delta_{w_2}b_2(\,\cdot\,,u_2)\cdot c_2(\,\cdot\,)$ is an
eigenfunction on the subgroup $K_1$ (noting that $c_2$ does not
depend on $u_2$ by assumption).
\end{itemize}

Thus the measurable map $\xi:u_2\mapsto
\Delta_{w_2}b_2(\,\cdot\,,u_2)\cdot c_2(\,\cdot\,)$ from $K_2$ to
the Polish Abelian group $\C(K_1)$ of isomorphism classes of Borel
maps up to almost-everywhere agreement actually almost surely takes
values in the closed subgroup $\E(K_1)$. Let us also define another
measurable map $\b:K_2\to \C(K_1)$ by $\b(u_2) :=
b_2(\,\cdot\,,u_2)$.

If we now choose any $\theta \in K_2$ and take the difference under
$\theta$ of the definition of $\xi$, then since $c_2$ is
$K_2$-invariant we obtain
\[\Delta_{w_2}\Delta_\theta \b = \Delta_\theta\xi.\]
This tells us that as members of $\C(K_1)$, $\Delta_\theta \b(u_2)$
and $\Delta_\theta \b(u_2w_2)$ almost surely differ only by a member
of $\E(K_1)$. Since $\E(K_1) \leq \C(K_1)$ is a closed subgroup and
so the quotient group carries a smooth Borel structure, and since
$R_{w_2}$ is ergodic on $K_2$, it follows that there are some fixed
Borel map $f_\theta \in \C(K_2)$ and a Borel selection of
eigenfunctions $u_2\mapsto \zeta_\theta(u_2)\in \E(K_1)$ such that
$\Delta_\theta \b(u_2) = f_\theta\cdot \zeta_\theta(u_2)$, and
moreover a simple measurable selection argument ensures that we can
take these to vary Borel measurably in $\theta$ while still
guaranteeing that this equation hold Haar-almost everywhere, so we
may write instead $\Delta_\theta \b(u_2) = f(\theta)\cdot
\zeta(\theta,u_2)$.

It follows that if we define $\overline{\b}:K_2\to \C(K_1)/\E(K_1)$
to be the quotient of $\b$ and similarly for $\overline{f}$, then
$\Delta_\theta \overline{\b}(u_2) = \overline{f}(\theta)$. Therefore
$\overline{f}$ is a homomorphism, since given $\theta$ and $\theta'$
we know that for almost every $u_2 \in K_2$ we have
\[\overline{f}(\theta)\cdot \overline{f}(\theta') = \Delta_\theta\overline{\b}(u_2)\cdot \Delta_{\theta'}\overline{\b}(u_2 \cdot \theta) = \Delta_{\theta\cdot \theta'}\overline{\b}(u_2) = \overline{f}(\theta\cdot \theta'),\]
and hence $\overline{\b}$ is an affine homomorphism (each up to
modification on a negligible set).

We may therefore find some fixed function $h \in \C(K_1)$ such that
if we write $\overline{h}$ for the image of $h$ in
$\C(K_1)/\E(K_1)$, define $\t{\b}(u_2) := \b(u_2)\cdot h$ and let
$\overline{\t{\b}}$ be its image under composition with the quotient
map $\C(K_2)\to \C(K_2)/\E(K_1)$, then this $\overline{\t{\b}}$ is a
true homomorphism.  Hence regarding it as a member of
$\Z^1(K_1,\C(K_1)/\E(K_1))$ we have $d\overline{\t{\b}} = 0$.
However, this in turn tells us that the $2$-cocycle $d\t{\b}$ takes
values in the closed subgroup $\E(K_1)$, endowed with the trivial
action of $K_2$, which we note is continuously isomorphic to
$\Sone\times \widehat{K_1}$ under the multiplication map
$(t,\chi)\mapsto t\cdot \chi$, so that $d\t{\b}$ may be identified
with a pair of $2$-cocycles, one taking values in $\bbT$ and the
other in $\widehat{K_1}$.

\quad\textbf{Step 2}\quad We now bring
Lemma~\ref{lem:more-explicit-2-cocycs} to bear on this cocycle
$d\t{\b}$. Each $K_i$ can be represented as an inverse limit of
finite-dimensional groups, say as
\[(K_i,(q_{(m),i})_{m\geq 0}) = \lim_{m\leftarrow}\big((K_{(m),i})_{m\geq 0},(q^{(m)}_{(k),i})_{m\geq k \geq 0}\big),\]
and correspondingly the group $\widehat{K_i}$ is the direct limit of
the groups $\widehat{K_{(m),i}}$ under the embeddings given by
composition with $q_{(m),i}$.  From the continuity of
$\rmH^2(\cdot,\cdot)$ given by Proposition~\ref{prop:inv-lim-cohom}
it follows that $d\t{\b}$ is cohomologous to a $2$-cocycle that
depends only on a finite-dimensional quotient group $K_{(m),2}$ of
$K_2$, and takes values in the lift of some $\widehat{K_{(m),1}}$:
that is, we can write
\[d\t{\b} = d\rho_2 \cdot \k\circ q_{(m),2}^{\times 2}\]
for some $\rho_2:K_2\to \E(K_1)$ and $2$-cocycle $\k:K_{(m),2}\times
K_{(m),2}\to \E(K_{(m),1})$.

As the dual of a finite-dimensional Abelian group,
$\widehat{K_{(m),1}}$ is finitely-generated and so the Structure
Theorem for these identifies it with some direct product
$\bbZ^D\times (\bbZ/n_1\bbZ)\times\cdots\times (\bbZ/n_r\bbZ)$.
Hence we obtain similarly $\E(K_{(m),1}) \cong \bbT\times
\bbZ^D\times (\bbZ/n_1\bbZ)\times\cdots\times (\bbZ/n_r\bbZ)$ with
trivial $K_2$-action, and so applying the relevant parts of
Lemma~\ref{lem:more-explicit-2-cocycs} to each coordinate we obtain
that, by a further adjustment of $\rho_2$ if necessary, we can
assume that $\k$ takes the form
\[\k(u_2,v_2) = \k'(u_2,v_2)\cdot \prod_{j=1}^k \g_j^{\lfr\{\chi_j(u_2)\} + \{\chi_j(v_2)\}\rfr}\]
for some $2$-cocycle $\k':K_{(m),2}\times K_{(m),2}\to \Sone\cdot
(\widehat{K_{(m),1}})_{\rm{tor}}$ that depends only on a finite
group quotient of $K_{(m),2}$ (where we write
$(\widehat{K_{(m),1}})_{\rm{tor}}$ for the torsion subgroup of
$\widehat{K_{(m),1}}$, which must in turn consist of those
characters that are lifted from the maximal finite group quotient of
$K_{(m),1}$), and finite lists $\g_1$, $\g_2$, \ldots, $\g_J \in
\widehat{K_{(m),1}}$, $\chi_1$, $\chi_2$, \ldots, $\chi_J \in
\widehat{K_{(m),2}}$.

\quad\textbf{Step 3}\quad Consider the $2$-cocycle
\[\prod_{j=1}^J \g_j^{\lfr\{\chi_j(u_2)\} + \{\chi_j(v_2)\}\rfr}\]
appearing in the above factorization.  An explicit computation shows
that this can be represented as the coboundary $d\b'$ of the
following $\C(K_1)$-valued $1$-cochain:
\[\b'(u_2)(u_1) = \prod_{j=1}^J\exp(2\pi\rm{i}\{\chi_j(u_2)\}\{\g_j(u_1)\}).\]

It follows that
\[\k'\circ q_{(m),2}^{\times 2} = \k\circ q_{(m),2}^{\times 2}\cdot \overline{d(\b'\circ q_{(m),2})} = d(\t{\b}\cdot\overline{\rho_2}\cdot \overline{\b'\circ q_{(m),2}}),\]
so the lift of $\k'$ to $K_2\times K_2$ is a $\C(K_1)$-valued
coboundary.

\quad\textbf{Step 4}\quad Let us now write $r_{(m),i}:K_{(m),i}\onto
F_{(m),i}$ for the maximal finite group quotient of $K_{(m),i}$,
whose kernel is just the identity connected component in
$K_{(m),i}$.  We have seen that $\k'$ factorizes through
$r_{(m),2}\times r_{(m),2}$ and takes values in $\E(r_{(m),1})$.

Also, from the above we have that $\k'\circ q_{(m),2}^{\times 2}$ is
a $\C(K_1)$-valued coboundary.  Since on the one hand
$r_{(m),1}\circ q_{(m),1}:K_1\onto F_{(m),1}$ has finite image, and
so its fibres all have individually positive measure, and on the
other hand our action of $K_2$ on $\C(K_1)$ is trivial, simply by
choosing a representative point from each fibre of $r_{(m),1}\circ
q_{(m),1}$ at random and sampling $\t{\b}\cdot\overline{\rho_2}\cdot
\overline{\b'\circ q_{(m),2}}$ at those points we deduce that
$\k'\circ q_{(m),2}^{\times 2}$ is actually the coboundary of some
$\C(r_{(m),1})$-valued $1$-cochain.

We will now argue further that, possibly after a finite further
increase in $m$, it must be the $\C(r_{(m),1})$-valued coboundary of
some $1$-cochain that depends only on coordinates in $F_{(m),2}$.
Indeed, this also follows directly from
Lemma~\ref{lem:more-explicit-2-cocycs}, since in view of the
triviality of the action we can simply write $\C(r_{(m),1}) \cong
\bbT^{\oplus F_{(m),1}}$ as $K_2$-modules, and for each of these
finitely many copies of $\bbT$ Part 3 of
Lemma~\ref{lem:more-explicit-2-cocycs} gives some $m'\geq m$ such
that $\k'(\,\cdot\,,\,\cdot\,)(x)$ regarded as a $\bbT$-valued
cocycle is a coboundary upon lifting only up as far as
$F_{(m'),2}\times F_{(m'),2}$.  Taking the maximum of the $m'$ so
obtained for different $x\in F_{(m),1}$ gives the result.

Hence after passing to a suitably-enlarged value of $m$ if necessary
we can express $\k' = d(\a\circ r_{(m),2})$ for some
$\a:F_{(m),2}\to \C(r_{(m),1})$, which we may of course
alternatively interpret as a $\Sone$-valued function that factorizes
through $r_{(m),1}\times r_{(m),2}$.

\quad\textbf{Step 5}\quad We have now represented the whole of $\k$
as the $\C(K_{(m),1})$-valued coboundary: $d((\a\circ
(r_{(m),2}\circ q_{(m),2})) \cdot (\b'\circ q_{(m),2})))$ where
\[\b'(u_2)(u_1) = \prod_{j=1}^k \exp(2\pi\rm{i}\{\chi_j(u_2)\}\{\g_j(u_1)\})\]
and $\a$ takes values in $\C(r_{(m),1})$.

Let us now write $\a$ and $\b'$ for the lifts of these cochains to
$K_2$ to lighten notation, omitting the compositions with
$q_{(m),2}$. Putting this factorization together with the definition
of $\k$ we have $d\t{\b} = d(\rho_2\cdot\a\cdot\b')$, and hence
$d(\t{\b}\cdot \overline{\rho_2\cdot\a\cdot\b'}) = 0$ so that
$\t{\b}\cdot \overline{\rho_2\cdot\a\cdot\b'}:K_2\to \C(K_1)$ is a
Borel homomorphism.  From this a simple inspection of the behaviour
of the map $u_2 \mapsto (\t{\b}\cdot
\overline{\rho_2\cdot\a\cdot\b'})(u_2)(u_1)$ pointwise for almost
every $u_1$ (formally, we are using Moore's treatment of
direct-integral cohomology groups in Theorem 2
of~\cite{Moo76(gr-cohomIV)}) indicates that there is some
$\rho'_1:K_1\times K_2$ such that $\rho'_1(u_1,\,\cdot\,)$ is almost
always a member of $\E(K_2)$ and
\[(\t{\b}\cdot \overline{\rho_2\cdot\a\cdot\b'})(u_2)(u_1) = \rho'_1(u_1,u_2)\] almost
everywhere.

Re-arranging this and recalling that $\t{\b}(u_1,u_2) =
b_2(u_1,u_2)h(u_1)$, we have obtained a factorization
\begin{eqnarray*}
b_2(u_1,u_2) &=& \overline{h(u_1)}\cdot
\a(u_1,u_2)\cdot\rho'_1(u_1,u_2)\cdot\rho_2(u_1,u_2)\cdot
b_2'(u_1,u_2)\\
&=& \a(u_1,u_2)\cdot\rho_1(u_1,u_2)\cdot\rho_2(u_1,u_2)\cdot
b_2'(u_1,u_2)
\end{eqnarray*}
where $\a:K_1\times K_2\to \Sone$ factorizes through the finite
quotient $r_{(m),1}\times r_{(m),2}$, $\rho_1(u_1,u_2) :=
\overline{h(u_1)}\cdot \rho_1'(u_1,u_2)$ has the property that that
$\rho_1(u_1,\,\cdot\,)$ is a member of $\E(K_2)$ for Haar-almost
every $u_1 \in K_1$ (with each value $\overline{h(u_1)}$ interpreted
simply as a constant function of $u_2$), $\rho_2$ has the symmetric
property and $b'_2$ is of the form
\[b_2'(u_1,u_2) = \exp\Big(2\pi\rm{i}\sum_{j=1}^J\{\g_j(u_1)\}\{\chi_j(u_2)\}\Big)\]
This gives us the asserted factorization of $b_2$.

\quad\textbf{Step 6}\quad Our last step is to turn the above
factorization into a suitable cohomology for each of $c_1$ and
$c_2$.

To do this we now difference the factorization of $b_2$ obtained
above with respect to $w_1$ and $w_2$ and insert the result back
into our original commutativity equation for $\s$.  This becomes
\begin{multline*}
\Delta_{w_2}c_1(u_2)\\ =
(\Delta_{w_1}\Delta_{w_2}\a(u_1,u_2))(\Delta_{w_1}\rho_1(u_1,w_2))(\Delta_{w_2}\rho_2(w_1,u_2))(\Delta_{w_1}\Delta_{w_2}b_2(u_1,u_2))\cdot\Delta_{w_1}c_2(u_1).
\end{multline*}
On the other hand, we can compute explicitly that
\begin{eqnarray*}
&&\Delta_{w_1}\Delta_{w_2}b'_2(u_1,u_2) =
\Delta_{w_1}\Big(\prod_{j=1}^J\exp(2\pi\rm{i}\{\g_j(u_1)\}(\{\chi_j(u_2
+ w_2)\} - \{\chi_j(u_2)\})\Big)\\
&&=
\Delta_{w_1}\Big(\prod_{j=1}^J\exp(2\pi\rm{i}\{\g_j(u_1)\}(\{\chi_j(w_2)\}
- \lfr\{\chi_j(u_2)\} + \{\chi_j(w_2)\}\rfr))\Big)\\
&&=
\Delta_{w_1}\Big(\prod_{j=1}^J\exp(2\pi\rm{i}\{\g_j(u_1)\}\{\chi_j(w_2)\})\cdot\prod_{j=1}^J\g_j(u_1)^{-\lfr\{\chi_j(u_2)\}
+ \{\chi_j(w_2)\}\rfr}\Big)\\
&&=
\prod_{j=1}^J\rm{e}^{2\pi\rm{i}\{\g_j(w_1)\}\{\chi_j(w_2)\}}\cdot\prod_{j=1}^J\g_j(w_1)^{-\lfr\{\chi_j(u_2)\}
+ \{\chi_j(w_2)\}\rfr}\cdot
\prod_{j=1}^J\chi_j(w_2)^{-\lfr\{\g_j(u_1)\} + \{\g_j(w_1)\}\rfr}.
\end{eqnarray*}
Also, we have
\begin{multline*}
\prod_{j=1}^J\g_j(w_1)^{-\lfr\{\chi_j(u_2)\} + \{\chi_j(w_2)\}\rfr}\\
=
\exp\Big(-2\pi\rm{i}\sum_{j=1}^J\{\g_j(w_1)\}\{\chi_j(w_2)\}\Big)\exp\Big(2\pi\rm{i}\sum_{j=1}^J\g_j(w_1)(\{\chi_j(u_2
+ w_2)\} - \{\chi_j(u_2)\})\Big)
\end{multline*}
and similarly for $\prod_{j=1}^J\chi_j(w_2)^{-\lfr\{\g_j(u_1)\} +
\{\g_j(w_1)\}\rfr}$, so we can write the above factorization as
\[\Delta_{w_1}\Delta_{w_2}b'_2(u_1,u_2) = \big(\rm{constant}\big)\cdot\Delta_{w_2}f_1(u_2)\cdot \Delta_{w_1}f_2(u_1)\]
with
\begin{multline*}
f_1(u_2) :=
\exp\Big(2\pi\rm{i}\sum_{j=1}^J\{\g_j(w_1)\}\{\chi_j(u_2)\}\Big)\\
\hbox{and}\quad\quad f_2(u_1) :=
\exp\Big(2\pi\rm{i}\sum_{j=1}^J\{\chi_j(w_2)\}\{\g_j(u_1)\}\Big).
\end{multline*}

It follows that we may re-arrange the commutativity condition to
deduce that both
\[\Delta_{w_2}(c_1(u_2)\cdot\overline{\rho_2(w_1,u_2)}\cdot\overline{f_1(u_2)})\]
and
\[\Delta_{w_1}(c_2(u_1)\cdot\rho_1(u_1,w_2)\cdot f_2(u_1))\]
must actually factorize through the finite quotient of $K_1\times
K_2$ under $r_{(m),1}\times r_{(m),2}$.

Since for any $n\geq 1$ we can form
\[\Delta_{w_2^n}(c_1(u_2)\cdot\overline{\rho_2(w_1,u_2)}\cdot\overline{f_1(u_2)})\]
by multiplying translates of
\[\Delta_{w_2}(c_1(u_2)\cdot\overline{\rho_2(w_1,u_2)}\cdot\overline{f_1(u_2)}),\]
and $R_{w_2}$ acts ergodically on $K_2$, it follows that we can find
some $n \geq 1$ such that $r_{(m),2}(w_2^n) = 1$, and thus that the
above condition tells us that in each ergodic component of
$R_{w_2^n}$ acting on $K_2$ the function
\[\Delta_{w_2^n}(c_1(u_2)\cdot\overline{\rho_2(w_1,u_2)}\cdot\overline{f_1(u_2)})\]
is constant, and hence that
\[c_1(u_2)\cdot\overline{\rho_2(w_1,u_2)}\cdot\overline{f_1(u_2)}\]
must an eigenfunction within each of these ergodic components.
Calling this function $g_1(u_2)$, and obtaining similarly
$g_2(u_1)$, one last re-arrangement gives that
\[c_1(u_2) = \rho_2(w_1,u_2)\cdot f_1(u_2)\cdot g_1(u_2) = \Delta_{w_1}\rho_2(u_1,u_2)\cdot f_1(u_2)\cdot g_1(u_2)\]
and
\[c_2(u_1) = \overline{\rho_1(w_1,u_2)}\cdot \overline{f_2(u_1)}\cdot g_2(u_1) = \Delta_{w_2}\overline{\rho_1(u_1,u_2)}\cdot \overline{f_2(u_1)}\cdot g_2(u_1).\]

Since the function $g_1(u_2)$ is an eigenfunction within each coset
of some finite-index subgroup of $K_2$, it follows that we may write
$g_1$ in the form
\[g_1(u_2) = \a_i(u_2)\prod_{j'=1}^{J'}\chi_j'(u_2)^{a_j(u_2)} = \a_i(u_2)\exp\Big(2\pi\rm{i}\sum_{j' = 1}^{J'}a_j(u_2)\{\chi'_j(u_2)\}\Big)\]
for some maps $\a_i:K_2\to \Sone$ and $a_j:K_2\to\bbZ$ that
factorize through some finite quotient group of $K_2$, and some
additional characters $\chi'_j \in \widehat{K_2}$. Combining this
with the explicit form obtained above for $f_1(u_2)$ and noting that
$\Delta_{w_1}\rho_2(u_1,u_2)$ is an $R_{w_1}$-coboundary, we see
that we have put $c_1(u_2)$ explicitly into the desired form, and
similarly for $c_2(u_1)$. This completes the proof of
Proposition~\ref{prop:combining-indiv-cobdry-eqns}. \qed

\section{Proof of the main theorem}\label{sec:polys}

We now turn to Theorem~\ref{thm:polyconv}:

\begin{thm*}
If $T_1,T_2:\bbZ\curvearrowright (X,\mu)$ commute then the averages
\[S_N(f_1,f_2) := \frac{1}{N}\sum_{n=1}^N (f_1\circ T_1^{n^2})(f_2\circ T_1^{n^2}T_2^n)\]
converge in $L^2(\mu)$ as $N\to\infty$ for any $f_1,f_2 \in
L^\infty(\mu)$.
\end{thm*}

The proof proceeds through a sequence of three reductions to
progressively simpler classes of polynomial average, each obtained
by deriving different consequences from some invocation of the van
der Corput estimate.  After the third reduction we will reach a
family of averages to which known results can be applied
more-or-less directly.

In rough outline, our first reduction amounts to an identification
of characteristic factors for these polynomial averages in some
pleasant extension, so that we may assume the functions $f_1$ and
$f_2$ take a special form in terms of these factors. This use of
characteristic factors is another outing for what is now the
standard approach to such questions. It is for this first step that
we will need the result for linear averages of
Theorem~\ref{thm:char-three-lines-in-2D}.  In fact, we will need
just a little more versatility than is contained in
Theorem~\ref{thm:char-three-lines-in-2D} as stated, but which
follows at once from combining that theorem with the following
immediate consequence of the definition of a characteristic tuple of
factors (see Lemma 4.3 in~\cite{Aus--lindeppleasant1}):

\begin{lem}\label{lem:char-factors-symmetry}
For any factor $\xi:\bfX\to \bfY$ the triple $(\xi,\id_X,\id_X)$ is
characteristic for the nonconventional averages
\[\frac{1}{N}\sum_{n=1}^N(f_1\circ T^{n\bf{p}_1})(f_2\circ T^{n\bf{p}_2})(f_3\circ T^{n\bf{p}_3}),\quad\quad f_1,f_2,f_3\in L^\infty(\mu),\]
if and only if the triple $(\id_X,\xi,\id_X)$ is characteristic for
the nonconventional averages
\[\frac{1}{N}\sum_{n=1}^N(f_0\circ T^{-n\bf{p}_j})(f_1\circ T^{n(\bf{p}_1 -\bf{p}_j)})(f_k\circ T^{n(\bf{p}_k - \bf{p}_j)}),\quad\quad f_0,f_1,f_k \in L^\infty(\mu),\]
whenever $\{j,k\} = \{2,3\}$. \qed
\end{lem}

\begin{cor}\label{lem:symmetry-in-the-structthm}
In the statement of Theorem~\ref{thm:char-three-lines-in-2D} we may
instead let the target system of $\eta$ be a $(\bf{p}_1 -
\bf{p}_2,m_{13}(\bf{p}_1 - \bf{p}_3),m_1\bf{p}_1)$-directional
CL-system or a $(\bf{p}_1 - \bf{p}_3,m_{12}(\bf{p}_1 -
\bf{p}_2),m_1\bf{p}_1)$-directional CL-system. \qed
\end{cor}

Our use for Theorem~\ref{thm:char-three-lines-in-2D} and
Lemma~\ref{lem:symmetry-in-the-structthm} will be to prove an even
more precise description of a characteristic pair of factors for our
nonconventional quadratic averages, by considering a whole family of
triple linear averages that arise from those quadratic averages
through an appeal to the van der Corput estimate, and then examining
the possible joint distribution of the characteristic factors for
those different triple linear averages inside the overall system.
The result of this step will be Theorem~\ref{thm:char-poly}.

The second reduction then follows quite quickly and uses similar
ideas: after simplifying the averages $S_N$ for functions measurable
with respect to the new characteristic pair of factors and
re-arranging slightly, a new sequence of averages emerges to which
another appeal to Theorem~\ref{thm:char-three-lines-in-2D} and the
resulting description of the Furstenberg self-joining gives a
further simplification.

The proof is completed through a closer examination of some
functions measurable with respect to a
$(h\bbZ^2,(h,0),(0,h))$-directional CL-system for some $h$. This is
heavily based on an older approach of Host and Kra~\cite{HosKra01}
to the study of the triple linear nonconventional averages
associated to three powers of a single transformation that does not
need the exact picture in terms of nilsystems, which was not
available at the time of that paper. It amounts to a way of using
directly the combined cocycle equation arising from the Mackey data
inside the Furstenberg self-joining of our system. This leads to a
classification of the polynomial averages output by the second
reduction into two cases. In the first case we can show they tend to
$0$ in $L^2(\mu)$, and in the second we will find that they can
eventually be re-written simply as a more classical sequence of
weighted ergodic averages, for which mean convergence is known.

\subsection{Directional CL-systems}\label{subs:dCL-again}

We are now ready to introduce the `directional CL-systems' that are
the main new ingredient that appear in
Theorem~\ref{thm:char-three-lines-in-2D}.  In this subsection we
will define these systems and establish some of their basic
properties.

Directional CL-cocycles are characterized by the existence of
solutions to some natural `directional' analogs of the classic
Conze-Lesigne equations among cocycles~(\cite{ConLes84,Les93}).  Let
us first introduce these equations, and then the class of cocycles
that they specify.

\begin{dfn}[Directional Conze-Lesigne equations]
Suppose that $A$ and $Z$ are compact metrizable Abelian groups, $K
\leq Z$ a closed subgroup and $\tau:Z\to A$ a Borel map. Then
another Borel map $b:Z\to A$ \textbf{satisfies the directional
Conze-Lesigne equation E$(u,v,K,\tau)$} for some $u,v\in Z$ if there
is a Borel map $c:Z/K\to A$ such that
\[\Delta_u\tau(z) = \Delta_vb(z)\cdot c(z\cdot K)\quad\quad \hbox{for }m_Z\hbox{-a.e. }z.\]
It is clear that this $c$ is then uniquely determined. We refer to
$b$ as a \textbf{solution} of the equation E$(u,v,K,\tau)$ and to
$c$ as the \textbf{one-dimensional auxiliary} of $b$ in this
equation. This is the classical Conze-Lesigne equation in case $K =
G$.
\end{dfn}

Although we have formulated the above definition for cocycles into
an arbitrary compact Abelian target group $A$, for technical reasons
we will use this equation only for cocycles into $\Sone$.

\textbf{Remark on notation}\quad Similarly to~\cite{Aus--lindeppleasant2}, we will henceforth write
$(Z_\star,m_{Z_\star},\phi_\star)$ to denote a $\bbZ^2$-system whose
underlying space is the direct integral of some measurably-varying
family of compact Abelian groups $Z_\star$, indexed by some other
standard Borel probability space $(S,\nu)$ on which the action is
trivial, with the overall action a fibrewise rotation defined by a
measurable selection for each fibre $Z_s$ of a dense homomorphism
$\phi_s:\bbZ^2\to Z_s$: writing $R_\phi$ for this action, it is
given by
\[R_\phi^{\bf{n}}(s,z) := (s,z\cdot \phi_s(\bf{n}))\quad\quad\hbox{for\ }s\in S,\ z\in Z_s\ \hbox{and}\ \bf{n} \in \bbZ^2.\]
We will refer to such a system as a \textbf{direct integral of
ergodic group rotations} and to $(S,\nu)$ as its \textbf{invariant
base space}. Sometimes we omit the base space $(S,\nu)$ from mention
completely, since once again the forthcoming arguments will all
effectively be made fibrewise, just taking care that all
newly-constructed objects can still be selected measurably. In
particular, we will often write just $Z_\star$ in place of $S\ltimes
Z_\star$. \fin

\begin{dfn}[Directional CL-cocycles]\label{dfn:dCL-cocyc}
Suppose that $\bf{n}_1$, $\bf{n}_2$, $\bf{n}_3 \in \bbZ^2$, that
$(Z_\star,m_{Z_\star},\phi_\star)$ is a direct integral of ergodic
$\bbZ^2$-group rotations with invariant base space $(S,\nu)$, and
that $A_\star$ is motionless compact metrizable Abelian group data
over $(Z_\star,m_{Z_\star},\phi_\star)$.

A cocycle-section $\tau:\bbZ^2\times Z_\star\to A_\star$ over the
fibrewise rotation action $R_\phi$ is an
\textbf{$(\bf{n}_1,\bf{n}_2,\bf{n}_3)$-directional CL-cocycle over
$R_\phi$} if for every $R_\phi$-invariant measurable selection of
characters $\chi_\star\in \widehat{A_\star}$ we have that
\begin{itemize}
\item for every $R_\phi$-invariant measurable selection
$u_\star\in\overline{\phi_\star(\bbZ\bf{n}_2)}$ there is a Borel map
$b:S\ltimes Z_\star \to \Sone$, denoted by $b_\star$, such that
$b_s$ solves the equation
E$(u_s,\phi_s(\bf{n}_1),\overline{\phi_s(\bbZ\bf{n}_3)},\chi_s\circ\tau(\bf{n}_1,\,\cdot\,)|_{Z_s})$
for $\nu$-almost every $s$, and
\item for every $R_\phi$-invariant measurable selection
$v_\star\in\overline{\phi_\star(\bbZ\bf{n}_3)}$ there is a Borel map
$b_\star:S\ltimes Z_\star \to \Sone$ that solves the equation
E$(v_s,\phi_s(\bf{n}_1),\overline{\phi_s(\bbZ\bf{n}_2)},\chi_s\circ\tau(\bf{n}_1,\,\cdot\,)|_{Z_s})$
for $\nu$-almost every $s$.
\end{itemize}

Given a subgroup $\G \leq \bbZ^2$, $\tau$ is a
\textbf{$(\G,\bf{n}_2,\bf{n}_3)$-directional CL-cocycle over
$R_\phi$} if for every $R_\phi$-invariant measurable selection of
characters $\chi_\star\in \widehat{A_\star}$ we have that
\begin{itemize}
\item for every $R_\phi$-invariant measurable selection
$u_\star\in\overline{\phi_\star(\bbZ\bf{n}_2)}$ there is a Borel map
$b_\star:S\ltimes Z_\star \to \Sone$ that \emph{simultaneously}
solves the equations
E$(u_s,\phi_s(\bf{n}_1),\overline{\phi_s(\bbZ\bf{n}_3)},\chi_s\circ\tau(\bf{n}_1,\,\cdot\,)|_{Z_s})$,
$\bf{n}_1\in\G$, for $\nu$-almost every $s$, and
\item for every $R_\phi$-invariant measurable selection
$v_\star\in\overline{\phi_\star(\bbZ\bf{n}_3)}$ there is a Borel map
$b_\star:S\ltimes Z_\star \to \Sone$ that \emph{simultaneously}
solves the equations
E$(v_s,\phi_s(\bf{n}_1),\overline{\phi_s(\bbZ\bf{n}_2)},\chi_s\circ\tau(\bf{n}_1,\,\cdot\,)|_{Z_s})$,
$\bf{n}_1 \in \G$, for $\nu$-almost every $s$.
\end{itemize}
\end{dfn}

In the above situation we will usually write more briefly that
\begin{quote}
`for every $\chi_\star \in \widehat{A_\star}$ and $u_\star \in
\overline{\phi_\star(\bbZ\bf{n}_2)}$, the map $b_\star:Z_\star\to
\Sone$ is a solution to the equations
E$(u_\star,\phi_\star(\bf{n}_1),\overline{\phi_\star(\bbZ\bf{n}_3)},\chi_\star\circ\tau(\bf{n}_1,\,\cdot\,))$',
\end{quote}
and similarly for the other equations (note, in particular, that the
restriction of $\tau(\bf{n}_1,\,\cdot\,)$ to the relevant fibre
$Z_\star$ is left to the understanding).

\begin{lem}
If $\G\leq \bbZ^2$ is a subgroup generated by a subset $F \subset
\bbZ^2$ then a cocycle-section $\tau:\bbZ^2\times Z_\star \to
A_\star$ is a $(\G,\bf{n}_2,\bf{n}_3)$-directional CL-cocycle over
$R_\phi$ for every $\bf{n}_1\in F$ if the simultaneous solutions
required above exist only for all of the families of equations
\[\bigvee_{\bf{n}_1 \in F}\rm{E}(u_\star,\phi_\star(\bf{n}_1),\overline{\phi_\star(\bbZ\bf{n}_3)},\chi_\star\circ\tau(\bf{n}_1,\,\cdot\,))\]
and
\[\bigvee_{\bf{n}_1 \in F}\rm{E}(v_\star,\phi_\star(\bf{n}_1),\overline{\phi_\star(\bbZ\bf{n}_2)},\chi_\star\circ\tau(\bf{n}_1,\,\cdot\,)).\]
\end{lem}

\textbf{Proof}\quad This follows from the simple property of the
directional Conze-Lesigne equations that if, say, $u \in
\overline{\phi_s(\bbZ\bf{n}_2)}$, $\bf{n},\bf{n}' \in F$ and $b$
solves the equations
\[\rm{E}(u,\phi_s(\bf{n}_1),\overline{\phi_s(\bbZ\bf{n}_3)},\chi_s\circ\tau(\bf{n}_1,\,\cdot\,)|_{Z_s})\]
for both $\bf{n}_1 = \bf{n}$ and $\bf{n}'$ with respective
one-dimensional auxiliaries $c$ and $c'$, then
\begin{eqnarray*}
\Delta_u \tau(\bf{n} + \bf{n}',z) &=& \Delta_u\tau(\bf{n},z +
\phi_s(\bf{n}'))\cdot\Delta_u\tau(\bf{n}',z)\\ &=&
\Delta_{\bf{n}}b(z + \phi_s(\bf{n}'))\cdot \Delta_{\bf{n}'}b(z)\\
&&\quad\quad\quad\quad\cdot c((z + \phi_s(\bf{n}'))\cdot
\overline{\phi_s(\bbZ\bf{n}_3)})\cdot c'(z\cdot
\overline{\phi_s(\bbZ\bf{n}_3)})\\
&=& \Delta_{\bf{n} + \bf{n}'}b(z)\cdot c''(z\cdot
\overline{\phi_s(\bbZ\bf{n}_3)})
\end{eqnarray*}
at $m_{Z_s}$-a.e. $z$, where $c''$ is the obvious product function
formed from $c$ and $c'$.  Therefore $b$ is also a solution to
\[\rm{E}(u,\phi_s(\bf{n} + \bf{n}'),\overline{\phi_s(\bbZ\bf{n}_3)},\chi_s\circ\tau(\bf{n} + \bf{n}',\,\cdot\,)|_{Z_s}).\]
A similar argument shows that it also solves
\[\rm{E}(u,\phi_s(-\bf{n}),\overline{\phi_s(\bbZ\bf{n}_3)},\chi_s\circ\tau(-\bf{n},\,\cdot\,)|_{Z_s}),\]
and so in fact it applies to the whole subgroup $\G$, as required.
\qed

\textbf{Remark}\quad For the above proof it would clearly not be
enough to demand that the equations
E$(u_s,\phi_s(\bf{n}_1),\overline{\phi_s(\bbZ\bf{n}_3)},\chi_s\circ\tau(\bf{n}_1,\,\cdot\,)|_{Z_s})$
for different $\bf{n}_1 \in \G$ have solutions separately.  The
requirement of simultaneous solutions when working with
$(\G,\bf{n}_2,\bf{n}_3)$-directional CL-cocycles will be very
important later precisely so that we can use similar manipulations
again. \fin

With the above preparations behind us, we can now define our new
class of systems itself.

\begin{dfn}[Directional CL-extensions and systems]\label{dfn:dCL}
If $\bfX$ is a $\bbZ^2$-system, $(Z_\star,m_{Z_\star},\phi_\star)$
is a direct integral of ergodic $\bbZ^2$-group rotations and
$\pi:\bfX\to (Z_\star,m_{Z_\star},\phi_\star)$ is a factor map, then
$\bfX$ is an \textbf{$(\bf{n}_1,\bf{n}_2,\bf{n}_3)$-directional
CL-extension} of $(Z_\star,m_{Z_\star},\phi_\star)$ through $\pi$ if
it can be coordinatized as $(Z_\star,m_{Z_\star},\phi_\star)\ltimes
(A_\star,m_{A_\star},\tau)$ with $\pi$ the canonical factor and
$\tau$ an $(\bf{n}_1,\bf{n}_2,\bf{n}_3)$-directional CL-cocycle over
$R_\phi$. More loosely, $\bfX$ is an
\textbf{$(\bf{n}_1,\bf{n}_2,\bf{n}_3)$-directional CL-system} if it
is an $(\bf{n}_1,\bf{n}_2,\bf{n}_3)$-directional CL-extension of
some factor that is a direct integral of group rotations, and then
any suitable choice for this group-rotation factor is a
\textbf{base} for $\bfX$.

If $\G\leq \bbZ^2$ then $\bfX$ is a
\textbf{$(\G,\bf{n}_2,\bf{n}_3)$-directional CL-extension} of
$(Z_\star,m_{Z_\star},\phi_\star)$ if the above coordinatization is
possible with $\tau$ a $(\G,\bf{n}_2,\bf{n}_3)$-directional
CL-cocycle.

We will write $\sfZ_\dCL^{\G,\bf{n}_2,\bf{n}_3}$ for the class of
$(\G,\bf{n}_2,\bf{n}_3)$-directional CL-systems, and generally write
this as $\sfZ_\dCL^{\bf{n}_1,\bf{n}_2,\bf{n}_3}$ if $\G =
\bbZ\bf{n}_1$.
\end{dfn}

The elementary properties of directional CL-cocycles follow easily
from the directional Conze-Lesigne equations.

\begin{lem}\label{lem:dCL-mult-and-lift}
Suppose that $\pi:(\t{Z}_\star,m_{\t{Z}_\star},\t{\phi}_\star)\to
(Z_\star,m_{Z_\star},\phi_\star)$ is a tower of direct integrals of
$\bbZ^2$-group rotations. Then
\begin{enumerate}
\item[(1)] if $\tau_1:\bbZ^2\times Z_\star\to A_\star$ is a
$(\G,\bf{n}_2,\bf{n}_3)$-directional CL-cocycle over $R_\phi$ then
$\tau_1\circ\pi$ is a $(\G,\bf{n}_2,\bf{n}_3)$-directional
CL-cocycle over $R_{\t{\phi}}$;
\item[(2)] if $\tau_2:\bbZ^2\times Z_\star\to A_\star$ is another $(\G,\bf{n}_2,\bf{n}_3)$-directional CL-cocycle over
$R_\phi$ then $\tau_1\cdot \tau_2$ is also a
$(\G,\bf{n}_2,\bf{n}_3)$-directional CL-cocycle over $R_{\phi}$;
\item[(3)] $(A_{(m),\star})_{m\geq 1}$, $(\Phi^{(m)}_{(k),\star})_{m\geq k\geq
0}$ is a motionless measurable family of inverse sequences of
compact Abelian groups over $(Z_\star,m_{Z_\star},\phi_\star)$ with
inverse limit family $A_{(\infty),\star}$,
$(\Phi_{(m),\star})_{m\geq 0}$ (which is clearly still measurable),
and $\tau_{(m)}:\bbZ^2\times Z_\star\to A_{(m),\star}$ is a family
of $(\G,\bf{n}_2,\bf{n}_3)$-directional CL-cocycles over $R_\phi$
satisfying the consistency equations $\tau_{(k)} =
\Phi^{(m)}_{(k),\star}\circ\tau_{(m)}$ for $m\geq k\geq 0$, then the
resulting inverse limit cocycle $\tau_{(\infty)}:\bbZ^2\times
Z_\star\to A_{(\infty),\star}$ is also a
$(\G,\bf{n}_2,\bf{n}_3)$-directional CL-cocycle.
\end{enumerate}
\end{lem}

\textbf{Proof}\quad  The first two parts follow immediately from
lifting and multiplying solutions to the directional Conze-Lesigne
equations, since $\pi$ must map each group rotation fibre of
$(\t{Z}_\star,m_{\t{Z}_\star},\t{\phi}_\star)$ onto a group rotation
fibre of $(Z_\star,m_{Z_\star},\phi_\star)$ via a measurably-varying
continuous affine epimorphism.

For the third part, first recall that by construction any character
on an inverse limit of compact Abelian groups factorizes through
some finite level of the inverse sequence.  This implies that for
any measurable selection of characters $\chi_\star\in
\widehat{A_{(\infty),\star}}$ we can find a measurable selection of
positive integers $m_\star$ such that $\chi_\star$ factorizes
through $\Phi_{(m_\star),\star}:A_{(\infty),\star}\to
A_{(m_\star),\star}$ almost surely (so
$\chi_\star\circ\tau_{(\infty)} = \chi_\star'\circ\tau_{(m_\star)}$
for some measurable selection of characters satisfying $\chi_\star=
\chi_\star'\circ\Phi_{(m_\star),\star}$). Now we may simply call on
the solutions to the directional Conze-Lesigne equations for this
$\tau_{(m_\star)}$ within each level set of the map $m_\star$, to
see that these patch together to give solutions to the directional
Conze-Lesigne equations for $\tau_{(\infty)}$. Note that this last
step illustrates the usefulness of defining directional CL-cocycles
in terms of the behaviour of their compositions with characters,
rather than directly, as discussed above. \qed

Now suppose that $(Z_{i,\star},m_{Z_{i,\star}},\phi_{i,\star})$ are
direct integrals of ergodic $\bbZ^2$-group rotations for $i=1,2$ and
that $\theta$ is a joining of them.  Then we may form the
measurably-varying family of compact Abelian groups
$Z_{1,\star}\times Z_{2,\star}$ simply by taking the product of the
underlying invariant base spaces $(S_i,\nu_i)$, and then taking the
products of the two fibres of each pair of index points $(s_1,s_2)$
from those spaces; and similarly we can define the obvious
homomorphism $(\phi_{1,s_1},\phi_{2,s_2}):\bbZ^2\to Z_{1,s_1}\times
Z_{2,s_2}$ above each such pair of index points. Now a simple
application of the non-ergodic Mackey
Theorem (Theorem 2.1 in~\cite{Aus--lindeppleasant2}) shows that $\theta$ decomposes
further into a direct integral of Haar measures on the cosets of the
measurably-varying family of subgroups
\[\overline{\{(\phi_{1,s_1}(\bf{n}),\phi_{2,s_2}(\bf{n})):\ \bf{n}\in\bbZ^2\}} \leq Z_{1,s_1}\times Z_{2,s_2},\]
and so the joined system $(Z_{1,\star}\times
Z_{2,\star},\theta,(\phi_{1,\star},\phi_{2,\star}))$ can also be
expressed as a direct integral of ergodic $\bbZ^2$-group rotations
(although the ergodic fibres may be strictly smaller than
$Z_{1,\star}\times Z_{2,\star}$, and the underlying invariant index
space correspondingly larger).

Combined with the above lemma this implies that given two
$(\G,\bf{n}_2,\bf{n}_3)$-directional CL-extensions $\pi_i:\bfX_i\to
(Z_{i,\star},m_{Z_{i,\star}},\phi_{i,\star})$ and any joining
$\theta$ as above, the lift of $\theta$ to a relatively independent
joining $\l$ of $\bfX_1$ and $\bfX_2$ gives a joint system that is a
$(\G,\bf{n}_2,\bf{n}_3)$-directional CL-extension of
$(Z_{1,\star}\times
Z_{2,\star},\theta,(\phi_{1,\star},\phi_{2,\star}))$. This will be
an important observation for us when combined with the following
proposition.

\begin{prop}\label{prop:dCL-change-of-base}
Suppose that $\pi:\bfX = (X,\mu,T)\to
(Z_\star,m_{Z_\star},\phi_\star)$ is a
$(\G,\bf{n}_2,\bf{n}_3)$-directional CL-extension, and that
$(\t{Z}_\star,m_{\t{Z}_\star},\t{\phi}_\star)$ is another direct
integral of ergodic $\bbZ^2$-group rotations which can be located
into a tower of systems
\[\bfX \stackrel{\t{\pi}}{\longrightarrow}
(\t{Z}_\star,m_{\t{Z}_\star},\t{\phi}_\star)\stackrel{\a}{\longrightarrow}
(Z_\star,m_{Z_\star},\phi_\star)\] so that $\t{\pi}$ is a relatively
ergodic extension. Then $\bfX$ is also a
$(\G,\bf{n}_2,\bf{n}_3)$-directional CL-extension of
$(\t{Z}_\star,m_{\t{Z}_\star},\t{\phi}_\star)$.
\end{prop}

\textbf{Proof}\quad This breaks into two steps.

\quad\textbf{Step 1}\quad We first show that the result holds when
$\t{\pi} = \pi\vee\zeta_0^T$ (so
$(\t{Z}_\star,m_{\t{Z}_\star},\t{\phi}_\star)$ is simply a
coordinatization of the factor of $\bfX$ generated by the base copy
of $(Z_\star,m_{Z_\star},\phi_\star)$ and the overall invariant
factor --- this is easily seen to be another direct integral of
ergodic group rotations, with the same fibres as
$(Z_\star,m_{Z_\star},\phi_\star)$ but possibly an enlargement of
the invariant base system). This is the smallest possible choice
that gives $\t{\pi}$ relatively ergodic. Let $(S,\nu)$ be the
invariant base space underlying $(Z_\star,m_{Z_\star},\phi_\star)$.

Suppose that $\tau:\bbZ^2\times Z_\star\to A_\star$ is the
$(\G,\bf{n}_2,\bf{n}_3)$-directional CL-cocycle over $R_\phi$
corresponding to some coordinatization of $\pi$.  In this case the
non-ergodic Mackey Theorem gives a
precise coordinatization of $\t{\pi}$: there are a motionless family
of closed subgroups $K_\star \leq A_\star$ and a measurable section
$\rho:S\ltimes Z_\star\to A_\star$ such that $\t{\pi}$ can be
coordinatized by the factor map
\[(S\ltimes Z_\star)\ltimes A_\star \to S\ltimes (A_\star/K_\star):((s,z),a)\mapsto (s,a \cdot \rho(s,z) \cdot K_{(s,z)}),\]
and so $\pi\vee\zeta_0^T$ in turn is coordinatized by
\[(S\ltimes Z_\star)\ltimes A_\star \to (S\ltimes Z_\star)\ltimes (A_\star/K_\star):((s,z),a)\mapsto ((s,z),a \cdot \rho(s,z) \cdot K_{(s,z)}).\]

If we now simply re-coordinatize $\pi$ by fibrewise rotations by
$\rho$, then $\tau$ is replaced by $\tau' := \tau \cdot
\Delta_{\phi}\rho$ so this now almost surely takes values in
$K_\star$, and this leads to an explicit recoordinatization of the
extension $\pi\vee\zeta_0^T$ as
\begin{center}
$\phantom{i}$\xymatrix{ \bfX
\ar[dr]_-{\pi\vee\zeta_0^T}\ar@{<->}[r]^-\cong & (Z_\star\ltimes
(A_\star/K_\star),m_{Z_\star\ltimes
(A_\star/K_\star)},(\phi_\star,1_{A_\star/K_\star}))\ltimes
(K_\star,m_{K_\star},\tau')\ar[d]^{\rm{canonical}}\\
& (Z_\star\ltimes (A_\star/K_\star),m_{Z\ltimes
(A_\star/K_\star)},(\phi_\star,1_{A_\star/K_\star})) }
\end{center}
(where we again abbreviate $S\ltimes Z_\star$ to $Z_\star$). In this
diagram the base system $(Z_\star\ltimes
(A_\star/K_\star),m_{Z_\star\ltimes
(A_\star/K_\star)},(\phi_\star,1_{A_\star/K_\star}))$ is expressed
as a direct integral of not-necessarily ergodic group rotations ---
indeed, the homomorphisms $\bf{n}\mapsto
(\phi_s(\bf{n}),1_{A_s/K_s})$ cannot have dense image unless $K_s =
A_s$ --- but by cutting down the fibres and enlarging the invariant
base system as previously it may clearly be re-coordinatized as a
direct integral of ergodic group rotations with the same fibres
$Z_\star$ as originally.

Since $\tau'$ depends only on the factor $Z_\star\ltimes
(A_\star/K_\star)\to Z_\star$ (since this is true of $\tau$ and
$\rho$), it suffices to show that $\tau'$, like $\tau$, admits
solutions to all the relevant directional Conze-Lesigne equations.
If $\chi_\star \in \widehat{K_\star}$ is a measurable selection of
characters then we can extend each $\chi_s$ to a character on the
whole of $A_s$ which we also denote by $\chi_s$ (it is classical
that this is always possible; see, for instance, Theorem 24.12 of
Hewitt and Ross~\cite{HewRos79}), and a simple appeal to the
Measurable Selector Theorem promises that we can choose these
extensions so as still to form a measurable family. Now if $\bf{n}
\in \G$, $u \in \overline{\phi_s(\bbZ\bf{n}_2)}$ for some $s$ and
$b$ is a solution to the equation
E$(u,\phi_s(\bf{n}_1),\overline{\phi_s(\bbZ\bf{n}_3)},\chi_s\circ\tau(\bf{n},\,\cdot\,)|_{Z_s})$
with one-dimensional auxiliary $c$, then we check at once that $b'
:= b\cdot \Delta_u (\chi_s\circ\rho|_{Z_s})$ satisfies
\begin{eqnarray*}
&&\Delta_{\phi_s(\bf{n})}b'(z)\cdot c(z\cdot
\overline{\phi_s(\bbZ\bf{n}_3)})\\ &&=
\Delta_{\phi_s(\bf{n})}\Delta_u(\chi_s\circ\rho|_{Z_s})\cdot\Delta_{\phi_s(\bf{n})}b(z)\cdot
c(z\cdot \overline{\phi_s(\bbZ\bf{n}_3)})\\
&&=
\Delta_{\phi_s(\bf{n})}\Delta_u(\chi_s\circ\rho|_{Z_s})\cdot\Delta_u(\chi_s\circ\tau(\bf{n},\,\cdot\,)|_{Z_s})\\
&&= \Delta_u(\chi_s\circ\tau'(\bf{n},\,\cdot\,)|_{Z_s}).
\end{eqnarray*}
Performing this procedure fibrewise on the Borel map $b_\star$ that
gives a solution for a measurable selection $u_\star$ clearly gives
a new Borel map $b'_\star$ as the new solution, as required.

\quad\textbf{Step 2}\quad We now prove the general case. In fact
this makes very little appeal to the exact structure of the system
$(\t{Z}_\star,m_{\t{Z}_\star},\t{\phi}_\star)$.

By Step 1 we can replace $\pi:\bfX\to
(Z_\star,m_{Z_\star},\phi_\star)$ by a suitable coordinatization of
$\pi\vee \zeta_0^T$ if necessary, and so suppose that $\pi$ itself
is relatively ergodic. Suppose again that $\tau:\bbZ^2\times
Z_\star\to A_\star$ is the $(\G,\bf{n}_2,\bf{n}_3)$-directional
CL-cocycle over $R_\phi$ of a coordinatization of $\pi$.  Clearly
$\a: (\t{Z}_\star,m_{\t{Z}_\star},\t{\phi}_\star)\to
(Z_\star,m_{Z_\star},\phi_\star)$ is also a relatively ergodic
Abelian isometric extension, so these two direct integrals of
ergodic group rotations have the same underlying invariant base
space, and since now both $\pi$ and $\a$ are relatively ergodic the
Relative Factor Structure Theorem (Theorem 2.5 in~\cite{Aus--lindeppleasant2}) applied to the
triangle
\begin{center}
$\phantom{i}$\xymatrix{ \bfX \ar[dr]_{\pi}\ar[r]^-{\t{\pi}} &
(\t{Z}_\star,m_{\t{Z}_\star},\t{\phi}_\star)\ar[d]^{\a}\\ &
(Z_\star,m_{Z_\star},\phi_\star)}
\end{center}
gives that there is some $R_\phi$-invariant family of quotients of
Abelian groups $q_\star:A_\star\to A_{0,\star}$ such that
\begin{center}
$\phantom{i}$\xymatrix{\bfX\ar[d]_{\t{\pi}}\ar@{<->}[rr]^-\cong
& & (Z_\star,m_{Z_\star},\phi_\star)\ltimes (A_\star,m_{A_\star},\tau)\ar[d]^{\id_{Z_\star}\ltimes q_\star}\\
(\t{Z}_\star,m_{\t{Z}_\star},\t{\phi}_\star)\ar[dr]_{\a}\ar@{<->}[rr]^-\cong
& & (Z_\star,m_{Z_\star},\phi_\star)\ltimes
(A_{0,\star},m_{A_{0,\star}},q_\star\circ\tau)\ar[dl]^{\rm{canonical}}\\
& (Z_\star,m_{Z_\star},\phi_\star). }
\end{center}

Choosing a $R_\phi$-invariant measurable selector
$\eta_\star:A_{0,\star}\to A_\star$, we can now give an explicit
re-coordinatization of the extension $\t{\pi}:\bfX\to
(\t{Z}_\star,m_{\t{Z}_\star},\t{\phi}_\star)$ as
\begin{center}
$\phantom{i}$\xymatrix{ \bfX\ar[d]_{\t{\pi}}\ar@{<->}[rr]^-\cong & &
(Z_\star\ltimes A_{0,\star},m_{Z_\star\ltimes
A_{0,\star}},(\phi_\star\ltimes \l_\star))\ltimes (\ker
q_\star,m_{\ker
q_\star},\t{\tau})\ar[d]^{\rm{canonical}}\\
(\t{Z}_\star,m_{\t{Z}_\star},\t{\phi}_\star) \ar@{<->}[rr]^-\cong &
& (Z_\star\ltimes A_{0,\star},m_{Z_\star\ltimes
A_{0,\star}},(\phi_\star\ltimes \l_\star)) }
\end{center}
for a suitable measurable selection of dense homomorphisms
$\l_\star:\bbZ^2\longrightarrow A_{0,\star}$, where the top
isomorphism is obtained by composing the previous coordinatization
$\bfX\cong (Z_\star,m_{Z_\star},\phi_\star)\ltimes
(A_\star,m_{A_\star},\tau)$ with the map
\[((s,z),a) \mapsto ((s,z),\ q_s(a),\ a\cdot \eta_s(q_s(a))^{-1}).\]
This results in a cocycle
\[\t{\tau}(\bf{n},(s,z,a_0)) := \tau(\bf{n},(s,z)) \cdot
\big(\eta_s(a_0\cdot q_s(\tau(\bf{n},(s,z))))\cdot
\eta_s(a_0)^{-1}\big)^{-1}\in \ker q_s\] for $(s,z,a_0)\in S\ltimes
Z_\star\ltimes A_{0,\star}$.

As in Step 1, it remains simply to verify that for any
measurably-varying $\chi_\star \in \widehat{\ker q_\star}$ the
cocycle $\t{\tau}:Z_\star\ltimes A_{0,\star} \to \ker q_\star$
admits $\Sone$-valued solutions to the equations
\[\rm{E}(u_\star,\phi_\star(\bf{n}_1),\overline{\phi_\star(\bbZ\bf{n}_3)},\chi_\star\circ\t{\tau}(\bf{n},\,\cdot\,))\]
for every $\bf{n} \in \G$ and $u_\star \in
\overline{\phi_\star(\bbZ\bf{n}_2)}$, and
\[\rm{E}(v_\star,\phi_\star(\bf{n}_1),\overline{\phi_\star(\bbZ\bf{n}_2)},\chi_\star\circ\t{\tau}(\bf{n},\,\cdot\,))\]
for every $\bf{n} \in \G$ and $v_\star \in
\overline{\phi_\star(\bbZ\bf{n}_3)}$. We will treat the first of
these, the second being exactly similar. Suppose that $\bf{n}\in
\G$, that $\chi_\star \in \widehat{\ker q_\star}$ which we
arbitrarily extend to a measurable selection from
$\widehat{A_\star}$, that
$u_\star\in\overline{\phi_\star(\bbZ\bf{n}_2)}$ and that $b_\star$
is a solution to the corresponding equation:
\[\Delta_{u_s}(\chi_s\circ\tau)(\bf{n},z) = \Delta_{\phi_s(\bf{n})}b_s(z)\cdot c_s(z\cdot \overline{\phi_s(\bbZ\bf{n}_3)})\quad\quad\hbox{for }m_{Z_s}\hbox{-a.e. }z\in Z_s\]
for $\nu$-a.e. $s \in S$. Let $\t{u}_\star$ be any measurable lift
of $u_\star$ through $\a$ to a measurable selection from
$\overline{\t{\phi}_s(\bbZ\bf{n}_2)}\leq \t{Z}_s$. Then from the
definition of $\t{\tau}$ we have
\[
\Delta_{\t{u}_s}(\chi_s\circ\t{\tau})(\bf{n},\t{z}) =
\Delta_{u_s}(\chi_s\circ\tau)(\bf{n},z)\cdot
\Delta_{\t{u}_s}\Delta_{\t{\phi}_s(\bf{n})}b'_s(\t{z})\] where
$b'_s(\t{z})$ is the function $\t{Z}_s\to \Sone$ that corresponds to
the function
\[Z_s\ltimes A_{0,s}\to \Sone: (z,a_0)\mapsto \chi_s(\eta_s(a_0))^{-1}\]
under the above isomorphism $\t{Z}_s \leftrightarrow Z_s\ltimes
A_{0,s}$, simply because under this isomorphism the expression
$q_s(\tau(\bf{n},(s,z)))$ appearing in the definition of $\t{\tau}$
describes the lift of the rotation by $\phi_s(\bf{n})\in Z_s$ to the
rotation by $\t{\phi}_s(\bf{n})\in \t{Z}_s$.

Hence adjusting $b_\star$ to $\t{b}_\star:(s,\t{z}) \mapsto
b_s(\a(\t{z}))\cdot \Delta_{\t{u}_s}b'_s(\t{z})$ and letting
$\t{c}_s(\t{z}) := c_s(\a(\t{z}))$ we obtain a solution to the
equation
E$(u_\star,\t{\phi}_\star(\bf{n}),\overline{\t{\phi}_\star(\bbZ\bf{n}_3)},\chi_\star\circ\t{\tau}(\bf{n},\,\cdot\,))$
over the lifted system, as required. This completes the proof. \qed

\textbf{Remark}\quad We make the assumption that $\t{\pi}$ is
relatively ergodic because if we start with a non-ergodic
directional CL-extension $\bfX\to (Z_\star,m_{Z_\star},\phi_\star)$
then it will also admit many intermediate systems that are
relatively invariant over $(Z_\star,m_{Z_\star},\phi_\star)$ and are
given by some complicated combination of cosets of the Mackey group.
\fin

\begin{cor}\label{cor:joining-dCL}
Any joining of two $(\G,\bf{n}_2,\bf{n}_3)$-directional CL-systems
is a $(\G,\bf{n}_2,\bf{n}_3)$-directional CL-system.
\end{cor}

\textbf{Proof}\quad By the preceding proposition we may regard two
directional CL-systems as directional CL-extensions of their
Kronecker factors (that is, their maximal factors that are
expressible as direct integrals of ergodic group rotations). Now as
explained previously the joining of those is another direct integral
of ergodic group rotations, and over this the overall joining is
simply given as an Abelian group extension with measure supported by
some cosets of the Mackey group data inside the product of the fibre
data of the two original systems. Even if this Abelian extension is
not relatively ergodic, we can still multiply solutions to the
individual directional CL-equations to show that the directional
CL-equations for the combined cocycle also always admit solutions,
as required (once again, this is possible because we define
directional CL-cocycles by considering only their image under the
fibrewise application of an arbitrary measurable selection of fibre
group characters). \qed

Proposition~\ref{prop:dCL-change-of-base} also enables us to take
inverse limits of directional CL-systems.

\begin{cor}\label{cor:inv-lim-dCL}
Any inverse limit of $(\G,\bf{n}_2,\bf{n}_3)$-directional CL-systems
is a $(\G,\bf{n}_2,\bf{n}_3)$-directional CL-system.
\end{cor}

\textbf{Proof}\quad After using
Proposition~\ref{prop:dCL-change-of-base} to write each of our
contributing directional CL-systems as a directional CL-extension of
its Kronecker factor, this now follows from the Relative Factor
Structure Theorem by first adjoining the Kronecker
factor of the inverse limit to each individual system in the
sequence to give a new sequence expressed as an inverse limit of
directional CL-extensions of the same base Kronecker system, and
then applying the third part of Lemma~\ref{lem:dCL-mult-and-lift}.
\qed

The following is also an immediate consequence of the above
definition and results.

\begin{lem}
If $\bfX$ is a $(\G,\bf{n}_2,\bf{n}_3)$-directional CL-system then
so are almost all of its ergodic components.
\end{lem}

\textbf{Proof}\quad Indeed, upon expressing the system as
$(Z_\star,m_{Z_\star},\phi_\star)\ltimes (A_\star,m_{A_\star},\s)$
so that the invariant base space $S$ of this direct integral
coordinatizes the whole of the invariant factor, almost every
ergodic component is of the form $(Z_s,m_{Z_s},\phi_s)\ltimes
(A_s,m_{A_s},\s)$ and so is manifestly also a
$(\G,\bf{n}_2,\bf{n}_3)$-directional CL-system. \qed

With this in hand we can now prove the following useful addendum to
Theorem~\ref{thm:char-three-lines-in-2D}.

\begin{lem}\label{lem:ext-still-ergodic}
If $\bfX$ is ergodic, then the pleasant extension $\pi:\t{\bfX}\to
\bfX$ output by Theorem~\ref{thm:char-three-lines-in-2D} may also be
assumed to be ergodic.
\end{lem}

For the introduction of satedness and the definition of an FIS
system, see Subsection 3.1 of~\cite{Aus--lindeppleasant1}.

\textbf{Proof}\quad First we note that by alternately implementing
Theorem~\ref{thm:char-three-lines-in-2D} and constructing an FIS
extension and then taking an inverse limit, we may always assume
that the system output by that Theorem is FIS.

Now given an extension $\pi:\t{\bfX}\to \bfX$, if $\bfX$ is ergodic
then almost every ergodic component of $\t{\mu}$ must still push
down onto $\mu$ under $\pi$, so almost every ergodic component of
$\t{\bfX}$ still defines an extension of $\bfX$.  Let us write
$\t{\mu}_\omega$, $\omega \in \O$, for some standard Borel
parameterization of the ergodic components of $\t{\mu}$.

We next show that if $\xi_i:\bfX\to\bfY_i$ are the characteristic
factors of the original system and $\t{\xi}_i$ is the join of
isotropy and directional CL-systems appearing in the characteristic
triple for the system $\t{\bfX}$, then $\t{\xi}_i$ must still
contain $\xi_i$ for almost every $\t{\mu}_\omega$. Let $(A_m)_{m\geq
1}$ be a sequence of $\xi_i$-measurable subsets of $X$ that generate
the whole $\xi_i$-measurable $\s$-algebra up to $\mu$-negligible
sets.  Since almost every $\t{\mu}_\omega$ is still a lift of $\mu$
under $\pi$, it follows that $(\pi^{-1}(A_m))_{m\geq 1}$ still
generates the whole $(\xi_i\circ\pi)$-measurable $\s$-algebra up to
$\t{\mu}_\omega$-negligible sets for almost every $\t{\mu}_\omega$.
On the other hand, since $\xi_i \precsim \t{\xi}_i$ for $\t{\mu}$,
we know that there are corresponding $\t{\xi}_i$-measurable subsets
$B_m \subseteq \t{X}$ such that $\t{\mu}(\pi^{-1}(A_m)\triangle B_m)
= 0$ for all $m\geq 1$.  This must now also still hold for almost
every $\t{\mu}_\omega$, and so we have deduced that under almost
every $\t{\mu}_\omega$ the $\s$-algebra generated by $\t{\xi}_i$
contains that generated by $\xi_i\circ\pi$ up to negligible sets.

Finally, we observe that $(\t{\xi}_i)_\#\t{\mu}$ is a joining of
three isotropy systems and a directional CL-system, and so by the
previous lemma and its obvious analog for isotropy systems we deduce
that $\t{\xi}_i$ is also a joining of (ergodic) isotropy systems and
a directional CL-system for almost every $\t{\mu}_\omega$.

Thus we have shown that any ergodic $\bfX$ admits an ergodic
extension $(\t{X},\t{\mu}_\omega,\t{T})$ such that the
characteristic triple of factors in $\bfX$ is still determined by
the corresponding joins of systems given by
Theorem~\ref{thm:char-three-lines-in-2D}.  It is less clear that the
lifted characteristic factors $\t{\xi}_i$ are still generated by
isotropy and directional CL-systems up to negligible sets for almost
every $\t{\mu}_\omega$, but this problem can be easily repaired by
iterating this construction and then taking the (still ergodic)
inverse limit of the tower of extensions that results. \qed

By taking ergodic decompositions, it is clear that the norm
convergence asserted by Theorem~\ref{thm:polyconv} holds in general
if and only if it holds for every ergodic $\bbZ^2$-action, and given
this observation and the above lemma we will now restrict our
attention to ergodic systems for the rest of the paper.

\subsection{First reduction}

We now return to the consideration of the averages
$S_N(\cdot,\cdot)$. Our first simplification will follow from
Theorem~\ref{thm:char-poly}, giving an identification of a pair of
characteristic factors in a pleasant extension for our quadratic
averages of interest. Having obtained this, by manipulating the
classes of functions that result we will see how to simplify the
averages we need to consider even further.

\begin{thm*}
Any ergodic $\bbZ^2$-system $\bfX_0$ admits an ergodic extension
$\pi:\bfX\to \bfX_0$ in which some factor
\[\xi_1 = \xi_2 := \zeta_\pro^{T^{\bf{e}_1}}\vee \zeta_0^{T^{\bf{e}_2}}\vee \bigvee_{h\geq 1}\eta_h\]
is characteristic for the averages $S_N(\cdot,\cdot)$, where each
$\eta_h$ is a factor of $\bfX$ whose target is an
$(h\bbZ^2,h\bf{e}_1,h\bf{e}_2)$-directional CL-system for the
lattice $h\bbZ^2 := \{(hm,hn):\ m,n\in\bbZ\}$, and so
\[S_N(f_1,f_2)\sim S_N(\sfE_\mu(f_1\,|\,\xi_1),\sfE_\mu(f_2\,|\,\xi_2))\]
in $L^2(\mu)$ as $N\to\infty$ for any $f_1,f_2 \in L^\infty(\mu)$.
\end{thm*}

We will prove this in a number of steps.

\begin{lem}\label{lem:vdC-for-poly}
If \[\frac{1}{N}\sum_{n=1}^N (f_1\circ T_1^{n^2})(f_2\circ
T_1^{n^2}T_2^n) \not\to 0\] in $L^2(\mu)$ as $N\to\infty$ then there
are some $\eps > 0$ and an increasing sequence of integers $1 \leq
h_1 < h_2 < \ldots$ such that
\[\Big\|\lim_{N\to\infty}\frac{1}{N}\sum_{n=1}^N(f_1\circ T_1^{h_i^2}\circ T_1^{2h_in})(f_2\circ T_2^n)(f_2\circ (T_1^{h_i^2}T_2^{h_i})\circ (T_1^{2h_i}T_2)^n)\Big\|_2^2 \geq \eps\]
for each $i\geq 1$.
\end{lem}

\textbf{Proof}\quad Setting $u_n := (f_1\circ T_1^{n^2})(f_2\circ
T_1^{n^2}T_2^n) \in L^2(\mu)$, the version of the classical van der
Corput estimate for bounded Hilbert space sequences (see, for
instance, Section 1 of Furstenberg and Weiss~\cite{FurWei96}) shows
that
\[\frac{1}{N}\sum_{n=1}^N (f_1\circ T_1^{n^2})(f_2\circ T_1^{n^2}T_2^n) \not\to 0\]
in $L^2(\mu)$ as $N\to\infty$ only if
\begin{eqnarray*}
&&\frac{1}{H}\sum_{h=1}^H\frac{1}{N}\sum_{n=1}^N\langle
u_n,u_{n+1}\rangle\\
&&=\frac{1}{H}\sum_{h=1}^H\int_Xf_1\cdot\frac{1}{N}\sum_{n=1}^N
((f_1\circ T_1^{h^2})\circ T_1^{2hn}) (f_2\circ
T_2^n) ((f_2\circ T_1^{h^2}T_2^h)\circ T_1^{2hn}T_2^n)\,\d\mu\\
&&\not\to 0,
\end{eqnarray*}
and hence, by the Cauchy-Schwartz inequality, only if $f_1 \neq 0$
and for some $\eps
> 0$ there is an increasing sequence $1 \leq h_2 < h_2 < \ldots$
such that
\begin{eqnarray*}
&&\|f_1\|_2^2\Big\|\lim_{N\to\infty}\frac{1}{N}\sum_{n=1}^N(f_1\circ
T_1^{h_i^2}\circ T_1^{2h_in})(f_2\circ T_2^n)(f_2\circ
(T_1^{h_i^2}T_2^{h_i})\circ (T_1^{2h_i}T_2)^n)\Big\|_2^2\\
&&\geq
\Big|\int_Xf_1\cdot\Big(\lim_{N\to\infty}\frac{1}{N}\sum_{n=1}^N
(f_1\circ T_1^{h_i^2}\circ T_1^{2h_in}) (f_2\circ T_2^n) (f_2\circ
(T_1^{h_i^2}T_2^{h_i})\circ
(T_1^{2h_i}T_2)^n)\Big)\,\d\mu\Big|\\
&&\geq \|f_1\|_2^2\eps
\end{eqnarray*}
as required. \qed

In view of Theorem~\ref{thm:char-three-lines-in-2D} and a judicious
appeal to Lemma~\ref{lem:symmetry-in-the-structthm} this immediately
implies the following.

\begin{cor}\label{cor:many-condexps-large}
Any ergodic $\bbZ^2$-system $\bfX_0$ admits an ergodic extension
$\pi:\bfX\to \bfX_0$ such that if $S_N(f_1,f_2) \not\to 0$ in
$L^2(\mu)$ as $N\to\infty$ for some $f_1,f_2 \in L^\infty(\mu)$ then
there are some $\eps > 0$ and an increasing sequence of integers $1
\leq h_1 < h_2 < \ldots$ such that
\[\big\|\sfE_\mu(f_1\,|\,\zeta_0^{T_1^{2h_i}}\vee\zeta_0^{T_1^{2h_i}T_2^{-1}}\vee\zeta_0^{T_2^{-1}}\vee\eta_{1,h_i})\big\|_2^2 \geq \eps\]
and
\[\big\|\sfE_\mu(f_2\,|\,\zeta_0^{T_1^{2h_i}}\vee\zeta_0^{T_1^{2h_i}T_2}\vee\zeta_0^{T_2}\vee\eta_{2,h_i})\big\|_2^2 \geq \eps\]
for each $i\geq 1$, where each $\eta_{1,h_i}$ is a factor of $\bfX$
whose target is a $((2h_i,-1),(2h_i,0),(0,-1))$-directional
CL-system and each $\eta_{2,h_i}$ is a factor whose target is a
$((2h_i,1),(2h_i,0),(0,1))$-directional CL-system (noting that for
these triples of directions all of the values $m_{ij}$ appearing in
Theorem~\ref{thm:char-three-lines-in-2D} equal $\pm 1$). \qed
\end{cor}

This corollary tells us that if $S_N(f_1,f_2) \not\to 0$ then each
of $f_1$ and $f_2$ must enjoy a large conditional expectation onto
not just one factor of $\bfX$ with a special structure, but a whole
infinite sequence of these factors.  We will now use this to cut
down the characteristic factors we need for the averages $S_N$
further by examining the possible joint distributions of the members
of these infinite families of factors. For this we need to recall
the following special property of certain Kronecker systems,
introduced in Subsection 4.8 of~\cite{Aus--lindeppleasant2}.

\begin{dfn}[DIO system]\label{dfn:DIO}
A $\bbZ^d$-Kronecker system $(Z,m_Z,\phi)$, where
$\phi:\bbZ^d\longrightarrow Z$ is a homomorphism, has the
\textbf{disjointness of independent orbits property} or is \textbf{DIO} if for any
subgroups $\G_1,\G_2 \leq \bbZ^2$ we have
\[\G_1 \cap \G_2 = \{\bf{0}\}\quad\quad\Rightarrow\quad\quad \overline{\phi(\G_1)}\cap \overline{\phi(\G_2)} = \{1_Z\}.\]
\end{dfn}

The following was Proposition 4.32 in~\cite{Aus--lindeppleasant2}:

\begin{lem}\label{lem:make-Kron-DIO}
If a $\bbZ^2$-system is FIS then its Kronecker factor is DIO, and consequently any $\bbZ^2$-Kronecker system has a Kronecker extension that is DIO. \qed
\end{lem}

We will also need the following base result on factorizing transfer
functions, which appears as Lemma 10.3 in Furstenberg and
Weiss~\cite{FurWei96}.

\begin{lem}\label{lem:FW-factorizing-transfers}
If $\bfX_i$ for $i=1,2$ are ergodic $\bbZ$-systems and
$f_i:X_i\to\Sone$ are Borel maps for which there is some Borel
$g:X_1\times X_2\to\Sone$ with $f_1\otimes f_2 = \Delta_{T_1\times
T_2}g$, $(\mu_1\otimes \mu_2)$-a.s., then in fact there are
constants $c_i \in \Sone$ and Borel maps $g_i:X_i\to\Sone$ such that
$f_i = c_i\cdot \Delta_{T_i} g_i$. \qed
\end{lem}

\begin{lem}\label{lem:joining-different-dCL}
Suppose that $h_1\neq h_2$ are distinct nonzero integers and let
$h:= \rm{l.c.m.}(h_1,h_2,h_1 - h_2,h_1 + h_2)$. Suppose that $\bfX$
is an ergodic $\bbZ^2$-system with a pair of factors
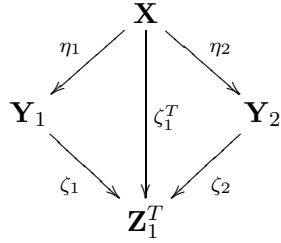
\begin{center}
$\phantom{i}$\xymatrix{ & \bfX\ar[dl]_{\eta_1}\ar[dd]^{\zeta_1^T}\ar[dr]^{\eta_2}\\
\bfY_1\ar[dr]_{\zeta_1} & & \bfY_2\ar[dl]^{\zeta_2}\\ & \bfZ_1^T }
\end{center}
such that each $\eta_i$ is an $((h_i,1),(h_i,0),(0,1))$-directional
CL-extension of $\zeta_i$, and that the Kronecker system $\bfZ_1^T$
is DIO. Then $\eta_1$ and $\eta_2$ are relatively independent under
$\mu$ over some further common factor $\eta:\bfX\to\bfY$ located as
in the diagram
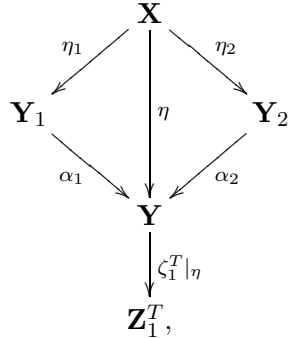
\begin{center}
$\phantom{i}$\xymatrix{ & \bfX\ar[dl]_{\eta_1}\ar[dd]^{\eta}\ar[dr]^{\eta_2}\\
\bfY_1\ar[dr]_{\a_1} & & \bfY_2\ar[dl]^{\a_2}\\ &
\bfY\ar[d]^{\zeta_1^T|_\eta}\\ & \bfZ_1^T,}
\end{center}
and where $\bfY$ is an $(h\bbZ^2,(h,0),(0,h))$-directional
CL-system.
\end{lem}

\textbf{Proof}\quad For $i=1,2$ let us pick a coordinatization
\begin{center}
$\phantom{i}$\xymatrix{ \bfY_i\ar[d]_{\zeta_i}\ar@{<->}[rr]^-\cong
& & (Z,m_Z,\phi)\ltimes (A_i,m_{A_i},\s_i)\ar[d]^{\rm{canonical}}\\
\bfZ_1^T\ar@{<->}[rr]^-\cong & & (Z,m_Z,\phi),}
\end{center}
so $\s_i$ is an $((h_i,1),(h_i,0),(0,1))$-directional CL-cocycle
over $R_\phi$.

These now combine to give a coordinatization of the target system of
the joint factor $\eta_1 \vee\eta_2$ of $\bfX$ as an extension of
$\bfZ_1^T \cong (Z,m_Z,\phi)$ by some $(R_\phi\ltimes
(\s_1,\s_2))$-invariant lift of $m_Z$ to the space $Z\ltimes
(A_1\times A_2)$. Calling this invariant lifted measure $\nu$, we
know that its two coordinate projections onto $Z\ltimes A_i$ must be
simply $m_Z\ltimes m_{A_i}$ (since this is just the measure on the
system $\bfZ_i$), and that it is relatively ergodic for the
$\bbZ^2$-action $R_\phi\ltimes (\s_1,\s_2)$ over the canonical
factor map onto $(Z,m_Z,\phi)$, simply because the whole of $\bfX$
is ergodic.

Therefore it follows from the Mackey Theorem describing ergodic
components of isometric extensions (see Proposition 4.7
in~\cite{Aus--ergdirint}) that $\nu$ takes the form $m_Z\ltimes
m_{b(\bullet)^{-1}M}$ for some section $b:Z\to A_1\times A_2$ and
some Mackey group $M \leq A_1\times A_2$ that has full
one-dimensional projections onto $A_1$ and $A_2$.

Now, in this description of $\nu$ we are free to alter $b$ pointwise
by any $M$-valued section, and so since $M$ has full one-dimensional
projections we may assume without loss of generality that $b$ takes
values in $\{1_{A_1}\}\times A_2$. Now simply identifying
$\{1_{A_1}\}\times A_2$ with a copy of the group $A_2$, if we adjust
our above coordinatization of the extension
$\bfY_i\stackrel{\zeta_i}{\longrightarrow} \bfZ_1^T$ by fibrewise
rotation by $b(\bullet)^{-1}$ we obtain a new coordinatization of
this extension by a compact Abelian group and cocycle with all the
properties of our initially-chosen coordinatization, and such that
the resulting Mackey data of the combined coordinatization has $b
\equiv 1_{A_1\times A_2}$.

Re-assigning our initial notation to this new coordinatization, we
now have $\nu = m_Z\ltimes m_M$ for some fixed $M \leq A_1\times
A_2$. It follows that the two coordinate-projection factors of the
joined system $(Z\ltimes (A_1\times A_2),\nu,R_\phi\ltimes
(\s_1,\s_2))$ onto $Z\ltimes A_i$ are relatively independent over
their further factors given by the maps
\[Z\ltimes A_i\to Z\ltimes (A_i/M_i):(z,a)\mapsto (z,aM_i)\]
where $M_i$ for $i=1,2$ are the one-dimensional slices of the Mackey
group $M$.  Moreover, the targets of these two factor maps are
identified within $(Z\ltimes (A_1\times A_2),\nu,R_\phi\ltimes
(\s_1,\s_2))$ (and hence within $\bfX$), because $M/(M_1\times M_2)$
is now a subgroup of $(A_1/M_1)\times (A_2/M_2)$ that has full
one-dimensional projections and trivial slices, and therefore
defines the graph of an isomorphism. This common target therefore
specifies some common Abelian subextension $\eta_1,\eta_2 \succsim
\eta \succsim \zeta_1^T$ over which the $\eta_i$ are relatively
independent.

This identifies the factor $\eta$ promised by the proposition; it
remains to show that its target is an
$(h\bbZ^2,(h,0),(0,h))$-directional CL-system.

First let $A \cong A_1/M_1\cong A_2/M_2$ be the fibre group of some
coordinatization of $\eta$ over $\zeta_1^T$, $q_i:A_i\onto A$ a
continuous epimorphism that corresponds to quotienting by the
subgroup $M_i$, and $\s:\bbZ^2\times Z\to A$ the cocycle over
$R_\phi$ of this coordinatization (so $\s = q_i\circ\s_i$ for
$i=1,2$). Now let $\chi \in \widehat{A}$, and let $\chi_i :=
\chi\circ q_i \in \widehat{A_i}$ for $i=1,2$.

For any $u \in \overline{\phi(\bbZ\cdot(0,1))}$ the equation
E$(u,\phi(h_i,1),\overline{\phi(\bbZ\cdot
(h_i,0))},\chi_i\circ\s_i((h_i,1),\,\cdot\,))$ gives a solution
$b_{i,u}:Z\to \Sone$ together with a one-dimensional auxiliary
$c_{i,u}:Z/\overline{\phi(\bbZ\cdot (h_i,0))}\to \Sone$ such that
\[\Delta_u\chi_i(\s_i((h_i,1),z)) = \Delta_{\phi(h_i,1)}b_{i,u}(z)\cdot c_{i,u}(z\cdot\overline{\phi(\bbZ\cdot(h_i,0))}),\]
and hence in fact
\[\Delta_u\chi(\s((h_i,1),z)) = \Delta_{\phi(h_i,1)}b_{i,u}(z)\cdot c_{i,u}(z\cdot\overline{\phi(\bbZ\cdot(h_i,0))})\]
for $i=1,2$, because $\chi\circ \s = \chi\circ q_i\circ \s_i =
\chi_i\circ \s_i$.  We will show that by modifying $b_{i,u}$ for
either $i=1$ or $i=2$ we can produce a map that simultaneously
satisfies the equations
E$(u,\phi(\bf{n}),\overline{\phi(\bbZ\cdot(h,0))},\chi\circ\s(\bf{n},\,\cdot\,))$
for all $\bf{n} \in h\bbZ^2$.  Since the case of any $v \in
\overline{\phi(\bbZ\cdot(h_1,0))}\cap\overline{\phi(\bbZ\cdot(h_2,0))}\supseteq
\overline{\phi(\bbZ\cdot(h,0))}$ is symmetrical, this will complete
the proof.

We can apply the differencing operator $\Delta_{\phi(h,0)}$ to the
above equation to obtain
\[\Delta_u\Delta_{\phi(h,0)}\chi(\s((h_i,1),z)) = \Delta_{\phi(h_i,1)}\Delta_{\phi(h,0)}b_{i,u}(z),\]
where we have used the commutativity of differencing and the fact
that $(h,0)\in \bbZ\cdot(h_i,0)$ and so
\[\Delta_{\phi(h,0)}c_{i,u}(z\cdot\overline{\phi(\bbZ\cdot(h_i,0))}) \equiv 1.\]

On the other hand, we can now appeal to the cocycle equation
$\Delta_{\phi(h,0)}\chi(\s((h_i,1),z)) =
\Delta_{\phi(h_i,1)}\chi(\s((h,0),z))$ to re-write the above as
\[\Delta_{\phi(h_i,1)}\big(\Delta_u\chi(\s((h,0),z))\cdot\Delta_{\phi(h,0)}b_{i,u}(z)^{-1}\big) \equiv 1,\]
and so we can write
\[\Delta_u\chi(\s((h,0),z))\cdot\Delta_{\phi(h,0)}b_{i,u}(z)^{-1} = f_{i,u}(z\cdot\overline{\phi(\bbZ\cdot(h_i,1))}),\]
for some $f_{i,u}:Z/\overline{\phi(\bbZ\cdot(h_i,1))}\to \Sone$.

Finally, taking the difference of these last equations for $i=1$ and
for $i=2$ we find
\[\Delta_{\phi(h,0)}(b_{2,u}\cdot b_{1,u}^{-1}) = (f_{1,u}\circ r_1)\cdot \overline{(f_{2,u}\circ r_2)}\]
where $r_i$ is the quotient epimorphism $Z\to
Z/\overline{\phi(\bbZ\cdot(h_i,1))}$.

Now, on the one hand $(h,0)\in \bbZ\cdot(h_1,1) + \bbZ\cdot(h_2,1)$,
and on the other we know that $\overline{\phi(\bbZ\cdot
(h_1,1))}\cap\overline{\phi(\bbZ\cdot (h_2,1))} = \{1_Z\}$ by the
DIO assumption.  Therefore we can analyze the above equation by
applying Lemma~\ref{lem:FW-factorizing-transfers} for each pair of
ergodic components of the restrictions $R_{\phi(h,0)}|_{r_i}$,
$i=1,2$, since the disjointness of the two orbit-closures tells us
that the above equation restricts to a combined coboundary equation
simply on the direct product of those two ergodic components. This
tells us that in fact the function $f_{i,u}(r_i(z))$ must take the
form
\[\Delta_{\phi(h,0)}(b'_{i,u}\circ r_i(z))\cdot g_{i,u}(z\cdot\overline{\phi(\G')})\]
for some Borel maps $b'_{i,u}:Z/\overline{\phi(\bbZ\cdot
(h_i,1))}\to \Sone$ and $g_{i,u}:Z/\overline{\phi(\G')}\to \Sone$
$\G':= \bbZ\cdot(h_1,1) + \bbZ\cdot(h,0)$.  Since $\G' \supseteq
h\bbZ^2$, we may instead regard $g_{i,u}$ as a map
$Z/\overline{\phi(h\bbZ^2)}\to\Sone$ and write the above function as
\[\Delta_{\phi(h,0)}(b'_{i,u}\circ r_i(z))\cdot g_{i,u}(z\cdot\overline{\phi(h\bbZ^2)}).\]
It also follows easily from the Measurable Selector Theorem that we
can take the above equations to hold for Haar-a.e. $u$ using Borel
selections $u\mapsto b_{i,u},g_{i,u}$.

Now, clearly $b'_{i,u}\circ r_i$ is invariant under
$R_{\phi(h_i,1)}$, and so $\Delta_{\phi(h_i,1)}(b_{i,u}\cdot
(b'_{i,u}\circ r_i)) = \Delta_{\phi(h_i,1)}b_{i,u}$.  This means we
can simply replace $b_{i,u}$ with $(b_{i,u}\cdot (b'_{i,u}\circ
r_i))$ in our original directional Conze-Lesigne equation, and hence
assume that the solutions we obtained for that equation also satisfy
\[\Delta_u\chi(\s((h,0),z))\cdot\Delta_{\phi(h,0)}b_{i,u}(z)^{-1} = g_{i,u}(z\cdot\overline{\phi(h\bbZ^2)}).\]

However, this now re-arranges into the form
\[\Delta_u\chi(\s((h,0),z)) = \Delta_{\phi(h,0)}b_{i,u}(z) \cdot g_{i,u}(z\cdot\overline{\phi(h\bbZ^2)}),\]
and so since $\bbZ\cdot(0,h)\subseteq h\bbZ^2$ and
$\bbZ\cdot(0,h)\subseteq \bbZ\cdot(0,h_i)$, this new version of
$b_{i,u}$ is a solution to both the originally-assumed equation
\[\rm{E}(u,\phi(h_i,1),\overline{\phi(\bbZ\cdot (h,0))},\chi\circ\s((h_i,1),\cdot))\]
and also the equation \[\rm{E}(u,\phi(h,0),\overline{\phi(\bbZ\cdot
(h,0))},\chi\circ\s((h,0),\cdot))\] (with different one-dimensional
auxiliaries), and so is actually a solution to
\[\rm{E}(u,\phi(\bf{n}),\overline{\phi(\bbZ\cdot (h,0))},\chi\circ\s(\bf{n},\cdot))\]
for every $\bf{n} \in \bbZ\cdot(h,0) + \bbZ\cdot(h_i,1) \supseteq
h\bbZ^2$, as required. \qed

We will shortly use the above lemma to examine the joint
distributions of the families of characteristic factors obtained
from Corollary~\ref{cor:many-condexps-large}.  However, before doing
so we record the following corollary of the above proof, which will
be useful later.

\begin{cor}\label{cor:str-dCL-nice}
If $\s:\bbZ^2\times Z\to\Sone$ is a
$(\G,\bf{n}_2,\bf{n}_3)$-directional CL-cocycle over a DIO system
where $\bf{n}_2,\bf{n}_3 \in \G$, then for $\{i,j\}=\{2,3\}$ there
are Borel maps $b_i:K_{\bf{n}_i}\times Z\to\Sone$ and
$c:K_{\bf{n}_i}\times Z/\overline{\phi(\G)}\to \Sone$ such that
\[\Delta_u \s(\bf{n}_j,\cdot\,) = b_i(u,z\cdot\phi(\bf{n}_j))\cdot\overline{b_i(u,z)}\cdot c(u,z\cdot\overline{\phi(\G)})\]
for Haar-almost every $(u,z)\in K_{\bf{n}_i}\times Z$.
\end{cor}

\textbf{Proof}\quad For a fixed $u \in K_{\bf{n}_i}$, the
construction of the new function $g_{i,u}$ in the previous proof
shows that we may find a solution together with a one-dimensional
auxiliary $c_u$ for the directional CL-equation
E$(u,\phi(\bf{n}_j),\overline{\phi(\G)},\chi\circ\s(\bf{n}_j,\,\cdot\,))$
--- in particular, such that $c_u(z)$ actually depends only on the
coset $z\cdot \overline{\phi(\G)}$.

It now follows from a simple measurable selection argument applied
to the collection
\begin{multline*}
\big\{(u,b',c')\in K_{\bf{n}_i}\times \C(K_{\bf{n}_i}\times
Z)\times \C(K_{\bf{n}_i}\times Z/\overline{\phi(\G)}):\\
\Delta_u \s(\bf{n}_j,\cdot\,) =
b'(z\cdot\phi(\bf{n}_j))\cdot\overline{b'(z)}\cdot
c'(z\cdot\overline{\phi(\G)})\big\}
\end{multline*}
(where as usual $\C(U)$ denotes the Polish group of equivalence
classes of Borel maps $U\to \Sone$ under $m_U$-a.e. agreement,
endowed with the topology of convergence in probability) that we may
take a selection of maps $b_{i,u}$ and $c_{i,u}$ that is Borel in
$u$ and satisfies this almost-sure equation for a.e. $u$.

It remains to obtain measurable functions $b_i$ on
$K_{\bf{n}_i}\times Z$ and $c$ on $K_{\bf{n}_i}\times
Z/\overline{\phi(\G)}$ such that $b_i(u,z) = b_{i,u}(z)$ and
$c_i(u,z\overline{\phi(\G)}) = c_{i,u}(z\overline{\phi(\G)})$ for
a.e. $(u,z)$ and hence that satisfy the desired equation Haar-almost
everywhere. This can be done, for example, by identifying $(Z,m_Z)$
with $([0,1),\rm{Lebesgue})$ as standard Borel probability spaces
and then defining $b_i(u,z)$ as the pointwise limit of the
(well-defined) averages of $b_{i,u}$ over increasingly short dyadic
intervals of values of $z$.  By the Lebesgue Density Theorem these
averages converge almost everywhere, and the resulting pointwise
limit function is clearly jointly measurable in $(u,z)$ and agrees
with $b_{i,u}$ almost surely for almost every $u$. A similar
construction applies to $c_{i,u}$, and we can make these functions
Borel by making one further modification on a negligible set. \qed

The immediate application we have for
Lemma~\ref{lem:joining-different-dCL} will require also some basic
results on the possible distributions of collections of
one-dimensional isotropy factors of a $\bbZ^2$-system.

\begin{lem}\label{lem:joining-Kron-to-isotropies}
Suppose that $\bf{n}_1$, $\bf{n}_2$, $\bf{n}_3 \in \bbZ^2 \setminus
\{\bs{0}\}$ are three directions no two of which are parallel, that
$\bfX_1 = (X_1,\mu_1,T_1)\in \sfZ_0^{\bf{n}_1}$, $\bfX_2 =
(X_2,\mu_2,T_2)\in \sfZ_0^{\bf{n}_2}$, $\bfX_3 = (X_3,\mu_3,T_3)\in
\sfZ_0^{\bf{n}_3}$ and that $\bfZ = (Z,\nu,S)$ is a group rotation
$\bbZ^2$-system. Suppose further that $\bfX = (X,\mu,T)$ is a
joining of these four systems through the factor maps
$\xi_i:\bfX\to\bfX_i$, $i=1,2,3$ and $\a:\bfX\to\bfZ$. Then
$(\xi_1,\xi_2,\xi_3,\a)$ are relatively independent under $\mu$ over
their further factors
$(\zeta_1^{T_1}\circ\xi_1,\zeta_1^{T_2}\circ\xi_2,\zeta_1^{T_3}\circ\xi_3,\a)$.
\end{lem}

\textbf{Proof}\quad  We will prove that under $\bfX$ the factors
$\xi_1$, $\xi_2$, $\xi_3$ and $\a$ are relatively independent over
$\zeta_1^{T_1}\circ\xi_1$, $\xi_2$, $\xi_3$ and $\a$; repeating this
argument to handle $\xi_2$ and $\xi_3$ then gives the full result.

Letting $\bfY = (\xi_3\vee\a)(\bfX)$ be the factor of $\bfX$
generated by $\xi_3$ (which is $T^{\bf{n}_3}$-invariant) and $\a$
(which is isometric for $T$, hence certainly for $T^{\bf{n}_3}$), we
see that this is a $T^{\bf{n}_3}$-isometric system.  This implies
that its joining to any other system is relatively independent over
the maximal $T^{\bf{n}_3}$-isometric factor of that other system.

On the other hand, $\xi_1$ and $\xi_2$ must be relatively
independent over $\xi_1\wedge \xi_2$ under $\mu$ (simply by
averaging with respect to $\bf{n}_2$), and that the subactions
generated by both $\bf{n}_1$ and $\bf{n}_2$ are trivial on this
meet, so $\xi_1\wedge\xi_2 \precsim
\zeta_0^{T^{\bf{n}_1},T^{\bf{n}_2}}$, whose target system is a
direct integral of finite group rotations factoring through the
quotient $\bbZ^2/(\bbZ\bf{n}_1 + \bbZ\bf{n}_2)$.

Since $\xi_1\vee \xi_2$ must be joined to $\xi_3\vee\a$ relatively
independently over the maximal $T^{\bf{n}_3}$-isometric factor of
$\xi_1\vee \xi_2$, it follows from the Furstenberg-Zimmer Structure
Theorem (recalled as Theorem 2.4 in~\cite{Aus--lindeppleasant2})
that $\xi_1\vee\xi_2$ is in particular joined to $\xi_3\vee\a$
relatively independently over the join of maximal isometric
subextensions
\[(\zeta_{1/(\xi_1\wedge\xi_2)|_{\xi_1}}^{T_1^{\bf{n}_3}}\circ\xi_1)\vee(\zeta_{1/(\xi_1\wedge\xi_2)|_{\xi_2}}^{T_2^{\bf{n}_3}}\circ\xi_2).\]
Since $\xi_1\wedge\xi_2$ has target a direct integral of
\emph{periodic} rotations, the maximal $T_i^{\bf{n}_3}$-isometric
subextension of $\xi_i\to (\xi_1\wedge\xi_2)|_{\xi_i}$ is simply the
maximal factor of $\xi_i$ that is coordinatizable as a direct
integral of group rotations for each $i=1,2$: that is, it is
$\zeta_1^{T_i}\circ\xi_i$. Hence we have shown that under $\mu$ the
factors $\xi_1\vee\xi_2$ and $\xi_3\vee\a$ are relatively
independent over
$(\zeta_1^{T_1}\circ\xi_1)\vee(\zeta_1^{T_2}\circ\xi_2)$ and
$\xi_3\vee\a$.  Thus whenever $f_i \in L^\infty(\mu_i)$ for
$i=1,2,3$ and $g\in L^\infty(\nu)$ we have
\begin{eqnarray*}
&&\int_X
(f_1\circ\xi_1)\cdot(f_2\circ\xi_2)\cdot(f_3\circ\xi_3)\cdot(g\circ\a)\,\d\mu\\
&&= \int_X
\sfE_\mu\big((f_1\circ\xi_1)\cdot(f_2\circ\xi_2)\,\big|\,(\zeta_1^{T_1}\circ\xi_1)\vee(\zeta_1^{T_2}\circ\xi_2)\big)\cdot(f_3\circ\xi_3)\cdot(g\circ\a)\,\d\mu\\
&&= \int_X
(\sfE_\mu(f_1\,|\,\zeta_1^{T_1})\circ\xi_1)\cdot(\sfE_\mu(f_2\,|\,\zeta_1^{T_2})\circ\xi_2)\cdot(f_3\circ\xi_3)\cdot(g\circ\a)\,\d\mu\\
&&= \int_X
(\sfE_\mu(f_1\,|\,\zeta_1^{T_1})\circ\xi_1)\cdot(f_2\circ\xi_2)\cdot(f_3\circ\xi_3)\cdot(g\circ\a)\,\d\mu,
\end{eqnarray*}
where the second equality follows from the relative independence of
$\xi_1$ and $\xi_2$ over $\xi_1\wedge \xi_2$, which is contained in
$\zeta_1^{T_i}\circ\xi_i$ for both $i=1,2$. This completes the
proof. \qed

\begin{lem}\label{lem:another-rel-indep-result-for-1Ds}
Suppose that $\bf{n}_1,\bf{n}_2,\bf{n}_3,\bf{n}_4 \in
\bbZ^2\setminus \{\bs{0}\}$ are directions no two of which are
parallel, that $\bfX_i = (X_i,\mu_i,T_i)\in \sfZ_0^{\bf{n}_i}$ for
$i=1,2,3,4$ and that $\bfY = (Y,\nu,S)$ is a two-step Abelian
isometric $\bbZ^2$-system. Suppose further that $\bfX = (X,\mu,T)$
is a joining of these five systems through the factor maps
$\xi_i:\bfX\to \bfX_i$, $i=1,2,3,4$ and $\eta:\bfX\to\bfY$, with the
maximality properties that $\xi_i = \zeta_0^{T^{\bf{n}_i}}$ for
$i=1,2,3,4$ and $\eta \succsim \zeta_1^T$. Then
$(\xi_1,\xi_2,\xi_3,\xi_4,\eta)$ are relatively independent under
$\mu$ over their further factors $(\zeta_{\Ab,2}^{T_1}\circ
\xi_1,\zeta_{\Ab,2}^{T_2}\circ \xi_2,\zeta_{\Ab,2}^{T_3}\circ
\xi_3,\zeta_{\Ab,2}^{T_4}\circ\xi_4,\eta)$.
\end{lem}

\textbf{Proof}\quad First set $\b_i := \zeta_2^{T_i}\circ\xi_i$ and
$\a_i := \zeta_{\Ab,2}^{T_i}\circ \xi_i$ for $i=1,2,3,4$, so each
$\a_i\succsim \zeta_1^{T_i}\circ\xi_i$ is the maximal Abelian
subextension of $\b_i\succsim \zeta_1^{T_i}\circ\xi_i$.

We need to prove that
\[\int_X f_1f_2f_3f_4g\,\d\mu = \int_X \sfE_\mu(f_1\,|\,\a_1)\sfE_\mu(f_2\,|\,\a_2)\sfE_\mu(f_3\,|\,\a_3)\sfE_\mu(f_4\,|\,\a_4)g\,\d\mu\] for any
$\xi_i$-measurable functions $f_i$ and $\eta$-measurable function
$g$.  In fact it will suffice to prove that
\[\int_X f_1f_2f_3f_4g\,\d\mu = \int_X f_1f_2f_3\sfE_\mu(f_4\,|\,\a_4)g\,\d\mu,\]
since then repeating the same argument for the other three isotropy
factors in turn completes the proof.

By Lemma~\ref{lem:joining-Kron-to-isotropies} the three factors
$\zeta_1^T\vee\xi_1$, $\zeta_1^T\vee \xi_2$ and $\zeta_1^T\vee\xi_3$
must be joined relatively independently over $\zeta_1^T$.  On the
other hand, the factor $\xi_4\vee\eta$ is an extension of
$\zeta_1^T$ that is certainly still an Abelian isometric extension
for the $(\bbZ\bf{n}_4)$-subaction, and so $\xi_1\vee
\xi_2\vee\xi_3\vee\zeta_1^T$ must be joined to it relatively
independently over
\[\zeta_2^{T^{\bf{n}_4}}\wedge\big(\xi_1\vee
\xi_2\vee\xi_3\vee\zeta_1^T\big).\]

However, now the Furstenberg-Zimmer Structure Theorem tells us that
this last factor must be contained in
\[(\zeta_2^{T^{\bf{n}_4}}\wedge\xi_1)\vee(\zeta_2^{T^{\bf{n}_4}}\wedge \xi_2)\vee(\zeta_2^{T^{\bf{n}_4}}\wedge\xi_3)\vee\zeta_1^T\]
(using that $\zeta_2^{T^{\bf{n}_4}}\wedge (\xi_i\vee \zeta_1^T) =
(\zeta_2^{T^{\bf{n}_4}}\wedge \xi_i)\vee \zeta_1^T$, because
$\zeta_1^T$ is already one-step distal). Here the factors
$\zeta_2^{T^{\bf{n}_4}}\wedge \xi_i$ are actually isometric
extensions of $\zeta_1^T\wedge \xi_i$ (not just of
$\zeta_1^{T^{\bf{n}_4}}\wedge\xi_i$), since in each case
isometricity for the $(\bbZ\bf{n}_4)$-subaction and
\emph{invariance} for the $(\bbZ\bf{n}_i)$-subaction together imply
isometricity for the whole $\bbZ^2$-system
$\zeta_1^{T^{\bf{n}_4}}\wedge\xi_i$, since $\bbZ\bf{n}_i +
\bbZ\bf{n}_4$ has finite index in $\bbZ^2$ by the non-parallel
assumption.

Overall this tells us that $\xi_4\vee\eta$ is relatively independent
from the factors $\xi_1$, $\xi_2$ and $\xi_3$ over their further
factors $\b_1$, $\b_2$ and $\b_3$; and now applying the same
argument with any of the other isotropy factors as the distinguished
factor in place of $\xi_4$, we deduce that this latter is relatively
independent from all our other factors over $\b_4$.

By reducing to the factor of $\bfX$ generated by the $\b_i$ and
$\eta$, we may therefore assume that each $\bfX_i$ is itself a
two-step distal system (since the join $\b_1\vee
\b_2\vee\b_3\vee\b_4\vee\eta$ is still two-step distal, and so its
maximal isotropy factor in each direction $\bf{n}_i$ is also
two-step distal and hence equal to $\b_i$).

To make the remaining reduction to have $\a_i$ in place of $\b_i$,
now let $\bfZ_1^T = (Z_\star,m_{Z_\star},\phi_\star)$ be some
coordinatization of the Kronecker factor $\zeta_1^T$ as a direct
integral of ergodic $\bbZ^2$-group rotations, and let us pick
coordinatizations
\begin{center}
$\phantom{x}$\xymatrix{\bfX_i\ar[dr]_{\zeta_1^T|_{\xi_i}}\ar@{<->}[rr]^-\cong
&& \bfZ_1^{T_i}\ltimes
(G_{i,\bullet}/H_{i,\bullet},m_{G_{i,\bullet}/H_{i,\bullet}},\s_i)\ar[dl]^{\rm{canonical}}\\
& \bfZ_1^{T_i}  }
\end{center}
and
\begin{center}
$\phantom{x}$\xymatrix{\bfY\ar[dr]_{\zeta_1^T|_\eta}\ar@{<->}[rr]^-\cong
&& \bfZ_1^T\ltimes
(A_\bullet,m_{A_\bullet},\tau)\ar[dl]^{\rm{canonical}}\\
& \bfZ_1^T.  }
\end{center}

We know this may be done so that the $\s_i$ and $\tau$ are
relatively ergodic, and so now replacing each $\bfX_i$ with its
covering group extension and joining these relatively independently
over the joining $\bfX$ of the $\bfX_i$'s and $\bfY$, we reduce the
problem to the case in which $H_{i,\bullet} =
\{1_{G_{i,\bullet}}\}$.

Given this we know that any joining of the above relatively ergodic
group extensions of $\bfZ_1^T$ is described by some
$T|_{\zeta_1^T}$-invariant measurable Mackey group data
\[M_z \leq \prod_{i=1}^4 G_{i,z_i}\times A_z\]
and a section $b:Z\to \prod_{i=1}^4 G_{i,z_i}\times A_z$, where $z
\in Z_\star$ and $z_i = \zeta_0^{T^{\bf{n}_i}}|_{\zeta_1^T}(z)$.  To
complete the proof we will show that
\[M_z \geq \prod_{i=1}^4[G_{i,z_i},G_{i,z_i}]\times \{1_{A_z}\}\]
almost surely, since in this case we may quotient out each extension
$\bfX_i\to \bfZ_1^{T_i}$ fibrewise by the normal subgroups
$[G_{i,\bullet},G_{i,\bullet}] \leq G_{i,\bullet}$ to obtain that
our joining is relatively independent over some Abelian
subextensions, as required.

The point is that for any three-subset $\{i_1,i_2,i_3\}\subset
\{1,2,3,4\}$ the projection of $M_\bullet$ onto the product of
factor groups $G_{i_j,z_{i_j}}$, $j=1,2,3$ is just the Mackey group
data of the joining of $\xi_{i_1}$, $\xi_{i_2}$, $\xi_{i_3}$ and
$\zeta_1^T$ as factors of $\bfX$.  By
Lemma~\ref{lem:joining-Kron-to-isotropies} these are relatively
independent over $\zeta_1^T$, so this coordinate projection of the
Mackey group must be the whole of $\prod_{j=1}^3G_{i_j,z_{i_j}}$.
Hence $M_\bullet$ has full projections onto any three of
$G_{i,z_i}$, and so for any $g_1,h_1 \in G_{1,z_1}$ (say) we can
find $g_2 \in G_{2,z_2}$, $h_3 \in G_{3,z_3}$ and $a,b \in A_z$ such
that
\begin{multline*}
(g_1,g_2,1,1,a),(h_1,1,h_3,1,b) \in M_z\\
\Rightarrow\quad\quad [(g_1,g_2,1,1,a),(h_1,1,h_3,1,b)] =
([g_1,h_1],1,1,1,1)\in M_z.
\end{multline*}

Arguing similarly for the other $G_{i,z_i}$, we deduce that
$M_\bullet$ contains the Cartesian product of commutator subgroups,
as required. \qed

\begin{prop}\label{prop:jointdists}
If $h_1$, $h_2$ and $h$ are integers as in
Lemma~\ref{lem:joining-different-dCL}, $\bfX$ is an ergodic
$\bbZ^2$-system whose Kronecker factor $\zeta_1^T:\bfX\to \bfZ_1^T$
is DIO and $\eta_i:\bfX\to\bfY_i$ is an
$((h_i,1),(h_i,0),(0,1))$-directional CL-extension of $\zeta_1^T$
for $i=1,2$, then the two factors
\[\zeta_0^{T_1^{h}}\vee \zeta_0^{T_2}\vee \zeta_0^{T_1^{h_i}T_2}\vee
\eta_i\quad\quad i=1,2\] of $\bfX$ are relatively independent over a
common further factor of the form $\zeta_0^{T_1^h}\vee
\zeta_0^{T_2}\vee \eta$ where $\eta$ has target an
$(h\bbZ^2,(h,0),(0,h))$-directional CL-system.
\end{prop}

\textbf{Proof}\quad Since $\eta_1\vee\eta_2$ still has target a
two-step Abelian isometric system, the preceding lemma shows that
$\zeta_0^{T_1^{h}}$, $\zeta_0^{T_2}$, $\zeta_0^{T_1^{h_1}T_2}$,
$\zeta_0^{T_1^{h_2}T_2}$ and $\eta_1\vee\eta_2$ are all relatively
independent over their maximal two-step Abelian factors. Denoting
the first four of these by $\a_1$, $\a_2$, $\a_{12,1}$ and
$\a_{12,2}$ respectively, it will therefore suffice to prove that
$\a_1\vee\a_{12,1}\vee\a_2\vee\eta_1$ and
$\a_1\vee\a_{12,2}\vee\a_2\vee\eta_2$ are relatively independent
over some further common factor $\a_1\vee\a_2\vee\eta$ with $\eta$ a
directional CL-factor of the kind asserted.

However, as described following the introduction of directional
CL-systems in Subsection 3.6 of~\cite{Aus--lindeppleasant2}, each
$\a_1\vee\a_{12,i}\vee\a_2\vee\eta_i$ is itself still an
$((h_i,1),(h_i,0),(0,1))$-directional CL-system, and so this latter
assertion follows at once from
Lemma~\ref{lem:joining-different-dCL}.  This completes the proof.
\qed

We can now make use of the above-found relative independence through
the following simple lemma.

\begin{lem}\label{lem:rel-ind-factor-seq} Suppose that $(X,\mu)$ is a
standard Borel probability space, $\pi_n:X\to Y_n$ is a sequence of
factor maps of $X$ and $\a_n:Y_n\to Z_n$ is a sequence of further
factor maps of $Y_n$ such that $(\pi_n,\pi_m)$ are relatively
independent over $(\a_n\circ\pi_n,\a_m\circ\pi_m)$ whenever $n \neq
m$ (note that we do not require such relative independence for more
than two of the $\pi_i$ at once).  If $f \in L^\infty(\mu)$ is such
that $\limsup_{n\to\infty}\|\sfE_\mu(f\,|\,\pi_n)\|_2 > 0$, then
also $\limsup_{n\to\infty}\|\sfE_\mu(f\,|\,\a_n)\|_2 > 0$.
\end{lem}

\textbf{Proof}\quad By thinning out our sequence if necessary, we
may assume that for some $\eta > 0$ we have
$\|\sfE_\mu(f\,|\,\pi_n)\|_2 \geq \eta$ for all $n$.  Suppose, for
the sake of contradiction, that $\sfE_\mu(f\,|\,\a_n) \to 0$ as
$n\to\infty$. Consider the sequence of Hilbert subspaces $L_n \leq
L^2(\mu)$ comprising those functions that are $\pi_n$-measurable and
the further subspaces $K_n \leq L_n$ comprising those that are
$\a_n$-measurable.  Then by assumption all the subspaces $L_n\ominus
K_n$ are mutually orthogonal, but $f$ has orthogonal projection of
norm at least $\eta/2$ onto all but finitely many of them, which is
clearly impossible. \qed

\textbf{Proof of Theorem~\ref{thm:char-poly}}\quad Letting
$\pi:\bfX\to \bfX_0$ be the ergodic pleasant extension for triple
linear averages in general position obtained by applying
Theorem~\ref{thm:char-three-lines-in-2D} and
Lemma~\ref{lem:ext-still-ergodic} and then making a further
extension of the Kronecker factor using
Lemma~\ref{lem:make-Kron-DIO} if necessary, now
Corollary~\ref{cor:many-condexps-large},
Proposition~\ref{prop:jointdists} and
Lemma~\ref{lem:rel-ind-factor-seq} show that whenever $f_1,f_2 \in
L^\infty(\mu)$ have $S_n(f_1,f_2)\not\to 0$, they also satisfy
$\sfE_\mu(f_i\,|\,\zeta_\pro^{T_1}\vee\zeta_0^{T_2}\vee\eta_\infty)\not=
0$ where $\zeta_\pro^{T_1}$ is the factor generated by all
$\zeta_0^{T_1^h}$, $h\geq 1$, and $\eta_\infty$ is a join over some
sequence of integers $h$ of $(h\bbZ^2,(h,0),(0,h))$-directional
CL-factors. Writing $\xi:=
\zeta_\pro^{T_1}\vee\zeta_0^{T_2}\vee\eta_\infty$, the proposition
follows at once by considering the decomposition
\begin{multline*}
S_N(f_1,f_2) = S_N(\sfE_\mu(f_1\,|\,\xi),\sfE_\mu(f_2\,|\,\xi))\\ +
S_N(f_1 - \sfE_\mu(f_1\,|\,\xi),\sfE_\mu(f_2\,|\,\xi)) + S_N(f_1,f_2
- \sfE_\mu(f_2\,|\,\xi)).
\end{multline*}
\qed

\subsection{Second reduction}

Theorem~\ref{thm:char-poly} shows that Theorem~\ref{thm:polyconv}
will follow if we prove that $S_N(f_1,f_2)$ converges whenever $f_i$
is $\xi_i$-measurable. By approximation in $L^2(\mu)$ and
multilinearity, it actually suffices to consider the averages
$S_N(f_{11} f_{12} g_1,f_{21}f_{22}g_2)$ in which each $f_{j1}$ is
$\zeta_0^{T^\ell_1}$-measurable for some large $\ell\geq 1$, each
$f_{j2}$ is $\zeta_0^{T_2}$-measurable and each $g_j$ is
$\eta$-measurable for some $(h\bbZ^2,(h,0),(0,h))$-directional
CL-factor $\eta$ for some large $h\geq 1$.

Next, writing
\begin{eqnarray*}
&&S_N(f_{11} f_{12}g_1,f_{21}f_{22}g_2) =
\frac{1}{N}\sum_{n=1}^N((f_{11}\cdot f_{12}\cdot g_1)\circ
T_1^{n^2})((f_{21}\cdot f_{22}\cdot g_2)\circ T_1^{n^2}T_2^n)\\
&&\sim \frac{1}{\ell}\sum_{k = 0}^{\ell-1}\frac{1}{(N/\ell)}\sum_{n
= 1}^{\lfloor N/\ell\rfloor} ((f_{11}\cdot f_{12}\cdot g_1)\circ
T_1^{(\ell n+k)^2})((f_{21}\cdot
f_{22}\cdot g_2)\circ T_1^{(\ell n+k)^2}T_2^{\ell n+ k})\\
&&= \frac{1}{\ell}\sum_{k = 0}^{\ell-1}(f_{11}\circ
T_1^{k^2})\Big(\frac{1}{(N/\ell)}\sum_{n = 1}^{\lfloor
N/\ell\rfloor} ((f_{12}\cdot g_1)\circ T_1^{(\ell
n+k)^2})((f_{21}\cdot f_{22}\cdot g_2)\circ T_1^{(\ell
n+k)^2}T_2^{\ell n+ k})\Big)\\
&&= \frac{1}{\ell}\sum_{k = 0}^{\ell-1}(f_{11}\circ
T_1^{k^2})\Big(\frac{1}{(N/\ell)}\sum_{n = 1}^{\lfloor
N/\ell\rfloor} ((f_{12}\cdot f_{22}\cdot g_1)\circ T_1^{(\ell
n+k)^2})\\
&&\quad\quad\quad\quad\quad\quad\quad\quad\quad\quad\quad\quad\quad\quad\quad\quad\cdot
(g_2\circ T_1^{(\ell n+k)^2}T_2^{\ell n+ k})(f_{21}\circ
T_1^{k^2}\circ T_2^{\ell n + k})\Big)
\end{eqnarray*}
(recalling that $\sim$ denotes asymptotic agreement in $L^2(\mu)$ as
$N\to\infty$), we see that it will suffice to prove convergence in
$L^2(\mu)$ for all averages along infinite arithmetic progressions
of the form
\[\frac{1}{(N/\ell)}\sum_{n = 1}^{\lfloor
N/\ell\rfloor} ((f_{12}\cdot f_{22}\cdot g_1)\circ T_1^{(\ell n)^2
+2k(\ell n)})(g_2\circ T_1^{(\ell n)^2 + 2k(\ell n)}T_2^{\ell
n})(f_{21}\circ T_2^{\ell n})\] for all $k \in
\{0,1,\ldots,\ell-1\}$, where for a fixed $k$ we have re-written
$(f_{12}\cdot f_{22}\cdot g_1)\circ T_1^{k^2}$ as simply
$f_{12}\cdot f_{22}\cdot g_1$ and similarly for the other factors,
and have discarded the initial multiplication by the $n$-independent
function $f_{11}\circ T_1^{k^2}T_2^k$.

If we now simply re-label $T_i^\ell$ as $T_i$ (and so effectively
restrict our attention to the subaction of $\ell\bbZ^2$), then the
above averages are modified to
\[\frac{1}{(N/\ell)}\sum_{n = 1}^{\lfloor
N/\ell\rfloor} ((f_{12}\cdot f_{22}\cdot g_1)\circ T_1^{\ell n^2
+2kn})(g_2\circ T_1^{\ell n^2 + 2kn}T_2^n)(f_{21}\circ T_2^n)\] and
now $f_{21}$ is simply $T_1$-invariant.  Moreover, it is clear that
any $(h\bbZ^2,(h,0),(0,h))$-directional CL-system for the action $T$
retains this property under this re-labeling (indeed, the same
property for the re-labeled system is potentially slightly weaker),
and also if we then restrict attention to any one of the (finitely
many) $\ell\bbZ^2$-ergodic components of the overall system.

Thus, we have now reduced our task to the proof of convergence for
averages of the form
\[\frac{1}{N}\sum_{n = 1}^{N} ((F_2 \cdot g_1)\circ T_1^{\ell n^2 + an})(g_2\circ
T_1^{\ell n^2 + an}T_2^n)(F_1\circ T_2^n),\] for any fixed integers
$\ell,a \geq 1$, where $F_2$ is $T_2$-invariant, $F_1$ is
$T_1$-invariant and $g_1$, $g_2$ are $\eta$-measurable.

This conclusion was obtained by simply re-writing the expression for
$S_N$ for the functions of interest to us (with a little sleight of
hand to deal with the rational spectrum of $T_1$).  However, it
turns out that we can do better still with just a little more work:
to wit, that we may also remove the function $F_2$ from
consideration, and so reduce Theorem~\ref{thm:polyconv} to
Proposition~\ref{prop:reduced-polyconv} below.  This will rely on
the following results
from~\cite{Aus--lindeppleasant1,Aus--lindeppleasant2}.

\begin{prop}[The Furstenberg self-joining controls nonconventional
averages]\label{prop:Fberg-controls-nonconv-aves} If $f_1,f_2,f_3
\in L^\infty(\mu)$ and
\[\frac{1}{N}\sum_{n=1}^N(f_1\circ T_1^{2\ell h n})(f_2\circ T_1^{2\ell h n}T_2^{-n})(f_3\circ T_2^{-n})\not\to 0\]
as $N\to\infty$, then there is some $(T_1^{2\ell h}\times T_1^{2\ell
h}T_2^{-1}\times T_2^{-1})$-invariant bounded Borel function $G:X^3
\to \bbR$ such that
\[\int_{X^3} (f_1\otimes f_2\otimes f_3)\cdot G\,\d\mu_h^\rm{F} \neq 0,\]
where $\mu_h^\rm{F} := \mu^\rm{F}_{T_1^{2\ell h},T_1^{2\ell
h}T_2^{-1}, T_2^{-1}}$ is the Furstenberg self-joining (see
Subsection 4.1 of~\cite{Aus--lindeppleasant1}).  This is a
three-fold self-joining of $(X,\mu,T_1^{2\ell h},T_2)$ that is also
invariant under the transformation $\vec{T}_h := T_1^{2\ell h}\times
T_1^{2\ell h}T_2^{-1}\times T_2^{-1}$, and has the following
properties:
\begin{itemize}
\item The restriction of $\mu_h^{\rm{F}}$ to $Z^3$ is the Haar measure $m_{Z_h}$
of some closed subgroup $Z_h \leq Z^3$, and if the Kronecker factor
$(Z,m_Z,\phi)$ of $\bfX$ is DIO then
\begin{eqnarray*}
Z_h &=& \{(z_1,z_2,z_3)\in Z^3:\ z_1z_2^{-1} \in
K_{(0,1)},\,z_1z_3^{-1} \in K_{(2\ell h,1)},\,z_2z_3^{-1}\in
K_{(2\ell h,0)}\}\\
&=& \{(zu,zuv,zv):\ z\in Z,\,u\in K_{(2\ell h,0)},\,v\in
K_{(0,1)},\,uv^{-1}\in K_{(2\ell h,1)}\},
\end{eqnarray*}
where as usual we write $K_\bf{n} := \overline{\phi(\bbZ\bf{n})}$.

\item The $\vec{T}_h$-ergodic components of the restriction of $\mu_h^{\rm{F}}$ to $(Z\ltimes A)^3$ are
almost all of the form \[m_{z_0\cdot\overline{(\phi(2\ell
h\bf{e}_1),\phi(2\ell h\bf{e}_1 -
\bf{e}_2),\phi(-\bf{e}_2))^\bbZ}}\ltimes m_{b_h(\bullet)^{-1}\cdot
M_h\cdot a}\] for some Mackey group $M_h \leq A^3$ on $Z_h$, some
Borel section $b_h:Z_h \to A^3$ and some fixed $a \in A^3$ and $z_0
\in Z_h$.
\end{itemize}
These last conclusions follow from the conjunction of Propositions
4.6 and 4.7 in~\cite{Aus--lindeppleasant2} and the discussion of
Subsection 4.8 of~\cite{Aus--lindeppleasant2}, except for the fact that the Mackey
group $M_h$ is constant which results from the presence of the
restrictions of the transformations $(T^{\bf{n}})^{\times 3}$ to
$(Z\ltimes A)^3$ that are described by $A$-valued cocycles and leave
$M_h$ invariant, as in the proof of Proposition 4.10
of~\cite{Aus--lindeppleasant2}. \qed
\end{prop}

\begin{prop}\label{prop:reduced-polyconv}
If $\bfX$ is a $\bbZ^2$-system as output by
Theorem~\ref{thm:char-poly} and $\ell,a \geq 1$ are fixed integers
then the nonconventional ergodic averages
\[S'_N(g_1,g_2,f) := \frac{1}{N}\sum_{n=1}^N(g_1\circ T_1^{\ell n^2 + an})(g_2\circ T_1^{\ell n^2 + an}T_2^n)(f\circ T_2^n)\]
converge in $L^2(\mu)$ as $N\to\infty$ whenever $g_1$, $g_2$ are
$\eta$-measurable and $f$ is $T_1$-invariant.
\end{prop}

\textbf{Proof of Theorem~\ref{thm:polyconv} from
Proposition~\ref{prop:reduced-polyconv}}\quad
Theorem~\ref{thm:char-poly} and the re-arrangement above show that
it suffices to prove convergence for averages of the form
\[\frac{1}{N}\sum_{n = 1}^{N} ((F_2 \cdot g_1)\circ T_1^{\ell n^2 + an})(g_2\circ
T_1^{\ell n^2 + an}T_2^n)(F_1\circ T_2^n)\] with $F_2$ being
$T_2$-invariant and $F_1$ being $T_1$-invariant. We will now show
that these tend to $0$ in $L^2(\mu)$ if $F_2$ is orthogonal to
$\zeta_{\Ab,2}^T$, which combined with the $T_2$-invariance of $F_2$
shows that it suffices to treat the case when $F_2$ is actually
measurable with respect to $\zeta_{\Ab,2}^T\wedge \zeta_0^{T_2}$,
which is another $(h\bbZ^2,(h,0),(0,h))$-directional CL-system and so may be subsumed into the
factor $\eta$.  The resulting averages will then be easily
re-arranged into the form $S'_N$.

By another appeal to the van der Corput estimate we know that the
above averages tend to zero in $L^2(\mu)$ unless also
\begin{eqnarray*}
&&\frac{1}{H}\frac{1}{N}\sum_{h=1}^H\sum_{n=1}^N\int_X(F_2\circ
T_1^{\ell n^2 + 2\ell hn + \ell h^2 + an + ah})(\overline{F_2}\circ
T_1^{\ell n^2 + an})(F_1\circ
T_2^{n+h})(\overline{F_1}\circ T_2^n)\\
&&\quad\quad\quad\quad\quad\quad\cdot(g_1\circ T_1^{\ell n^2+2\ell
hn + \ell h^2 + an + ah})(\overline{g_1}\circ T_1^{\ell n^2 + an})\\
&&\quad\quad\quad\quad\quad\quad\cdot(g_2\circ T_1^{\ell n^2+2\ell
hn+\ell h^2 + an + ah}T_2^{n+h})(\overline{g_2}\circ T_1^{\ell n^2 +
an}T_2^n)\,\d\mu \not\to 0
\end{eqnarray*}
as $N\to\infty$ and then $H\to\infty$.

Using the invariances of the $F_i$ we can change variables in each
of the integrals appearing above by $T_1^{-\ell n^2 - an}T_2^{-n}$
and find that the above conclusion simplifies to
\begin{eqnarray*}
&&\frac{1}{H}\frac{1}{N}\sum_{h=1}^H\sum_{n=1}^N\int_X(F_2\circ
T_1^{2\ell hn + \ell h^2 + ah})\cdot \overline{F_2}\cdot(F_1\circ
T_2^h)\cdot \overline{F_1}\\
&&\quad\quad\quad\quad\quad\quad\cdot (g_1\circ T_1^{2\ell hn + \ell
h^2 + ah}T_2^{-n})(\overline{g_1}\circ T_2^{-n})(g_2\circ
T_1^{2\ell hn+\ell h^2 + ah}T_2^h)\overline{g_2}\,\d\mu\\
&&=\frac{1}{H}\frac{1}{N}\sum_{h=1}^H\sum_{n=1}^N\int_X((F_2\cdot(g_2\circ
T_2^h))\circ T_1^{2\ell hn + \ell h^2 + ah})(g_1\circ T_1^{2\ell hn
+
\ell h^2 + ah}T_2^{-n})(\overline{g_1}\circ T_2^{-n})\\
&&\quad\quad\quad\quad\quad\quad\quad\quad\quad\quad\quad\quad\quad\quad\quad\quad\cdot
\overline{F_2}\cdot(F_1\circ T_2^h)\cdot \overline{F_1}\cdot
\overline{g_2}\,\d\mu \not\to 0.
\end{eqnarray*}
Hence, extracting the active part of the average over
$n\in\{1,2,\ldots,N\}$ it follows that for some $h\geq 1$ (here we
need only one such value) we have
\[\frac{1}{N}\sum_{n=1}^N((F_2\cdot(g_2\circ
T_2^h))\circ T_1^{2\ell hn + \ell h^2 + ah})(g_1\circ T_1^{2\ell hn
+ \ell h^2 + ah}T_2^{-n})(\overline{g_1}\circ T_2^{-n}) \not\to 0\]
in $L^2(\mu)$.

This is another instance of the kind of triple linear average that
we have considered previously, but now with functions $F_2\cdot
(g_2\circ T_2^h)$, $g_1\circ T_1^{\ell h^2 + ah}$ and
$\overline{g_1}$ that are measurable with respect to more restricted
factors of the overall system $\bfX$.  Applying
Proposition~\ref{prop:Fberg-controls-nonconv-aves} we obtain
\[\int_{X^3} \big(((F_2\cdot(g_2\circ
T_2^h))\circ T_1^{\ell h^2 + ah})\otimes (g_1\circ T_1^{\ell h^2 +
ah})\otimes \overline{g_1}\big)\cdot G\,\d\mu^{\rm{F}} \neq 0\] for
some function $G \in L^\infty(\mu^\rm{F})$ that is invariant under
$\vec{T} := T_1^{2\ell h}\times T_1^{2\ell h}T_2^{-1}\times
T_2^{-1}$.

Let $\pi_1$, $\pi_2$ and $\pi_3$ be the three coordinate projections
$X^3\to X$, and now consider on $(X^3,\mu^{\rm{F}})$ the two
$\mu^\rm{F}$-preserving transformations $\vec{T}$ and $T_2^{\times
3}$. The function $(F_2\circ T_1^{\ell h^2 + ah})\circ \pi_1$ is
$T_2^{\times 3}$-invariant (simply because $F_2$ was assumed
$T_2$-invariant), and the above nonvanishing integral asserts that
this function has a positive inner product with the function
\[(g_2\circ T_2^h\circ T_1^{\ell h^2 + ah}\circ \pi_1)\cdot (((g_1\circ T_1^{\ell h^2 + ah})\cdot \overline{g_1})\circ \pi_2)\cdot G,\]
where $g_2\circ T_2^h\circ T_1^{\ell h^2 + ah}\circ \pi_1$ and
$((g_1\circ T_1^{\ell h^2 + ah})\cdot \overline{g_1})\circ \pi_2$
are both measurable with respect to some two-step Abelian factor by
assumption and where $G$ is $\vec{T}$-invariant. Moreover $\vec{T}$
simply restricts to $T_1^{2\ell h}$ under $\pi_1$. Therefore
Lemma~\ref{lem:another-rel-indep-result-for-1Ds} above implies that
the factor $\zeta_0^{T_2}\circ \pi_1 \lesssim \zeta_0^{T_2^{\times
3}}$ of $\bfX' := (X^3,\mu^\rm{F},\vec{T},T_2^{\times 3})$ is
relatively independent from $\zeta_{\Ab,2}^{T_2^{\times 3},\vec{T}}
\vee \zeta_0^{\vec{T}}$ over the two-step Abelian factor
$\zeta_{\Ab,2}^{\vec{T},T_2^{\times 3}}\wedge \zeta_0^{T_2^{\times
3}}$. This, in turn, is a two-step Abelian isometric system on which
$T_2^{\times 3}$ is invariant, and so it must be joined to
$\zeta_0^{T_2}\circ\pi_1$ relatively independently over the maximal
two-step Abelian factor of $\zeta_0^{T_2}\circ\pi_1$. It follows
that $F_2$ must have nonzero conditional expectation onto the factor
$\zeta_{\Ab,2}^T\wedge \zeta_0^{T_2}$, as claimed.

Since this last factor is also a $(h\bbZ^2,(h,0),(0,h))$-directional
CL-system, we may assume that it is already contained in $\eta$, and
therefore we have shown that it suffices to prove convergence of our
averages when we write simply $g_1$ in place of $F_2\cdot g_1$.
These puts them into the form $S'_N(g_1,g_2,F_1)$ treated by
Proposition~\ref{prop:reduced-polyconv}, and so completes the proof.
\qed

By continuing in the vein of the above proof we could try to obtain
also a simplification of the function $F_1$. However, in fact these
methods do not seem to give a reduction of this function that is
strong enough to be useful. In the next subsections we will change
tack to give a different kind of simplification of the averages,
from which convergence can be proved given no further information
about the function $F_1$.

\subsection{Using the Mackey group of the Furstenberg
self-joining}\label{subs:backtotheMac}

The last subsection has left us to consider the averages
\[S'_N(g_1,g_2,f) := \frac{1}{N}\sum_{n=1}^N (g_1\circ T_1^{\ell n^2 + an})(g_2\circ T_1^{\ell n^2 + an}T_2^n)(f\circ T_2^n)\]
for $g_1,g_2$ that are measurable with respect to some
$(m\bbZ^2,(m,0),(0,m))$-directional CL-factor $\eta:\bfX\to\bfY$ and
$f$ that is $T_1$-invariant. Let us pick a coordinatization of the
directional CL-factor, say as $\eta:\bfX\to (Z,m_Z,\phi)\ltimes
(A,m_A,\s)$ for some compact metrizable Abelian groups $Z$ and $A$,
a dense homomorphism $\phi:\bbZ^2\to Z$ and a cocycle
$\s:\bbZ^2\times Z\to A$ over $R_\phi$, chosen so that the canonical
further factor onto $(Z,m_Z,\phi)$ is the whole Kronecker factor. By
Lemma~\ref{lem:make-Kron-DIO} we may assume that $(Z,m_Z,\phi)$ has
the DIO property.

In these terms, again by $L^2$-continuity and multilinearity, to
prove convergence of these averages it suffices to consider
functions $g_i(z,a)$ of the form $\k_i(z)\chi_i(a)$ with $\k_i\in
\widehat{Z}$ and $\chi_i \in \widehat{A}$ for $i=1,2$. We will refer
to functions of this form as \textbf{vertical eigenfunctions} of the
system $(Z,m_Z,\phi)\ltimes (A,m_A,\s)$, and will refer to the
characters $\chi_i$ appearing in their definition as their
associated \textbf{vertical characters}.  For these functions our
averages become
\begin{eqnarray*}
&&S'_N(g_1,g_2,f)(x)\\
&&= \frac{1}{N}\sum_{n=1}^N \k_1(\phi(\ell n^2 +
an,0)z)\cdot\chi_1(a)\cdot\chi_1(\s((\ell n^2 +
an,0),z))\\
&&\quad\quad\quad\quad\quad\quad\cdot\k_2(\phi(\ell n^2 +
an,n)z)\cdot\chi_2(a)\cdot\chi_2(\s((\ell n^2 + an,n),z))\cdot
f(T_2^n(x))\\
&&= \k_1(z)\chi_1(a)\k_2(z)\chi_2(a)\frac{1}{N}\sum_{n=1}^N
\k_1(\phi(\bf{e}_1))^{\ell n^2 + an}\k_2(\phi(\bf{e}_1))^{\ell
n^2 + an}\k_2(\phi(\bf{e}_2))^n\\
&&\quad\quad\quad\quad\quad\quad\quad\quad\cdot\chi_1(\s((\ell n^2 +
an,0),z))\cdot\chi_2(\s((\ell n^2 + an,n),z))\cdot f(T_2^n(x))
\end{eqnarray*}
where we write $(z,a) := \eta(x)$ and have used that $\k_i$ and
$\chi_i$ are characters.  Writing $\theta_1 := (\k_1\cdot
\k_2)(\phi(\bf{e}_1))$ and $\theta_2 := \k_2(\phi(\bf{e}_2))$, we
immediately deduce the following.

\begin{lem}\label{lem:further-reduced-polyconv}
The averages $S_N'(g_1,g_2,f)$ of
Proposition~\ref{prop:reduced-polyconv} all converge in $L^2(\mu)$
as $N\to\infty$ if and only if this is true of the averages
\[\frac{1}{N}\sum_{n=1}^N
\theta_1^{\ell n^2 + an}\theta_2^n\cdot\chi_1(\s((\ell n^2 +
an,0),z))\cdot\chi_2(\s((\ell n^2 + an,n),z))\cdot f(T_2^n(x))\] for
any $\theta_1,\theta_2 \in \Sone$. \qed
\end{lem}

In the conclusion of this lemma it is clear that the remaining
`awkwardness' for the purposes of proving norm convergence resides
in the expression
\[\chi_1(\s((\ell n^2 +
an,0),z))\cdot\chi_2(\s((\ell n^2 + an,n),z)).\] This is a sequence
of functions on the group rotation factor $Z$ whose behaviour as $n$
varies we have yet to control with much precision. Most of the
remainder of the proof will be directed towards exerting such
control.  In our approach to this we will follow the basic strategy
used by Host and Kra in~\cite{HosKra01} of arguing that if our
averages do not tend to $0$ in $L^2(\mu)$, then the cocycle $\s$
must give rise to some nontrivial Mackey data, and hence a
nontrivial combined cocycle equation, inside the Furstenberg
self-joining; and then using that equation itself to analyze the
behaviour of expressions such as our product of cocycles above.
However, the details of our implementation of this approach are
rather different from Host and Kra's, and in particular will rest on
much of our earlier study of directional CL-systems.

To begin the next stage of our analysis, we once again apply the van
der Corput estimate. Letting $u_n:= (g_1\circ T_1^{\ell n^2 +
an})(g_2\circ T_1^{\ell n^2 + an}T_2^n)(f\circ T_2^n)$, we deduce as
before that either $S'_N(g_1,g_2,f) \to 0$ in $L^2(\mu)$ or else we
also have
\begin{eqnarray*}
&&\frac{1}{H}\sum_{h=1}^H\frac{1}{N}\sum_{n=1}^N\int_X (g_1\circ
T_1^{\ell n^2 + 2\ell nh + \ell h^2 + an + ah})(\overline{g_1}\circ
T_1^{\ell n^2 + an})\\
&&\quad\quad\quad\quad\quad\quad\quad\quad\cdot(g_2\circ T_1^{\ell
n^2 + 2\ell nh + \ell h^2 + an + ah}T_2^{n+h})(\overline{g_2}\circ
T_1^{\ell n^2
+ an}T_2^n)\\
&&\quad\quad\quad\quad\quad\quad\quad\quad\cdot (f\circ
T_2^{n+h})(\overline{f}\circ T_2^n)\,\d\mu \not\to 0
\end{eqnarray*}
as $N\to\infty$ and then $H\to\infty$; and now, still as in the
previous section, using the $T_1$-invariance of $f$ we can change
variables in these integrals by $T_1^{-\ell n^2 - an}T_2^{-n}$ (and
change the order of some of the factors) to obtain
\begin{multline*}
\frac{1}{H}\sum_{h=1}^H\frac{1}{N}\sum_{n=1}^N\int_X (g_2\circ
T_1^{2\ell nh + \ell h^2 + ah}T_2^h)(g_1\circ T_1^{2\ell nh + \ell
h^2 + ah}T_2^{-n})(\overline{g_1}\circ
T_2^{-n})\\
\cdot \overline{g_2}\cdot(f\circ T_2^h)\cdot\overline{f}\,\d\mu
\not\to 0,
\end{multline*}
and this implies that for some $\eps > 0$ we have
\[\Big\|\lim_{N\to\infty}\frac{1}{N}\sum_{n=1}^N(g_2\circ T_1^{2\ell nh + \ell h^2 + ah}T_2^h)(g_1\circ
T_1^{2\ell nh + \ell h^2 + ah}T_2^{-n})(\overline{g_1}\circ
T_2^{-n})\Big\|_2^2 \geq \eps\] for infinitely many integers $h\geq
1$.

At this point another appeal to
Proposition~\ref{prop:Fberg-controls-nonconv-aves} implies that for
infinitely many integers $h\geq 1$ the function $(g_2\circ T_1^{\ell
h^2 + ah}T_2^h)\otimes(g_1\circ T_1^{\ell h^2 + ah})\otimes
\overline{g_1}$ has non-zero conditional expectation onto the
$\vec{T}_h$-invariant factor of $(X^3,\mu_h^{\rm{F}})$, where
$\vec{T}_h := T_1^{2\ell h}\times (T_1^{2\ell h}T_2^{-1})\times
T_2^{-1}$ as in that proposition. This is essentially the same
conclusion that was used for our first reduction above, except that
our change-of-variables above was slightly different this time
(there we changed by $T_1^{- \ell n^2 - an}$, rather than $T_1^{-
\ell n^2 - an}T_2^{-n}$), and this has led here to a different
triple of directions.

Nevertheless, they are still in general position with the origin,
and so we can make use of the description of the restriction of
$\mu^{\rm{F}}_h$ to $(Z\ltimes A)^3$ given in
Proposition~\ref{prop:Fberg-controls-nonconv-aves}. Observe also
that
\[g_2\circ T_1^{\ell h^2 + ah}T_2^h(z,a) = \k_2(\phi(\ell h^2 +
ah,h)z)\cdot \chi_2(\s((\ell h^2 + ah,h),z))\cdot\chi_2(a)\] is
still a vertical eigenfunction with vertical character $\chi_2$, and
similarly $g_1\circ T_1^{\ell h^2 + ah}$ and $\overline{g_1}$.
Combining this with the description of the $\vec{T}_h$-ergodic
components of $\mu^{\rm{F}}_h$ given in
Proposition~\ref{prop:Fberg-controls-nonconv-aves}, it follows that
if $(g_2\circ T_1^{\ell h^2 + ah}T_2^h)\otimes(g_1\circ T_1^{\ell
h^2 + ah})\otimes \overline{g_1}$ has nontrivial conditional
expectation onto the $\vec{T}_h$-invariant factor then the character
$\chi_2\otimes \chi_1\otimes \overline{\chi_1}$ must have nonzero
average over the Mackey group $M_h \leq A^3$. Combining this with
our other conclusions leads to the following.

\begin{lem}\label{lem:final-choice-of-hs}
For any $h$ for which the above averages do not tend to zero we must
have
\[M_h \leq \ker (\chi_2\otimes \chi_1\otimes \overline{\chi_1})\]
where $M_h$ is the Mackey group given by
Proposition~\ref{prop:Fberg-controls-nonconv-aves}, and so its
Mackey section quotients to give a Borel function $b_h:Z_h \to
\Sone$ such that
\begin{multline*}
\chi_2\circ\s((2\ell h,0),z_1)\cdot\chi_1\circ\s((2\ell h,-1),z_2)\cdot\overline{\chi_1}\circ\s((0,-1),z_3)\\
= \Delta_{(\phi(2\ell h\bf{e}_1),\phi(2\ell h\bf{e}_1 -
\bf{e}_2),\phi(-\bf{e}_2))}b_h(z_1,z_2,z_3)
\end{multline*}
for Haar-a.e. $(z_1,z_2,z_3)\in Z_h$. \qed
\end{lem}

We will soon argue that given any two different values of $h$, say
$h_1$ and $h_2$, for which the conclusion of
Lemma~\ref{lem:final-choice-of-hs} holds, we can use the structure
of directional CL-systems in conjunction with the above combined
coboundary equations to give some useful information for our
combined cocycle on a subgroup of $Z^3$ that is `effectively' much
larger than either of $Z_{h_1}$ or $Z_{h_2}$ individually, and for a
whole finite-index subgroup $\G\leq \bbZ^2$.

\subsection{Using several combined coboundary equations}

The following is another useful consequence of the DIO property.

\begin{lem}\label{lem:re-alignment-isos}
If $(Z,m_Z,\phi)$ has the DIO property and $\bf{n}_1$, $\bf{n}_2 \in
\bbZ^2$ are linearly independent then there is a unique continuous
isomorphism $\g_{\bf{n}_1,\bf{n}_2}:K_{\bf{n}_1}\to K_{\bf{n}_2}$
such that the map \[u\mapsto u\cdot\g_{\bf{n}_1,\bf{n}_2}(u)\] is an
isomorphism $K_{\bf{n}_1}\to K_{\bf{n}_1 + \bf{n}_2}$.
\end{lem}

\textbf{Proof}\quad Since $\bf{n}_1 = (\bf{n}_1 + \bf{n}_2) -
\bf{n}_2$ it follows that $K_{\bf{n}_1} \leq K_{\bf{n}_1 +
\bf{n}_2}\cdot K_{\bf{n}_2}$.  Hence for any $u \in K_{\bf{n}_1}$
there are $w\in K_{\bf{n}_1 + \bf{n}_2}$ and $v\in K_{\bf{n}_2}$
such that $u = wv^{-1}$, and moreover the DIO property implies that
$K_{\bf{n}_1 + \bf{n}_2}\cap K_{\bf{n}_2} = \{1_Z\}$ and so these
$w$ and $v$ are uniquely determined.  Now setting
$\g_{\bf{n}_1,\bf{n}_2}(u) := v$ it follows easily from uniqueness
that this is a continuous homomorphism, and that it has the
analogously-defined map $\g_{\bf{n}_2,\bf{n}_1}$ for an inverse and
so is an isomorphism.  Finally, we can check similarly that the map
\[u\mapsto u\cdot\g_{\bf{n}_1,\bf{n}_2}(u)\]
simply gives the analogously-defined map $\g_{\bf{n}_1,\bf{n}_1 +
\bf{n}_2}$ so it is also a continuous isomorphism.  This completes
the proof. \qed

We now introduce the `essentially larger' subgroup of $Z^3$ where we
will still be able to establish some useful structure to our
combined cocycle. Recalling that the target of $\eta$ is an
$(m\bbZ^2,m\bf{e}_1,m\bf{e}_2)$-directional CL-system for some
$m\geq 1$, and given two distinct integers $h_1$ and $h_2$
satisfying the conclusions of Lemma~\ref{lem:final-choice-of-hs},
let $h := 2\ell\cdot\rm{l.c.m.}(m,h_1,h_2,h_1 + h_2,h_1 - h_2)$, and
let
\[\t{Z}_0 := \{(zu,zuv,zv):\ z\in Z,\,u\in K_{h\bf{e}_1},\,v\in K_{h\bf{e}_2}\}.\]

It is easy to see that $\t{Z}_0\cap Z_{h_i}$ is always of finite
index in $Z_{h_i}$ for $i=1,2$: indeed, if $(zu,zuv,zv)\in Z_{h_i}$
then there is always some $k\in\{0,1,\ldots,h\}$ for which
$(zu,zuv,zv)\cdot(\phi(2\ell h_i\bf{e}_1),\phi(2\ell h_i \bf{e}_1 +
\bf{e}_2),\phi(\bf{e}_2))^k \in \t{Z}_0$.  On the other hand, this
intersection can be of infinite index in $\t{Z}_0$.

Let $\psi:\bbZ^2\to Z^3$ be the homomorphism $(n_1, n_2) \mapsto
(\phi(n_1\bf{e}_1),\phi(n_1\bf{e}_1 +
n_2\bf{e}_2),\phi(n_2\bf{e}_2))$. Also, by restricting from our
$\bbZ^2$-action to any of the (finitely many) ergodic components of
the subaction of $h\bbZ^2$, and observing that all of the structural
information we have accrued so far is preserved, we may assume that
the subaction of $h\bbZ^2$ is ergodic.

We will show that given the two combined coboundary equations from
Lemma~\ref{lem:final-choice-of-hs} for $h_1$ and $h_2$ and also the
previously-obtained structure of a directional CL-system, we can
actually obtain some useful information on the combined cocycle over
$R_\psi$ for the whole of the further finite-index subgroup $\G :=
\bbZ(2\ell h_1 h,h) + \bbZ(2\ell h_2 h,h)\leq h\bbZ^2$.

\begin{lem}\label{lem:big-combined-cobdry-eqn}
For any integers $h_1$, $h_2$ and $h$ satisfying the conclusion of
Lemma~\ref{lem:final-choice-of-hs} there are a Borel maps
$\t{b}_i:\t{Z}_0\to\Sone$ and $\t{c}_i:K_{h\bf{e}_1}\times
K_{h\bf{e}_2} \to \Sone$ for $i=1,2$ such that $\t{c}_i$ takes the
special form of the functions output by
Proposition~\ref{prop:combining-indiv-cobdry-eqns}, and
\begin{multline*}
\chi_2\circ\s(2\ell h_i h\bf{e}_1,zu)\cdot\chi_1\circ\s(2\ell
h_ih\bf{e}_1 -
h\bf{e}_2,zuv)\cdot\overline{\chi_1}\circ\s(-h\bf{e}_2,zu)\\ =
\Delta_{\psi(2\ell h_ih,h)}\t{b}_i(zu,zuv,zv)\cdot \t{c}_i(u,v)
\end{multline*}
for Haar-a.e. $(zu,zuv,zv)\in \t{Z}_0$.
\end{lem}

\textbf{Proof}\quad First note that \[R_{\psi(n_1,n_2)}(zu,zuv,zv) =
\big(z(u\phi(n_1\bf{e}_2)),z(u\phi(n_1\bf{e}_1))(v\phi(n_2\bf{e}_2)),z(v\phi(n_2\bf{e}_2))\big).\]
As a result, the above combined cocycle equation can be regarded
separately for each fixed value of $z$ as an equation involving only
the variables $u$ and $v$.  Therefore it suffices to prove instead
the existence of maps $\t{c}_i$ satisfying the above equations that
are simply Borel, $R_{(\phi(2\ell h_ih,0),\phi(0,-h))}$-invariant
and do not depend on $z$, since we can then choose some generic
$z\in Z$ and apply
Proposition~\ref{prop:combining-indiv-cobdry-eqns} to the resulting
combined cocycle equations for that fixed $z$ to modify each
$\t{c}_i$ into the desired special form.

Having observed this, the proof that there are Borel maps $\t{b}_i$
and $\t{c}_i$ of this form satisfying the above equation will not
involve the fact that we are assuming ourselves given two distinct
values of $h_i$ as output by be Lemma~\ref{lem:final-choice-of-hs};
the only appeal we make to this fact is in this initial application
of Proposition~\ref{prop:combining-indiv-cobdry-eqns}.

Let us write
\begin{multline*}
\tau_i(zu,zuv,zv)\\ := \chi_2\circ\s(2\ell
h_ih\bf{e}_1,zu)\cdot\chi_1\circ\s(2\ell h_ih\bf{e}_1 -
h\bf{e}_2,zuv)\cdot\overline{\chi_1}\circ\s(-h\bf{e}_2,zv).
\end{multline*}

We will need the isomorphisms given by
Lemma~\ref{lem:re-alignment-isos}. In particular, let
$\g_i:K_{(0,h)}\to K_{(2\ell hh_i,0)}$ be such that
$v\g_i(v)^{-1}\in K_{(2\ell hh_i,-h)}$ for all $v\in K_{(0,h)}$.

For any $(zu,zuv,zv)\in \t{Z}_0$ consider the decomposition
\begin{multline*}
\tau_i(zu,zuv,zv)\\ = \tau_i(zu,zuv,zu\g_i(v)^{-1}v) \cdot
\overline{\big(\overline{\tau_i(zu,zuv,zv)}\cdot\tau_i(zu,zuv,zu\g_i(v)^{-1}v)\big)}.
\end{multline*}
We will examine the two factors on the right-hand side of this
decomposition separately.

On the one hand, by the construction of $\g_i$ we know that
$(zu,zuv,zu\g_i(v)^{-1}v) \in Z_{h_i}$ and that the map $\t{Z}_0\to
Z_{h_i}:(zu,zuv,zv)\mapsto (zu,zuv,zu\g_i(v)^{-1}v)$ is a
homomorphism that covers a finite-index (and so positive-measure)
subgroup of $Z_{h_i}$, because by the uniqueness of $\g_i$ it must
be the identity on $\t{Z}_0\cap Z_{h_i}$. Hence by
Lemma~\ref{lem:final-choice-of-hs} we have
\[\tau_i(zu,zuv,zu\g_i(v)^{-1}v) = (\Delta_{\psi(2\ell h_ih,-h)}b_{h_i})(zu,zuv,zu\g_i(v)^{-1}v)\] for
$m_{\t{Z}_0}$-a.e. $(zu,zuv,zv)$. Since we must have
$\g_i(\phi(0,-h)) = \phi(2\ell h_ih,0)$, again by the uniqueness of
$\g_i$, and therefore $\phi(2\ell h_ih,0)\g_i(\phi(0,-h))^{-1} = 1$,
if we define
\[b'_i(zu,zuv,zv) := b_{h_i}(zu,zuv,zu\g_i(v)^{-1}v)\]
then it follows that
\begin{eqnarray*}
&&\Delta_{\psi(2\ell h_ih,-h)}b'_i(zu,zuv,zv)\\ &&=
b'_i(zu\cdot\phi(2\ell h_ih,0),zuv\cdot\phi(2\ell h_i
h,-h),zv\cdot\phi(0,-h))\cdot\overline{b'_i(zu,zuv,zv)}\\
&&= b_{h_i}(zu\cdot\phi(2\ell h_ih,0),zuv\cdot\phi(2\ell h_i
h,-h),zu\g_i(v)^{-1}v\cdot\phi(0,-h))\\
&&\quad\quad\quad\quad\quad\quad\quad\quad \cdot
\overline{b_{h_i}(zu,zuv,zu\g_i(v)^{-1}v)}\\
&&= (\Delta_{\psi(2\ell h_ih,-h)}b_{h_i})(zu,zuv,zu\g_i(v)^{-1}v),
\end{eqnarray*}
and so we can re-express the above coboundary equation as
\[\tau_i(zu,zuv,zu\g_i(v)^{-1}v) = \Delta_{\psi(2\ell
h_ih,-h)}b'_i(zu,zuv,zv).\]

On the other hand, recalling the consequences of the directional
CL-structure obtained in Corollary~\ref{cor:str-dCL-nice}, we know
that there are Borel maps $b_i^\circ:K_{(2\ell h_i,0)}\times Z\to
\Sone$ and $c_i^\circ:K_{(2\ell h_i,0)}\times
Z/\overline{\phi(h\bbZ^2)}\to \Sone$ such that
\begin{eqnarray*}
&&\overline{\tau_i(zu,zuv,zv)}\cdot \tau_i(zu,zuv,zu\g_i(v)^{-1}v)\\
&&=
\Delta_{u\g_i(v)^{-1}}\overline{\chi_1}\circ\s(-h\bf{e}_2,zv)\\
&&= b_i^\circ(u\g_i(v)^{-1},zv\cdot \phi(0,-h))\cdot
\overline{b_i^\circ(u\g_i(v)^{-1},zv)}\cdot
c_i^\circ(u\g_i(v)^{-1},zv\cdot \overline{\phi(h\bbZ^2)}).
\end{eqnarray*}
Moreover, recalling that we have reduced to the case in which
$h\bbZ^2$ acts ergodically through $R_\phi$, the dependence on the
coset $zv\cdot\overline{\phi(h\bbZ^2)}$ above may be dropped.

Since the map $(zu,zuv,zv)\mapsto (u\g_i(v)^{-1},zv)$ is also easily
seen to be a homomorphism onto a finite-index (and hence
positive-measure) subgroup of $K_{(2\ell h_i,0)}\times Z$, the above
holds $m_{\t{Z}_0}$-almost everywhere. In addition, if we now define
\[b''_i(zu,zuv,zv) := b^\circ_i(u\g_i(v)^{-1},zv)\]
then using again that fact that $\phi(2\ell
h_ih,0)\g(\phi(0,-h))^{-1} = 1$ we can compute directly that
\begin{eqnarray*}
&&\Delta_{\psi(2\ell h_ih,-h)}b''_i(zu,zuv,zv)\\ &&=
b''_i(zu\cdot\phi(2\ell h_ih,0),zuv\cdot\phi(2\ell h_i
h,-h),zv\cdot\phi(0,-h))\cdot\overline{b''_i(zu,zuv,zv)}\\
&&= b^\circ_i(u\g_i(v)^{-1},zv\cdot\phi(0,-h))\cdot
\overline{b^\circ_i(u\g_i(v)^{-1},zv)},
\end{eqnarray*}
and so we can re-express the above coboundary equation as
\begin{multline*}
\overline{\tau_i(zu,zuv,zv)}\cdot \tau_i(zu,zuv,zu\g_i(v)^{-1}v)\\
=\Delta_{\psi(2\ell h_ih,-h)}b''_i(zu,zuv,zv)\cdot
c^\circ_i(u\g_i(v)^{-1}).
\end{multline*}

Finally we can put the coboundary equations obtained above together
by setting $\t{b}_i:= b'_i\cdot\overline{b''_i}$ and
\[\t{c}_i(u,v) :=
\overline{c_i^\circ(u\g_i(v)^{-1})}\] to obtain
\begin{eqnarray*}
&&\tau_i(zu,zuv,zv)\\
&&= \tau_i(zu,zuv,zu\g_i(v)^{-1}v)\cdot
\overline{\big(\overline{\tau_i(zu,zuv,zv)}\cdot\tau_i(zu,zuv,zu\g_i(v)^{-1}v)\big)}\\
&&= \Delta_{\psi(2\ell h_1h,-h)}\t{b}_i(zu,zuv,zv)\cdot \t{c}_i(u,v)
\end{eqnarray*}
$m_{\t{Z}_0}$-almost everywhere, where $\t{c}_i(u,v)$ is
$R_{(\phi(2\ell h_1h,0),\phi(0,-h))}$-invariant, as required. \qed

The remaining steps in the proof of
Proposition~\ref{prop:reduced-polyconv} follow quite closely the
ideas of Host and Kra's neat approach in~\cite{HosKra01} to the
convergence of triple linear averages associated to three powers of
a single ergodic transformation.

The main technical result we need is the `compactification' result
for the family of functions
\[\chi_1\circ\s(n_1\bf{e}_1,\,\cdot\,)\cdot\chi_2\circ\s(n_1\bf{e}_1 +
n_2\bf{e}_2,\,\cdot\,)\cdot\overline{\chi_2}\circ\s(n_2\bf{e}_2,\,\cdot\,)\quad\quad
(n_1,n_2)\in \G\] given in the next proposition. This will serve as
our analog of Lemma 4.2 of~\cite{HosKra01}, but it differs from that
result in certain important details.  Most notably, our proposition
is a little more `quantitative', as a result of the introduction of
an additional `phase function' given by a generalized polynomial.
Generalized polynomials have been objects of interest among ergodic
theorists for some time, and so we recall their definition here for
completeness but will refer elsewhere for their properties that we
need.

\begin{dfn}[Gen-polynomials]
A map $p:\bbZ^2\to \bbR$ is a \textbf{generalized polynomial
(`gen-polynomial')} if it can be expressed using repeated
composition of ordinary real-valued polynomials and the operations
of taking the integer part, addition and multiplication.
\end{dfn}

For the basic properties of gen-polynomials we refer to Bergelson
and Leibman~\cite{BerLei07}, Leibman~\cite{Lei09} and the references
given there.  Recall that we have now restricted our attention to
the subgroup $\G := \bbZ(2\ell h_1 h,h) + \bbZ(2\ell h_2 h,h)$, and
let us henceforth write $\bf{q}_i = (q_{i1},q_{i2}) := (2\ell h_i
h,h)$ for brevity.  We also now abbreviate $K_1 := K_{(h,0)}$ and
$K_2 := K_{(0,h)}$, and observe from the DIO property that $K_1\cdot
K_2 \cong K_1\times K_2$ in $Z$, so in particular for any $\bf{n}\in
\G \leq h\bbZ^2$ we may interpret each $\phi(\bf{n})$ uniquely as a
member of $K_1\times K_2$.

\begin{prop}~\label{prop:compactification}
There is a gen-polynomial $p:\bbZ^2\to \bbR$ for which the following
holds.  For any $\a > 0$ there are
\begin{itemize}
\item a Borel function $C_\a:Z\times (K_1\times K_2)^2\to \Sone$ such that the
family of slices
\[Z\mapsto \Sone:z\mapsto C_\a(z,u_1,v_1,u_2,v_2)\]
indexed by $(u_1,v_1,u_2,v_2) \in (K_1\times K_2)^2$ all lie in
$L^2(m_Z)$ and vary continuously with $(u_1,v_1,u_2,v_2)$, and
\item an open subset $U_\a \subseteq K_1\times K_2$ of the form
\[U_\a = \bigcap_{\g \in \F}\{(u,v)\in K_1\times K_2:\ \delta < \{\g(u,v)\} < 1 - \delta\}\]
for some $\delta > 0$ and some finite subset $\F \subseteq
\widehat{K_1\times K_2}$ such that $\g(\phi(\bf{q}_1))\in\Sone$ is
irrational for every $\g\in \F$ and $m_{K_1\times K_2}(U_\a)
> 1 - \a$
\end{itemize}
such that
\begin{eqnarray*}
&&\exp(2\pi\rm{i}p(m,n))\cdot\chi_2\circ\s((m q_{11} +
nq_{21},0),z)\\
&&\quad\quad\quad\quad\quad\quad\cdot\chi_1\circ\s((mq_{11} + nq_{21},mq_{12} + nq_{22}),z)\cdot\overline{\chi_1}\circ\s((0,mq_{12} + nq_{22}),z)\\
&&\quad\quad\quad\quad\quad\quad\quad\quad=
C_\a(z,\phi(m\bf{q}_1),\phi(n\bf{q}_2))\quad\quad\hbox{for
Haar-a.e.}\ z \in Z
\end{eqnarray*}
for every $m \in \bbZ$ such that $\phi(m\bf{q}_1) \in U_\a$, where
we use our identification of $\phi(m\bf{q}_i) \in K_1\cdot K_2\leq
Z$ with a member of $K_1\times  K_2$.
\end{prop}

We will prove this proposition following a couple of preparatory
lemmas. The first of these is a simple calculation from
Lemma~\ref{lem:big-combined-cobdry-eqn}.

\begin{lem}\label{lem:finally-putting-the-combd-cocyc-together}
For any point $m\bf{q}_1 + n\bf{q}_2\in \G$ we have
\begin{eqnarray*}
&&\chi_2\circ\s((2\ell h_1 hm + 2\ell h_2hn,0),zu)\\
&&\quad\quad\cdot\chi_1\circ\s((2\ell
h_1 hm + 2\ell h_2 hn, - hm - hn),zuv)\\
&&\quad\quad\quad\quad\cdot\overline{\chi_1}\circ\s((0,- hm - hn),zu)\\
&&= \Delta_{\psi(2\ell h_1
hm,-hm)}\t{b}_1(zu,zuv,zv)\\
&&\quad\quad\quad\quad\cdot\Delta_{\psi(2\ell h_2
hn,-hn)}\t{b}_2(zu\phi(2\ell h_1hm,0),zuv\phi(2\ell
h_1hm,-hm),zv\phi(-hm))\\
&&\quad\quad\quad\quad\cdot\t{c}_1(u,v)^m\cdot \t{c}_2(u\phi(2\ell
h_1hm,0),v\phi(0,-hm))^n.
\end{eqnarray*}
\end{lem}

\textbf{Proof}\quad This follows immediately from the separate
conclusions of Lemma~\ref{lem:big-combined-cobdry-eqn} for $i=1$ and
$i=2$ by observing the consequences of the defining equations for a
cocycle over a $\bbZ^2$-action that
\begin{multline*}
\s((2\ell h_1 hm + 2\ell h_2hn,0),zu)\\ = \s((2\ell h_1
hm,0),zu)\cdot\s((2\ell h_2hn,0),zu\cdot\phi(2\ell h_1 hm,0)),
\end{multline*}
\begin{multline*}
\s((2\ell h_1 hm + 2\ell h_2hn,-hm -hn),zuv)\\ = \s((2\ell h_1
hm,-hm),zuv) \cdot\s((2\ell h_2hn,-hn),zuv\cdot\phi(2\ell h_1
hm,-hm)),
\end{multline*}
and
\[\s((-hm-hn),zv) =
\s((0,-hm),zv)\cdot\s((0,-hn),zv\cdot\phi(0,-hm)),\] and then
multiplying these together. \qed

The proof of Proposition~\ref{prop:compactification} will also
require the following analog of an enabling lemma from Host and
Kra~\cite{HosKra01}.

\begin{lem}[C.f. Lemma 3.3
in~\cite{HosKra01}]\label{lem:HK-factorizing-almostzero} Suppose
that $\delta < 1/100$ and that $f_i:Z\to\Sone$, $i=1,2,3$, and
$h:K_1\times K_2 \to \Sone$ are Borel functions such that
\[f_1(zu)f_2(zuv)f_3(zv)h(u,v) \approx_{\delta} 1\quad\quad\hbox{in}\ L^2(m_{Z\times K_1\times K_2})\]
Then there are Borel functions $g_1:Z/K_1\to\Sone$ and $g_2:Z/K_2\to
\Sone$, characters $\g_1\in K_2^\perp$ and $\g_2 \in K_1^\perp$ and
constants $\a_1,\a_2,\a_3\in \Sone$ and $\b \in \Sone$ satisfying
$\a_1\a_2\a_3\b = 1$ such that
\[f_1(z) \approx_{\rm{o}_\delta(1)} \a_1g_1(z K_1)\g_2(z),\]
\[f_2(z) \approx_{\rm{o}_\delta(1)} \a_2\overline{g_1(z K_1)\g_1(z)g_2(z K_2)\g_2(z)},\]
\[f_3(z) \approx_{\rm{o}_\delta(1)} \a_3g_2(z K_2)\g_1(z)\]
and
\[h(u,v) \approx_{\rm{o}_\delta(1)} \b\g_1(u)\g_2(v),\]
where all approximations hold in the norm of the relevant $L^2$
space.  Consequently we also have
\begin{multline*}
f_1(z)f_2(z)f_3(z)\\ \approx_{\rm{o}_\delta(1)} \a_1g_1(z
K_1)\chi_2(z)\cdot \a_2\overline{g_1(z K_1)\g_1(z)g_2(z
K_2)\g_2(z)}\cdot \a_3g_2(z K_2)\g_2(z) \equiv \b
\end{multline*}
in $L^2(m_Z)$.
\end{lem}

\textbf{Proof}\quad Recalling that the system $(Z,m_Z,\phi)$ is DIO
and that $K_1 \cdot K_2$ has finite index in $Z$, by restriction to
a coset we may assume that $Z = K_1\times K_2$, and so write the
given equation as
\[f_1(z_1u,z_2)f_2(z_1u,z_2v)f_3(z_1,z_2v)h(u,v) \approx_{\delta} 1\quad\quad\hbox{in}\ L^2(m_{K_1\times K_2\times K_1\times K_2}).\]

In the argument below all approximations $\approx$ will implicitly
refer to an error of the form $\rm{o}_\delta(1)$.

Changing variables so that $z'_1 := z_1u$ and $v' := z_2v$, this
becomes
\[f_1(z'_1,z_2)f_2(z'_1,v')f_3(z'_1u^{-1},v')h(u,v'z_2^{-1}) \approx_{\delta} 1\quad\quad\hbox{in}\ L^2(m_{K_1\times K_2\times K_1\times K_2}),\]
and so for most fixed choices of $u$ and $v'$ we have
\[f_1(z'_1,z_2)\approx_\delta \overline{f_2(z'_1,v')f_3(z'_1u^{-1},v')h(u,v'z_2^{-1})}\quad\quad\hbox{in}\ L^2(m_{K_1\times K_2}),\]
which is manifestly a product of functions each of which depends
only on $z_1'$ (or, equivalently, on $z_1$) or only on $z_2$. We may
therefore approximate
\[f_1(z_1,z_2)\approx g_{11}(z_1)g_{12}(z_2)\]
for some $g_{1i}:K_i\to\Sone$, and exactly similarly we can
approximate
\[f_3(z_1,z_2)\approx g_{31}(z_1)g_{32}(z_2).\]
Substituting these right-hand sides into our original approximation
we obtain
\[g_{11}(z_1u)g_{13}(z_2)f_2(z_1u,z_2v)g_{31}(z_1)g_{32}(z_2v)h(u,v) \approx 1\quad\quad\hbox{in}\ L^2(m_{K_1\times K_2\times K_1\times K_2}),\]
or, changing variables to $z'_1 := z_1u$ and $z'_2 := z_2v$,
\[g_{11}(z'_1)g_{13}(z'_2v^{-1})f_2(z'_1,z'_2)g_{31}(z_1u^{-1})g_{32}(z'_2)h(u,v) \approx 1\quad\quad\hbox{in}\ L^2(m_{K_1\times K_2\times K_1\times K_2}).\]

Again fixing some $u$ and $v$ for which this is true for most $z_1'$
and $z_2'$, we find that $f_2$ must also take an approximate product
form,
\[f_2(z_1,z_2) \approx g_{21}(z_1)g_{22}(z_2),\]
while fixing instead $z'_1$ and $z'_2$ and allowing $u$ and $v$ to
vary we obtain the same conclusion for $h$:
\[h(u,v)\approx h_1(u)h_2(v).\]
Now we substitute all these approximate factorizations back into our
original approximation one last time to obtain
\begin{multline*}
\big((g_{11}\cdot g_{21})(z_1u)\cdot g_{31}(z_1)\cdot
h_1(u)\big)\cdot \big(g_{12}(z_2)\cdot (g_{22}\cdot
g_{32})(z_2v)\cdot h_2(v)\big)\approx 1\\ \hbox{in}\
L^2(m_{K_1\times K_2\times K_1\times K_2}),
\end{multline*}
and so in fact we must have that $\big((g_{11}\cdot
g_{21})(z_1u)\cdot g_{31}(z_1)\cdot h_1(u)\big)$ is close to a
constant-valued map in $L^2(m_{K_1\times K_1})$ and similarly that
$\big(g_{12}(z_2)\cdot (g_{22}\cdot g_{32})(z_2v)\cdot h_2(v)\big)$
is close in $L^2(m_{K_2\times K_2})$ to a map with value the inverse
of that constant.

Calling this constant $\g \in \Sone$ and writing $h_1':= \g\cdot
h_1$, we are left with the approximate equation
\[(g_{11}\cdot
g_{21})(z_1u)\cdot g_{31}(z_1)\cdot h'_1(u) \approx 1\] in
$L^2(m_{K_1\times K_1})$.  Since the functions $g_{11}\cdot g_{21}$,
$g_{31}$ and $h_1'$ take values in $\Sone$, they all have norm $1$
in $L^2(m_{K_1})$.  On the other hand, averaging over $z_1$ in the
above approximation gives that
\[\overline{h_1'} \approx (g_{11}\cdot g_{31})\ast g'_{31}\]
where we define $g'_{31}(z_1) := \overline{g_{31}(z_1^{-1})}$.
Hence, taking the Fourier transform of this approximation gives
\[\|\overline{h_1'} - (g_{11}\cdot g_{31})\ast g'_{31}\|_2^2 = \sum_{\g \in \widehat{K_1}}|\widehat{\overline{h_1'}}(\g) - \widehat{(g_{11}\cdot g_{31})}(\g)\cdot\widehat{g'_{31}}(\g)|^2 \approx 0\]
and hence also
\[1 = \|h_1'\|_2^2 = \sum_{\g \in \widehat{K_1}}|\widehat{h_1'}(\g)|^2 \approx \sum_{\g \in \widehat{K_1}}|\widehat{(g_{11}\cdot g_{31})}(\g)\cdot\widehat{g'_{31}}(\g)|^2.\]
On the other hand, simply by the non-negativity of all the terms
involved we have \[\sum_{\g \in \widehat{K_1}}|\widehat{(g_{11}\cdot
g_{31})}(\g)|^2|\widehat{g'_{31}}(\g)|^2 \leq \Big(\sum_{\g \in
\widehat{K_1}}|\widehat{(g_{11}\cdot
g_{31})}(\g)|^2\Big)\Big(\sum_{\g \in
\widehat{K_1}}|\widehat{g'_{31}}(\g)|^2\Big) \leq 1\] with
approximate equality only if $\widehat{g'_{31}}$ and
$\widehat{(g_{11}\cdot g_{31})}$ are both concentrated on a single
character.

Thus the above approximation in $L^2(m_{K_1\times K_1})$ is possible
only if there are some character $\g_1 \in \widehat{K_1}$ and some
constants $\a_3,\eta_1 \in \Sone$ such that $g_{31} \approx
\a_3\g_1$, $g_{11}\cdot g_{31}\approx \eta_1\overline{\g_1}$ and
$h'_1 \approx \overline{\a_3\eta_1}\g_1$.  Exactly similarly we
obtain a character $\g_2 \in \widehat{K_2}$ and constants
$\a_1,\eta_2 \in \Sone$ such that $g_{12} \approx \a_1\g_2$,
$g_{22}\cdot g_{32}\approx \eta_2\overline{\g_2}$ and $h'_2 \approx
\overline{\a_1\eta_2}\g_2$.  Setting $\a_2 :=
\overline{\eta_1\eta_2}$ and $\b :=
\overline{\a_1\eta_2}\cdot\overline{\a_3\eta_1} =
\overline{\a_1\a_2\a_3}$, we see that combining these resulting
approximants gives the result. The final assertion that
\begin{multline*}
f_1(z)f_2(z)f_3(z)\\ \approx_{\rm{o}_\delta(1)} \a_1g_1(z
K_1)\g_2(z)\cdot \a_2\overline{g_1(z K_1)\g_1(z)g_2(z
K_2)\g_2(z)}\cdot \a_3g_2(z K_2)\g_2(z) \equiv \b
\end{multline*}
in $L^2(m_Z)$ follows immediately.  \qed

\textbf{Proof of Proposition~\ref{prop:compactification}}\quad This
will rest on the special form of the functions $\t{c}_i$ obtained
from Proposition~\ref{prop:combining-indiv-cobdry-eqns} and its
consequence Lemma~\ref{lem:big-combined-cobdry-eqn}.  Those results
tell us that these functions are of the form
\[\t{c}_i(u,v) = \a_i(u,v)\exp\Big(2\pi\rm{i}\sum_{j=1}^{J_i}a_{i,j}(u,v)\{\chi_{i,j}(\phi(2\ell h_i h\bf{e}_1),\phi(-h\bf{e}_2))\}\{\g_{i,j}(u,v)\}\Big)\]
for some maps $\a_1,\a_2:K_1\times K_2 \in \Sone$ and
$a_{i,j}:K_1\times K_2\to \bbZ$ that factorize through some finite
quotient group and some characters
$\g_{i,j},\chi_{i,j}\in\widehat{K_1\times K_2}$.

In this expression, we note that if for some $j$ the character
$\g_{i,j}$ has image a finite subgroup of $\Sone$, rather than the
whole of $\Sone$, then we can simply replace $\a_i(u,v)$ by
\[\a_i(u,v)\cdot \exp(2\pi\rm{i}a_{i,j}(u,v)\{\chi_{i,j}(\phi(2\ell h_i
h\bf{e}_1),\phi(h\bf{e}_2))\}\{\g_{i,j}(u,v)\})\] and remove the
term $a_{i,j}(u,v)\{\chi_{i,j}(\phi(2\ell h_i
r\bf{e}_1),\phi(r\bf{e}_2))\}\{\g_{i,j}(u,v)\}$ from the sum inside
the main exponential. Therefore we may assume further that in this
expression the characters $\g_{i,j}$ all map $K_1\times K_2$ onto
the whole of $\Sone$.  Having made these arrangements, we may now
choose some large integer $r\geq 1$ for which each $\a_i$ and
$a{i,j}$ is actually constant on each coset of
$\overline{\psi(r\G)}$. Replacing $h$ with $rh$, each $\bf{q}_i$
with $r\bf{q}_i$, and thus $\G$ with the further finite-index
sublattice $r\G \leq \G$, we may now simply assume that each $\a_i$
and $a_{i,j}$ is constant.

Now let $\cal{J}\subseteq \{1,2,\ldots,J_2\}$ be the subset of
indices for which $\g_{2,j}(\bf{q}_1)$ is an irrational element of
the circle group $\Sone$. From the condition that each $\g_{i,j}$
have range equal to the whole of $\Sone$ it follows that for any $\a
> 0$ there is some $\delta(\a)
> 0$ such that the open set
\[U_\a:= \{(u,v)\in K_1\times K_2:\ \delta(\a) < \{\g_{2,j}(u,v)\} < 1 - \delta(\a)\ \forall j \in \cal{J}\}\]
has $m_{K_1\times K_2}(U_\a) > 1 - \a$.  In addition, we may take
$\a \mapsto \delta(\a)$ to be strictly increasing for sufficiently
small $\a$, so that $\overline{U_\a} \subseteq U_{\a/2}$. We will
obtain the function $C_\a$ by showing that for a suitably chosen
generalized polynomial $p$, for any sequence $(m_k,n_k)_{k \geq 1}$
in $\bbZ^2$ such that
\[\phi(m_k\bf{q}_1) \in U_{\a/2}\quad\quad\forall k,\]
\[\phi(m_k\bf{q}_1)\to (u^\circ_1,v^\circ_1) \in K_1\times K_2\quad\quad\hbox{as}\ k\to\infty\]
and \[\phi(n_k\bf{q}_2)\to (u^\circ_2,v^\circ_2) \in K_1\times
K_2\quad\quad\hbox{as}\ k\to\infty,\] we have that the sequence of
functions
\begin{multline*}
z\mapsto \exp(2\pi\rm{i}p(m_k,n_k))\cdot\chi_2\circ \s(((m_kq_{11} +
n_kq_{21},0),z)\\ \cdot \chi_1\circ \s((m_kq_{11} +
n_kq_{21},m_kq_{12} + n_kq_{22}),z)\cdot \overline{\chi_1}\circ
\s((0,m_kq_{21} + n_kq_{22}),z)
\end{multline*}
on $Z$ converges in $L^2(m_Z)$.  From this it follows that for any
$(u_1,v_1)\in U_{\a/2}$ we may unambiguously define a function
$z\mapsto C'_\a(z,u_1,v_1,u_2,v_2)$ to be the limit of these
functions when $u_i = u^\circ_i$ and $v_i = v^\circ_i$, and this
defines a Borel map $C'_\a$ on $Z\times U_{\a/2}\times (K_1\times
K_2)$ such that $(u_1,v_1,u_2,v_2)\mapsto
C'_\a(\,\cdot\,,u_1,v_1,u_2,v_2)$ is a continuous map from
$U_{\a/2}\times (K_1\times K_2)$ to $L^2(m_Z)$. Having done this we
can simply choose any continuous function $\varphi$ satisfying
$1_{U_\a} \leq \varphi \leq 1_{U_{\a/2}}$ and define
\[C_\a(z,u_1,v_1,u_2,v_2) :=
\left\{\begin{array}{ll}\varphi(u_1,v_1)C'_\a(z,u_1,v_1,u_2,v_2)&\quad\hbox{if}\
(u_1,v_1)\in U_{\a/2}\\
0&\quad\hbox{else:}\end{array}\right.\] it is now clear that this
function has the desired properties in conjunction with the set
$U_\a$.

Thus it remains to show this convergence for an arbitrary such
sequence $(m_k,n_k)$. Letting
\begin{eqnarray*}
f_{1,k}(z) &:=& \chi_2\circ\s((m_kq_{11} + n_kq_{21},0),z),\\
f_{2,k}(z) &:=& \chi_1\circ \s((m_kq_{11} + n_kq_{21},m_kq_{12} +
n_kq_{22}),z)\\
\hbox{and}\ f_{3,k} &:=& \overline{\chi_1}\circ\s((0,m_kq_{21} +
n_kq_{22}),z),
\end{eqnarray*}
from Lemma~\ref{lem:finally-putting-the-combd-cocyc-together} we
have
\begin{eqnarray*}
&&f_{1,k}(zu)f_{2,k}(zuv)f_{2,k}(zv)\\
&&= \Delta_{\psi(2\ell h_1
hm_k,-hm_k)}\t{b}_1(zu,zuv,zv)\\
&&\quad\quad\quad\quad\cdot\Delta_{\psi(2\ell h_2
hn_k,-hn_k)}\t{b}_2(zu\phi(2\ell h_1hm_k,0),zuv\phi(2\ell
h_1hm_k,-hm_k),zv\phi(-hm_k))\\
&&\quad\quad\quad\quad\cdot\t{c}_1(u,v)^{m_k}\cdot
\t{c}_2(u\phi(2\ell h_1hm_k,0),v\phi(0,-hm_k))^{n_k}.
\end{eqnarray*}
Re-arranging, we deduce that
\begin{eqnarray*}
&&f_{1,k}(zu)f_{2,k}(zuv)f_{2,k}(zv)\cdot\big(\t{c}_1(u,v)^{m_k}\cdot
\t{c}_2(u\phi(2\ell h_1hm_k,0),v\phi(0,-hm_k))^{n_k}\big)\\
&&= \Delta_{\psi(2\ell h_1
hm_k,-hm_k)}\t{b}_1(zu,zuv,zv)\\
&&\quad\quad\quad\quad\cdot\Delta_{\psi(2\ell h_2
hn_k,-hn_k)}\t{b}_2(zu\phi(2\ell h_1hm_k,0),zuv\phi(2\ell
h_1hm_k,-hm_k),zv\phi(-hm_k))\\
&&\to
\Delta_{(u^\circ_1,u^\circ_1v^\circ_1,v^\circ_1)}\t{b}_1(zu,zuv,zv)\cdot\Delta_{(u^\circ_2,u^\circ_2v^\circ_2,v^\circ_2)}\t{b}_2(zuu_{11},zuu^\circ_1vv^\circ_1,zvv^\circ_1)
\end{eqnarray*}
in $L^2(m_{Z_0})$ as $k\to\infty$, and hence that
\begin{eqnarray*}
&&f_{1,k}(zu)\overline{f_{1,\ell}(zu)}f_{2,k}(zuv)\overline{f_{2,\ell}(zuv)}f_{3,k}(zv)\overline{f_{3,\ell}(zv)}\\
&&\quad\quad\quad\quad \cdot \t{c}_1(u,v)^{m_k}\cdot
\t{c}_2(u\phi(2\ell h_1hm_k,0),v\phi(0,-hm_k))^{n_k}\\
&&\quad\quad\quad\quad \cdot\overline{\t{c}_1(u,v)^{m_\ell}\cdot
\t{c}_2(u\phi(2\ell h_1hm_\ell,0),v\phi(0,-hm_\ell))^{n_\ell}}\\
&&\to 0
\end{eqnarray*}
in $L^2(m_{Z_0})$ as $k,\ell\to \infty$.

It now follows from Lemma~\ref{lem:HK-factorizing-almostzero} that
in $L^2(m_Z)$ the $\Sone$-valued function
\[f_{1,k}(z)\overline{f_{1,\ell}(z)}f_{2,k}(z)\overline{f_{2,\ell}(z)}f_{3,k}(z)\overline{f_{3,\ell}(z)}\]
approaches the subset of constant $\Sone$-valued functions in
$L^2(m_Z)$ as $k,\ell \to\infty$, and that $\b_{k,\ell}\in\Sone$ is
a family of constants to which the above functions are
asymptotically equal if and only if the function
\begin{multline*}
\b_{k,\ell}\cdot \t{c}_1(u,v)^{m_k}\cdot
\t{c}_2(u\phi(2\ell h_1hm_k,0),v\phi(0,-hm_k))^{n_k}\\
\cdot\overline{\t{c}_1(u,v)^{m_\ell}\cdot \t{c}_2(u\phi(2\ell
h_1hm_\ell,0),v\phi(0,-hm_\ell))^{n_\ell}}
\end{multline*}
is close in $L^2(m_{K_1\times K_2})$ to a character (which is
necessarily unique once this approximation is sufficiently good,
since all distinct characters are separated by a distance of
$\sqrt{2}$ in $L^2(m_{K_1\times K_2})$).

To complete the proof, it will therefore suffice to find some
gen-polynomial $p(m,n)$ (not depending on the choices we made above
for a particular $\a$) such that the constants $\b_{k,\ell} =
\exp(2\pi\rm{i}(p(m_k,n_k) - p(m_\ell,n_\ell)))$ satisfy this latter
condition. We will now see that such a gen-polynomial can simply be
read off from the special form of the functions $\t{c}_1$ and
$\t{c}_2$ guaranteed by Lemma~\ref{lem:big-combined-cobdry-eqn} and
recalled above.

Indeed, having replaced $\G$ with the sufficiently small
finite-index subgroup $r\G$ and re-assigned our notation, these
functions are of the form
\[\t{c}_i(u,v) = \a_i\exp\Big(2\pi\rm{i}\sum_{j=1}^{J_i}a_{i,j}\{\chi_{i,j}(\phi(2\ell h_i h\bf{e}_1),\phi(-h\bf{e}_2))\}\{\g_{i,j}(u,v)\}\Big)\]
and for some $\a_1,\a_2 \in \Sone$, $a_{i,j}\in\bbZ$ and characters
$\g_{i,j},\chi_{i,j}\in\widehat{K_1\times K_2}$ whose images are the
whole circle group $\Sone$. In terms of these expressions we can now
write
\begin{eqnarray*}
&&\t{c}_1(u,v)^{m_k}\cdot \t{c}_2(u\phi(2\ell
h_1hm_k,0),v\phi(0,-hm_k))^{n_k}\\
&&=
\a_1^{m_k}\exp\Big(2\pi\rm{i}\sum_{j=1}^{J_1}m_ka_{1,j}\{\chi_{1,j}(\phi(2\ell h_1 h\bf{e}_1),\phi(-h\bf{e}_2))\}\{\g_{1,j}(u,v)\}\Big)\\
&&\quad\cdot
\a_2^{n_k}\exp\Big(2\pi\rm{i}\sum_{j=1}^{J_2}n_ka_{2,j}\{\chi_{2,j}(\phi(2\ell
h_2 h\bf{e}_1),\phi(-h\bf{e}_2))\}\\
&&\quad\quad\quad\quad\quad\quad\quad\quad\quad\quad\quad\quad\quad\quad\cdot\{\g_{2,j}(u\phi(2\ell
h_1 hm_k\bf{e}_1),v\phi(-hm_k\bf{e}_2))\}\Big).
\end{eqnarray*}
In order to use this expression we next note the elementary identity
\begin{eqnarray*}
&&\{\g_{2,j}(u\phi(2\ell h_1 hm_k\bf{e}_1),v\phi(-hm_k\bf{e}_2))\}\\
&&=
\{\g_{2,j}(u,v)\} + \{\g_{2,j}(\phi(2\ell h_1 hm_k\bf{e}_1),\phi(-hm_k\bf{e}_2))\}\\
&&\quad - \lfr\{\g_{2,j}(u,v)\} + \{\g_{2,j}(\phi(2\ell h_1
hm_k\bf{e}_1),\phi(-hm_k\bf{e}_2))\}\rfr.
\end{eqnarray*}
Substituting this identity and its partner for $(m_\ell,n_\ell)$ and
taking the difference of the results we obtain
\begin{eqnarray*}
&&\t{c}_1(u,v)^{m_k}\cdot \t{c}_2(u\phi(2\ell
h_1hm_k,0),v\phi(0,-hm_k))^{n_k}\\
&&\quad\quad\quad\quad\quad\quad\quad\quad\quad\quad\quad\quad\cdot\overline{\t{c}_1(u,v)^{m_\ell}\cdot
\t{c}_2(u\phi(2\ell
h_1hm_\ell,0),v\phi(0,-hm_\ell))^{n_\ell}} \\
&&=
\a_1^{m_k - m_\ell}\exp\Big(2\pi\rm{i}\sum_{j=1}^{J_1}a_{1,j}(m_k - m_\ell)\{\chi_{1,j}(\phi(2\ell h_1h\bf{e}_1),\phi(-h\bf{e}_2))\}\{\g_{1,j}(u,v)\}\Big)\\
&&\quad\cdot \a_2^{n_k -
n_\ell}\exp\Big(2\pi\rm{i}\sum_{j=1}^{J_2}a_{2,j}(n_k -
n_\ell)\{\chi_{2,j}(\phi(2\ell h_2
h\bf{e}_1),\phi(-h\bf{e}_2))\}\{\g_{2,j}(u,v)\}\Big)\\
&&\quad\cdot
\exp\Big(2\pi\rm{i}\sum_{j=1}^{J_2}a_{2,j}\{\chi_{2,j}(\phi(2\ell
h_2 h\bf{e}_1),\phi(-h\bf{e}_2))\}\big(n_k\{\g_{2,j}(\phi(2\ell h_1
hm_k\bf{e}_1),\phi(-hm_k\bf{e}_2))\}\\
&&\quad\quad\quad\quad\quad\quad\quad\quad\quad\quad\quad\quad -
n_\ell\{\g_{2,j}(\phi(2\ell h_1
hm_\ell\bf{e}_1),\phi(-hm_\ell\bf{e}_2))\}\big)\Big)\\
&&\quad\cdot\exp\Big(-2\pi\rm{i}\sum_{j=1}^{J_2}a_{2,j}\{\chi_{2,j}(\phi(2\ell
h_2 h\bf{e}_1),\phi(-h\bf{e}_2))\}\\
&&\quad\quad\quad\quad\quad\quad\quad\quad\cdot\big(n_k\lfr\{\g_{2,j}(u,v)\}
+
\{\g_{2,j}(\phi(2\ell h_1 hm_k\bf{e}_1),\phi(-hm_k\bf{e}_2))\}\rfr\\
&&\quad\quad\quad\quad\quad\quad\quad\quad\quad\quad -
n_\ell\lfr\{\g_{2,j}(u,v)\} + \{\g_{2,j}(\phi(2\ell h_1
hm_\ell\bf{e}_1),\phi(-hm_\ell\bf{e}_2))\}\rfr\big)\Big).
\end{eqnarray*}

Let us now consider some of the factors in this product in turn.

\begin{itemize}
\item First, we have by assumption that $\phi(2\ell h_1hm_k\bf{e}_1) \to
u^\circ_1$ and $\phi(-hm_k\bf{e}_2) \to v^\circ_1$ as $k\to\infty$.
Since $\chi_{1,j}$ is a character on $K_1\times K_2$, it follows
that
\[\rm{dist}\big((m_k - m_\ell)\{\chi_{1,j}(\phi(2\ell h_1h\bf{e}_1),\phi(-h\bf{e}_2))\}\,,\,\bbZ\big)\to 0\]
as $k,\ell \to\infty$.  Let us here write $I(r)\in \bbZ$ for the
closest integer to any $r\in\bbR$, rounding down when $r$ is a
proper half-integer, so that $I(r) \in \{\lfr r\rfr,\lfr r\rfr
+1\}$. From the above it follows that as $k,\ell \to\infty$ the
distance in $L^2(m_{K_1\times K_2})$ between the function
\[(u,v)\mapsto \exp\Big(2\pi\rm{i}\sum_{j=1}^{J_1}a_{1,j}(m_k -
m_\ell)\{\chi_{1,j}(\phi(2\ell
h_1h\bf{e}_1),\phi(-h\bf{e}_2))\}\{\g_{1,j}(u,v)\}\Big)\] and the
character
\begin{multline*}
\exp\Big(2\pi\rm{i}\sum_{j=1}^{J_1}a_{1,j}I\big((m_k -
m_\ell)\{\chi_{1,j}(\phi(2\ell
h_1h\bf{e}_1),\phi(-h\bf{e}_2))\}\big)\{\g_{1,j}(u,v)\}\Big)\\
= \prod_{j=1}^{J_1}\g_{1,j}(u,v)^{a_{1,j}I\big((m_k -
m_\ell)\{\chi_{1,j}(\phi(2\ell
h_1h\bf{e}_1),\phi(-h\bf{e}_2))\}\big)}
\end{multline*}
tends to $0$.  Exactly similarly the functions
\[\exp\Big(2\pi\rm{i}\sum_{j=1}^{J_2}a_{2,j}(n_k -
n_\ell)\{\chi_{2,j}(\phi(2\ell h_2
h\bf{e}_1),\phi(-h\bf{e}_2)\}\{\g_{2,j}(u,v)\}\Big)\] are also
asymptotically close to characters as $k,\ell\to\infty$, and hence
the same is true of the product of these two exponential functions.

\item Now consider the last factor above,
\begin{eqnarray*}
&&\exp\Big(-2\pi\rm{i}\sum_{j=1}^{J_2}a_{2,j}\{\chi_{2,j}(\phi(2\ell
h_2
h\bf{e}_1),\phi(-h\bf{e}_2)\}\\
&&\quad\quad\quad\quad\quad\quad\cdot\big(n_k\lfr\{\g_{2,j}(u,v)\} +
\{\g_{2,j}(\phi(2\ell h_1 hm_k\bf{e}_1),\phi(-hm_k\bf{e}_2))\}\rfr\\
&&\quad\quad\quad\quad\quad\quad\quad\quad\quad\quad -
n_\ell\lfr\{\g_{2,j}(u,v)\} + \{\g_{2,j}(\phi(2\ell h_1
hm_\ell\bf{e}_1),\phi(-hm_\ell\bf{e}_2))\}\rfr\big)\Big)\\
&&=\prod_{j=1}^{J_2}\exp\Big(-2\pi\rm{i}a_{2,j}\{\chi_{2,j}(\phi(2\ell
h_2 h\bf{e}_1),\phi(-h\bf{e}_2)\}\\
&&\quad\quad\quad\quad\quad\quad\cdot\big(n_k\lfr\{\g_{2,j}(u,v)\} +
\{\g_{2,j}(\phi(2\ell h_1 hm_k\bf{e}_1),\phi(-hm_k\bf{e}_2))\}\rfr\\
&&\quad\quad\quad\quad\quad\quad\quad\quad\quad\quad -
n_\ell\lfr\{\g_{2,j}(u,v)\} + \{\g_{2,j}(\phi(2\ell h_1
hm_\ell\bf{e}_1),\phi(-hm_\ell\bf{e}_2))\}\rfr\big)\Big)
\end{eqnarray*}
We will argue that each of the individual factors of this product
over $j$ is asymptotically close to the constant function $1$ in
$L^2(m_{K_1\times K_2})$, using again the fact that
\[\g_{2,j}(\phi(2\ell h_1 hm_k\bf{e}_1),\phi(-hm_k\bf{e}_2))\,,\,\g_{2,j}(\phi(2\ell h_1
hm_\ell\bf{e}_1),\phi(-hm_\ell\bf{e}_2))\to
\g_{2,j}(u^\circ_1,v^\circ_1)\] as $k,\ell \to\infty$.  For this
argument we must treat the cases $j \in \cal{J}$ and
$j\not\in\cal{J}$ separately.

If $j\in \cal{J}$, then we know that $\delta(\a/2) \leq
\{\g_{2,j}(u^\circ_1,v^\circ_1)\} \leq 1 - \delta(\a/2)$ from the
restriction $\phi(m_k\bf{q}_1) \in U_{\a/2}$ and continuity. This
implies that once $k$ and $\ell$ are sufficiently large then we have
that
\[\{\g_{2,j}(\phi(2\ell h_1 hm_k\bf{e}_1),\phi(-hm_k\bf{e}_2))\}\quad\hbox{and}\quad\{\g_{2,j}(\phi(2\ell h_1
hm_\ell\bf{e}_1),\phi(-hm_\ell\bf{e}_2))\}\] lie close together and
both inside $(0,1)$.  From this we deduce that
\begin{eqnarray*}
&&m_{K_1\times K_2}\big\{(u,v)\in K_1\times K_2:\
\lfr\{\g_{2,j}(u,v)\} + \{\g_{2,j}(\phi(2\ell h_1
hm_k\bf{e}_1),\phi(-hm_k\bf{e}_2))\}\rfr\\
&&\quad\quad\quad\quad\quad\quad\quad\quad\quad\quad\quad\neq
\lfr\{\g_{2,j}(u,v)\} + \{\g_{2,j}(\phi(2\ell h_1
hm_\ell\bf{e}_1),\phi(-hm_\ell\bf{e}_2))\}\rfr\big\}\\
&& \to 0
\end{eqnarray*}
as $k,\ell \to\infty$, and so in this case the $j^{\rm{th}}$
function in the above product is asymptotically close in
$L^2(m_{K_1\times K_2})$ to the function
\begin{eqnarray*}
&&\exp\big(-2\pi\rm{i}a_{2,j}\{\chi_{2,j}(\phi(2\ell h_2
h\bf{e}_1),\phi(-h\bf{e}_2)\}\\
&&\quad\quad\quad\quad\quad\quad\cdot(n_k -
n_\ell)\lfr\{\g_{2,j}(u,v)\} + \{\g_{2,j}(\phi(2\ell h_1
hm_k\bf{e}_1),\phi(-hm_k\bf{e}_2))\}\rfr\big)\\
&&= \chi_{2,j}(\phi(2\ell h_2
h\bf{e}_1),\phi(-h\bf{e}_2))^{a_{2,j}(n_\ell -
n_k)\lfr\{\g_{2,j}(u,v)\} + \{\g_{2,j}(\phi(2\ell h_1
hm_k\bf{e}_1),\phi(-hm_k\bf{e}_2))\}\rfr},
\end{eqnarray*}
and this is close to $1$ for either of the possible values ($0$ or
$1$) of $\lfr\{\g_{2,j}(u,v)\} + \{\g_{2,j}(\phi(2\ell h_1
hm_k\bf{e}_1),\phi(-hm_k\bf{e}_2))\}\rfr$, because $a_{2,j}$ is a
fixed integer and
\[(\phi(2\ell h_2 h\bf{e}_1),\phi(-h\bf{e}_2))^{n_k}\approx
(\phi(2\ell h_2 h\bf{e}_1),\phi(-h\bf{e}_2))^{n_\ell}\] when $k$ and
$\ell$ are large.

On the other hand, if $j \in \{1,2,\ldots,J_2\}\setminus \cal{J}$
then $\g_{2,j}(\phi(2\ell h_1h,-h))$ is a root of unity, and so
since the sequence $\phi(m_k\bf{q}_1)$ converges the values
$\g_{2,j}(\phi(2\ell h_1hm_k,-hm_k))$ are eventually constant. Once
this is so, of course we have
\begin{multline*}
\lfr\{\g_{2,j}(u,v)\} + \{\g_{2,j}(\phi(2\ell h_1
hm_k\bf{e}_1),\phi(-hm_k\bf{e}_2))\}\rfr\\
= \lfr\{\g_{2,j}(u,v)\} + \{\g_{2,j}(\phi(2\ell h_1
hm_\ell\bf{e}_1),\phi(-hm_\ell\bf{e}_2))\}\rfr,
\end{multline*}
for all $(u,v)\in K_1\times K_2$ and we may complete the proof of
this case as above.

\textbf{Remark}\quad It is for the above argument that we must make
a restriction such as $\phi(m_k\bf{q}_1) \in U_{\a/2}$.  Indeed,
without this we might have chosen a limit point
$(u^\circ_1,v^\circ_1)$ for which $\g_{2,j}(u^\circ_1,v^\circ_1) =
0$ for some $j \in \cal{J}$, and in this case it will generally
happen that there are large $k$ and $\ell$ for which, say,
$\{\g_{2,j}(\phi(2\ell h_1 hm_k\bf{e}_1),\phi(-hm_k\bf{e}_2))\}$ is
very slightly more than $0$ but $\{\g_{2,j}(\phi(2\ell h_1
hm_\ell\bf{e}_1),\phi(-hm_\ell\bf{e}_2))\}$ is very slightly less
than $1$.  This disrupts the above argument that the last factor in
our large product is close to $1$, and we find instead that it might
be close to some other constant, which seems to be hard to account
for in the desired expression $p(m_k,n_k) - p(m_\ell,n_\ell)$. \fin
\end{itemize}

Putting the above approximations together we obtain that for $k$ and
$\ell$ sufficiently large we have
\begin{eqnarray*}
&&\t{c}_1(u,v)^{m_k}\cdot \t{c}_2(u\phi(2\ell
h_1hm_k,0),v\phi(0,-hm_k))^{n_k}\\
&&\quad\quad\quad\quad\quad\quad\quad\quad\quad\quad\quad\quad\cdot\overline{\t{c}_1(u,v)^{m_\ell}\cdot
\t{c}_2(u\phi(2\ell
h_1hm_\ell,0),v\phi(0,-hm_\ell))^{n_\ell}} \\
&&\approx \a_1^{m_k - m_\ell}\cdot\big(\rm{character}\big)\cdot
\a_2^{n_k -
n_\ell}\cdot\big(\rm{character}\big)\\
&&\quad\cdot
\exp\Big(2\pi\rm{i}\sum_{j=1}^{J_2}a_{2,j}\{\chi_{2,j}(\phi(2\ell
h_2 h\bf{e}_1),\phi(-h\bf{e}_2))\}\big(n_k\{\g_{2,j}(\phi(2\ell h_1
hm_k\bf{e}_1),\phi(-hm_k\bf{e}_2))\}\\
&&\quad\quad\quad\quad\quad\quad\quad\quad\quad\quad\quad\quad -
n_\ell\{\g_{2,j}(\phi(2\ell h_1
hm_\ell\bf{e}_1),\phi(-hm_\ell\bf{e}_2))\}\big)\Big),
\end{eqnarray*}
so defining
\begin{multline*}
p(m,n) = \{\a_1^m\} + \{\a_2^n\}\\ +
\sum_{j=1}^{J_2}a_{2,j}\{\chi_{2,j}(\phi(2\ell h_2
h\bf{e}_1),\phi(-h\bf{e}_2))\}n_k\{\g_{2,j}(\phi(2\ell h_1
hm_k\bf{e}_1),\phi(-hm_k\bf{e}_2))\}
\end{multline*}
we see that this is a gen-polynomial not depending on $\a$ that has
the desired property. \qed

In Proposition~\ref{prop:compactification} we begin to see the
makings of the simplification of the expressions
\[\chi_1(\s((\ell n^2 +
an,0),z))\cdot\chi_2(\s((\ell n^2 + an,n),z)),\] that was promised
immediately after the proof of
Lemma~\ref{lem:further-reduced-polyconv}, although it will require
some more manipulation before the above proposition bears on this
expression directly.

\begin{cor}\label{cor:of-compactification}
If $p:\bbZ^2\to \bbR$ is the gen-polynomial of
Proposition~\ref{prop:compactification} then for any $\eps > 0$
there are some $K \geq 1$, functions $\xi_1$, $\xi_2$, \ldots,
$\xi_K \in L^2(m_Z)$ and characters $\chi_{i,1}$, $\chi_{i,2}$,
\ldots, $\chi_{i,K} \in \widehat{K_1\times K_2}$ for $i=1,2$ such
that
\begin{multline*}
\chi_2\circ\s((mq_{11} + nq_{21},0),z)\cdot\chi_1\circ\s((mq_{11}
+ nq_{21},mq_{12} + nq_{22}),z)\cdot\overline{\chi_1}\circ\s((0,mq_{12} + nq_{22}),z)\\
\approx_\eps \exp(-2\pi\rm{i}p(m,n))\cdot \sum_{k=1}^K
\chi_{1,k}(\phi(m\bf{q}_1))\chi_{2,k}(\phi(n\bf{q}_2))\cdot \xi_k(z)
\end{multline*}
in $L^2(m_Z)$ for every $m \in \bbZ$ such that $\phi(m\bf{q}_1) \in
U_\a$.
\end{cor}

\textbf{Proof}\quad Letting $C_\a$ be the Borel function $Z\times
(K_1\times K_2)^2 \to \Sone$ output by
Proposition~\ref{prop:compactification}, it will suffice to prove
that there are $\xi_1$, $\xi_2$, \ldots, $\xi_K \in L^2(m_Z)$ and
characters $\chi_{i,1}$, $\chi_{i,2}$, \ldots, $\chi_{i,K} \in
\widehat{K_1\times K_2}$ as above such that
\[C_\a(\,\cdot\,,u_1,v_1,u_2,v_2) \approx_\eps \sum_{k=1}^K
\chi_{1,k}(u_1,v_1)\chi_{2,k}(u_2,v_2)\cdot
\xi_k\quad\quad\hbox{in}\ L^2(m_Z)\] for all $(u,v)\in (K_1\times
K_2)^2$.

Proposition~\ref{prop:compactification} gives us that the map
$(u_1,v_1,u_2,v_2)\mapsto C_\a(\,\cdot\,,u_1,v_1,u_2,v_2)$ is
continuous from $(K_1\times K_2)^2$ into $L^2(m_Z)$. This implies
that its image is compact, and so lies within the
$(\eps/2)$-neighbourhood of some finite-dimensional subspace of
$L^2(m_Z)$; let $\xi_1$, $\xi_2$, \ldots, $\xi_K$ be a basis for
that subspace.  Simply by projecting onto this subspace it follows
that we can approximate the map $(u_1,v_1,u_2,v_2)\mapsto
C_\a(\,\cdot\,,u_1,v_1,u_2,v_2)$ uniformly in $(u_1,v_1,u_2,v_2)$ by
some map of the form
\[\sum_{m=1}^M C_{\a,m}(u_1,v_1,u_2,v_2)\cdot \xi_m\]
with each $C_{\a,m}:(K_1\times K_2)^2\to \bbC$ a continuous
function.

However, now the Stone-Weierstrass Theorem gives for each $C_{\a,m}$
a trigonometric polynomial $(K_1\times K_2)^2\to \bbC$ that
approximates $C_{\a,m}$ uniformly to within $\eps/(2(\|\xi_1\|_2 +
\ldots + \|\xi_K\|_2))$. Replacing each $c_m$ by this trigonometric
polynomial in our first approximant to $C_\a$ and re-arranging the
terms gives the result. \qed

\subsection{Completion of the proof}

We are finally ready to prove
Proposition~\ref{prop:reduced-polyconv}.

\textbf{Proof of Proposition~\ref{prop:reduced-polyconv}}\quad By
Lemma~\ref{lem:further-reduced-polyconv} we need only prove
convergence of the averages
\[\frac{1}{N}\sum_{n=1}^N
\theta_1^{\ell n^2 + an}\theta_2^n\cdot\chi_1(\s((\ell n^2 +
an,0),z))\cdot\chi_2(\s((\ell n^2 + an,n),z))\cdot f(T_2^n(x))\] for
any $\theta_1,\theta_2 \in \Sone$, and by
Lemma~\ref{lem:final-choice-of-hs} we may restrict our attention to
the case covered by the above results, and in particular
Corollary~\ref{cor:of-compactification}.  We will handle this case
in two steps.

\quad\textbf{Step 1}\quad We first need a simple but slightly fiddly
re-arrangement in order to bring
Corollary~\ref{cor:of-compactification} to bear, because it applies
only to the sublattice $\G = \bbZ\bf{q}_1 + \bbZ\bf{q}_2$ of
$\bbZ^2$. To do this, let us choose an integer $\ell_1 \geq 1$ so
that $\ell_1\bbZ^2 \leq \G$ and break up the above average as
\begin{eqnarray*}
&&\frac{1}{\ell_1}\sum_{j=1}^{\ell_1}\frac{1}{\lfr
N/\ell_1\rfr}\sum_{n=0}^{\lfr N/\ell_1\rfr} \theta_1^{\ell (\ell_1 n
+ j)^2 + a(\ell_1 n + j)}\theta_2^{\ell_1 n + j}\cdot\chi_1(\s((\ell
(\ell_1 n + j)^2 + a(\ell_1 n + j),0),z))\\
&&\quad\quad\quad\quad\quad\quad\quad\quad\cdot\chi_2(\s((\ell
(\ell_1 n + j)^2 + a(\ell_1 n + j),\ell_1 n + j),z))\cdot
f(T_2^{\ell_1 n}(T_2^j(x)))\\
&&\quad\quad\quad\quad\quad\quad\quad\quad + R\\
&&= \frac{1}{\ell_1}\sum_{j=1}^{\ell_1}\theta_1^{\ell j^2 +
aj}\theta_2^j\frac{1}{\lfr N/\ell_1\rfr}\sum_{n=0}^{\lfr
N/\ell_1\rfr} \theta_1^{\ell_1(\ell \ell_1 n^2 + 2\ell j n + a
n)}\theta_2^{\ell_1 n}\cdot\chi_1(\s((\ell
(\ell_1 n + j)^2 + a(\ell_1 n + j),0),z))\\
&&\quad\quad\quad\quad\quad\quad\quad\quad\cdot\chi_2(\s((\ell
(\ell_1 n + j)^2 + a(\ell_1 n + j),\ell_1 n + j),z))\cdot
f(T_2^{\ell_1 n}(T_2^j(x)))\\
&&\quad\quad\quad\quad\quad\quad\quad\quad + R
\end{eqnarray*}
where the remainder term satisfies $\|R\|_2 = \rm{O}(1/N)$, and so
may henceforth be ignored. It will suffice to prove that for each $j
\in \{1,2,\ldots,\ell_1\}$ the inner average over $0 \leq n \leq
\lfr N/\ell_1\rfr$ converges in $L^2(\mu)$.

To simplify these inner averages, let us recall the consequence of
the defining equation for the cocycle $\s$ that we have
factorizations
\begin{multline*}
\chi_1(\s((\ell (\ell_1 n + j)^2 + a(\ell_1 n + j),0),z))\\ =
\chi_1\big(\s((\ell_1(\ell \ell_1 n^2 + 2\ell jn +
an),0),z)\big)\cdot\chi_1\big(\s((\ell j^2 + aj,0),z\cdot
\phi(\ell_1(\ell \ell_1 n^2 + 2\ell jn + an)\bf{e}_1))\big)
\end{multline*}
and similarly
\begin{eqnarray*}
&&\chi_2(\s((\ell (\ell_1 n + j)^2 + a(\ell_1 n + j),(\ell_1 n +
j)),z))\\ &&= \chi_2\big(\s((\ell_1 (\ell \ell_1 n^2 + 2\ell jn +
an),\ell_1 n),z)\big)\\
&&\quad\quad\quad\quad\cdot\chi_2\big(\s((\ell j^2 + aj,j),z\cdot
\phi(\ell_1 (\ell \ell_1 n^2 + 2\ell jn + an)\bf{e}_1 + \ell_1
n\bf{e}_2))\big).
\end{eqnarray*}

Now, for fixed integers $\ell_1$ and $j$ the second factors in the
factorizations above correspond to the functions
\[h_1:z\mapsto \chi_1(\s((\ell j^2 + aj,0),z))\]
and
\[h_2:z\mapsto \chi_2(\s((\ell j^2 + aj,j),z)),\]
so that we can write
\begin{multline*}
\chi_1\big(\s((\ell j^2 + aj,0),z\cdot \phi(\ell_1(\ell \ell_1 n^2 +
2\ell jn + an)\bf{e}_1))\big)\\
\cdot\chi_2\big(\s((\ell j^2 + aj,j),z\cdot \phi(\ell_1 (\ell \ell_1
n^2 + 2\ell jn + an)\bf{e}_1 + \ell_1 n\bf{e}_2))\big)
\end{multline*}
as
\[h_1(R_{\phi(\ell_1(\ell \ell_1 n^2 + 2\ell jn +
an)\bf{e}_1)}z)\cdot h_2(R_{\phi(\ell_1 (\ell \ell_1 n^2 + 2\ell jn
+ an)\bf{e}_1 + \ell_1 n\bf{e}_2)}z).\] Since we may approximate
each of $h_1$ and $h_2$ arbitrarily well in $L^2(m_Z)$ by a
trigonometric polynomial on $Z$, it follows by continuity and
multilinearity that the desired convergence will follow if we prove
it instead for the averages
\begin{eqnarray*}
&&\frac{1}{N}\sum_{n=0}^{N} \theta_1^{\ell_1(\ell \ell_1 n^2 + 2\ell
j n + a n)}\theta_2^{\ell_1 n}\cdot\chi_1(\s((\ell_1(\ell \ell_1 n^2
+
2\ell jn + an),0),z))\\
&&\quad\quad\quad\quad\quad\quad\quad\quad
\cdot\chi_2(\s((\ell_1(\ell \ell_1 n^2 + 2\ell jn + an),\ell_1
n),z))\\
&&\quad\quad\quad\quad\quad\quad\quad\quad\cdot
h_1(R_{\phi(\ell_1(\ell \ell_1 n^2 + 2\ell jn + an)\bf{e}_1)}z)\cdot
h_2(R_{\phi(\ell_1 (\ell \ell_1 n^2 + 2\ell jn + an)\bf{e}_1 +
\ell_1 n\bf{e}_2)}z)\\
&&\quad\quad\quad\quad\quad\quad\quad\quad\cdot f(T_2^{\ell_1
n}(T_2^j(x)))
\end{eqnarray*}
where each of $h_1$ and $h_2$ is a character.  In that case
\[h_1(R_{\phi(\ell_1(\ell \ell_1 n^2 + 2\ell jn + an)\bf{e}_1)}z) = h_1(\phi(\bf{e}_1))^{\ell_1(\ell \ell_1 n^2 + 2\ell jn + an)}h_1(z)\]
and similarly for $h_2$, so by taking the $n$-independent functions
$h_1(z)$ and $h_2(z)$ outside the average and adjusting the values
of $\theta_1$ and $\theta_2$ we can now drop the mention of these
functions $h_i$ altogether to leave the averages
\begin{multline*}
\frac{1}{N}\sum_{n=1}^{N} \theta_1^{\ell_1(\ell \ell_1 n^2 + 2\ell
jn + a n)}\theta_2^{\ell_1 n}\cdot\chi_1(\s((\ell_1(\ell \ell_1 n^2
+
2\ell jn + an),0),z))\\
\cdot\chi_2(\s((\ell_1(\ell \ell_1 n^2 + 2\ell jn + an),\ell_1
n),z))\cdot f(T_2^{\ell_1 n}(T_2^j(x))).
\end{multline*}

\quad\textbf{Step 2}\quad The value of the simplification achieved
in Step 1 above is that now by our choice of $\ell_1$ we have
$(\ell_1(\ell \ell_1 n^2 + 2\ell jn + an),-\ell_1 n) \in \G$ for all
$n\geq 1$.  In particular, it follows that there are independent
linear forms $L_1,L_2:\G\to \bbZ$ such that
\begin{multline*}
(\ell_1(\ell \ell_1 n^2 + 2\ell jn + an),-\ell_1 n) =
L_1(\ell_1(\ell \ell_1 n^2 + 2\ell jn + an),-\ell_1 n)\bf{q}_1\\ +
L_2(\ell_1(\ell \ell_1 n^2 + 2\ell jn + an),-\ell_1 n)\bf{q}_2
\end{multline*}
for all $n$.  Let us abbreviate $\vec{L} := (L_1,L_2)$ and
\[Q_i(n) := L_i(\ell_1(\ell \ell_1 n^2 + 2\ell jn + an),-\ell_1 n),\]
so that $Q_1$ and $Q_2$ are two non-constant, linearly independent
quadratic functions $\bbZ\to\bbZ$.

Now recall the open subsets $U_\a \subseteq K_1\times K_2$
introduced in Proposition~\ref{prop:compactification}.  The set
\[\{m \in \bbZ:\ \phi(m\bf{q}_1) \in U_\a\}\]
is a Bohr set in $\bbZ$, and by construction it is defined by
irrational phases. Consequently, the multidimensional version of
Weyl's Equidistribution Theorem (see, for instance, Theorem 1.6.4 in
Kuipers and Niederreiter~\cite{KuiNie74}) gives that the set
\[E_\a := \{n \geq 1:\ \phi(Q_1(n)\bf{q}_1) \in U_\a\}\] has asymptotic density equal to
$m_{K_1\times K_2}(U_\a) > 1-\a$.  Since the terms of our average
\begin{multline*}
\frac{1}{N}\sum_{n=1}^{N} \theta_1^{\ell_1(\ell \ell_1 n^2 + 2\ell
jn + a n)}\theta_2^{\ell_1 n}\cdot\chi_1(\s((\ell_1(\ell \ell_1 n^2
+ 2\ell kn + an),0),z))\\
\cdot\chi_2(\s((\ell_1(\ell \ell_1 n^2 + 2\ell jn + an),\ell_1
n),z))\cdot f(T_2^{\ell_1 n}(T_2^j(x))).
\end{multline*}
are uniformly bounded in $L^\infty$, to prove norm convergence it
suffices to prove it for the related averages in which we restrict
the sum to those $n$ that lie inside some subset of $\bbN$, provided
we can choose that set to have arbitrarily high asymptotic density.
Hence, in particular, it will suffice to prove for every $\a > 0$
the convergence of the averages in which we restrict the summation
to $n \in \{1,2,\ldots,N\}\cap E_\a$.

Now, Corollary~\ref{cor:of-compactification} gives a gen-polynomial
$p:\bbZ^2\to\bbR$, and for any $\a
> 0$ and $\eps > 0$ some functions $\xi_1$,
$\xi_2$, \ldots, $\xi_K \in L^2(m_Z)$ and characters $\chi_{i,1}$,
$\chi_{i,2}$, \ldots, $\chi_{i,K}\in \widehat{K_1\times K_2}$ for
$i=1,2$ such that
\begin{eqnarray*}
&&\chi_2(\s((\ell_1(\ell \ell_1 n^2 + 2\ell jn +
an),0),z))\cdot\chi_1(\s((\ell_1(\ell \ell_1 n^2 + 2\ell jn +
an),-\ell_1 n),z))\\
&&\quad\quad\quad\quad\quad\quad\quad\quad \cdot\overline{\chi_1(\s((0,-\ell_1 n),z))}\\
&& \approx_\eps \exp(-2\pi\rm{i}\,p\circ\vec{L}(\ell_1(\ell \ell_1
n^2 + 2\ell jn
+ an),-\ell_1 n))\\
&&\quad\quad\quad\quad\quad\quad\quad\quad\cdot
\sum_{k=1}^K\chi_{1,k}(\phi(Q_1(n)\bf{q}_1))\chi_{2,k}(\phi(Q_2(n)\bf{q}_2))\cdot
\xi_k(z)
\end{eqnarray*}
in $L^2(m_Z)$ for all $n\geq 1$ with $n\in E_\a$. Using the cocycle
equation we can re-write
\begin{eqnarray*}
&&\chi_1(\s((\ell_1(\ell \ell_1 n^2 + 2\ell jn +
an),0),z))\cdot\chi_2(\s((\ell_1(\ell \ell_1 n^2 + 2\ell jn +
an),\ell_1 n),z))\\
&&= \chi_1(\s((\ell_1(\ell \ell_1 n^2 + 2\ell jn + an),-\ell_1
n),z\phi(0,\ell_1 n)))\cdot \chi_1(\s((0,\ell_1 n),z))\\
&&\quad\quad\quad\quad\cdot \chi_2(\s((\ell_1(\ell \ell_1 n^2 +
2\ell jn + an),0),z\phi(0,\ell_1 n)))\cdot \chi_2(\s((0,\ell_1
n),z)),
\end{eqnarray*}
and now substituting from the above approximation we see that for
all $n\in \bbN\cap E_\a$ this lies within $\eps$ in $L^2(m_Z)$ of
\begin{eqnarray*}
&&\exp(-2\pi\rm{i}\,p\circ\vec{L}(\ell_1(\ell \ell_1 n^2 + 2\ell jn
+ an),-\ell_1 n))\\
&&\quad\quad\quad\quad\cdot
\sum_{k=1}^K\chi_{1,k}(\phi(Q_1(n)\bf{q}_1))\chi_{2,k}(\phi(Q_2(n)\bf{q}_2))\cdot
\xi_k(z\phi(0,\ell_1 n))\\
&&\quad\quad\quad\quad\cdot \chi_2(\s((0,\ell_1
n),z))\cdot\chi_1(\s((0,\ell_1 n),z))\cdot \chi_1(\s((0,-\ell_1
n),z\phi(0,\ell_1 n)))\\
&&= \exp(-2\pi\rm{i}\,p\circ\vec{L}(\ell_1(\ell \ell_1 n^2 + 2\ell
jn
+ an),-\ell_1 n))\\
&&\quad\quad\quad\quad\cdot
\sum_{k=1}^K\chi_{1,k}(\phi(Q_1(n)\bf{q}_1))\chi_{2,k}(\phi(Q_2(n)\bf{q}_2))\cdot
\xi_k(z\phi(0,\ell_1 n))\chi_2(\s((0,\ell_1 n),z)),
\end{eqnarray*}
using that the cocycle equation also gives
\[\s((0,\ell_1 n),z)\cdot \s((0,-\ell_1 n),z\phi(0,\ell_1 n)) = \s((0,0),z) = 1.\]

Since $\eps > 0$ was arbitrary we may substitute this approximation
into our averages above and appeal again to multilinearity to deduce
that it suffices to prove instead the norm convergence of the
averages
\begin{eqnarray*}
&&\frac{1}{N}\sum_{1 \leq n \leq N,\,n\in E_\a}
\theta_1^{\ell_1(\ell \ell_1 n^2 + 2\ell j n + a n)}\theta_2^{\ell_1
n}\cdot\exp(-2\pi\rm{i}\,p\circ\vec{L}(\ell_1(\ell \ell_1 n^2 +
2\ell jn
+ an),-\ell_1 n))\\
&&\quad\quad\quad\quad\quad\quad\quad\quad\quad\quad\quad\quad\cdot\t{\chi}_1(\phi(Q_1(n)\bf{q}_1))\t{\chi}_2(\phi(Q_2(n)\bf{q}_2))\\
&&\quad\quad\quad\quad\quad\quad\quad\quad\quad\quad\quad\quad\cdot
\xi(z\phi(0,\ell_1 n))\cdot\chi_2(\s((0,\ell_1 n),z))\cdot
f(T_2^{\ell_1 n}(T_2^k(x)))
\end{eqnarray*}
for any two characters $\t{\chi}_1,\t{\chi}_2 \in \widehat{K_1\times
K_2}$ and any fixed function $\xi \in L^2(m_Z)$.

Finally, in order to prove convergence we may freely insert the
$n$-independent function $(z,a)\mapsto \chi_2(a)$ into these
averages, because this function is bounded away from zero. This
trick now leads to the simplification
\[\xi(z\phi(0,\ell_1 n))\cdot \chi_2(a)\cdot\chi_2(\s((0,\ell_1 n),z))\cdot f(T_2^{\ell_1
n}(T_2^k(x))) = F(T_2^{\ell_1 n}(x))\] where $F(x) :=
\xi(z)f(T_2^k(x))\chi_2(a)$ (remembering that $(z,a) = \eta(x)$).
On the other hand, the expression
\begin{multline*}
\theta_1^{\ell_1(\ell \ell_1 n^2 + 2\ell j n + a n)}\theta_2^{\ell_1
n}\cdot\exp(-2\pi\rm{i}\,p\circ\vec{L}(\ell_1(\ell \ell_1 n^2 +
2\ell jn
+ an),-\ell_1 n))\\
\cdot\t{\chi}_1(\phi(Q_1(n)\bf{q}_1))\t{\chi}_2(\phi(Q_2(n)\bf{q}_2))
\end{multline*}
clearly just defines an expression of the form $\exp(\rm{i}Q_3(n))$
for $Q_3:\bbZ\to\bbR$ a new gen-polynomial, and so the rather
unwieldy averages above can be written in the simple form
\[\frac{1}{N}\sum_{1 \leq n \leq N,\,n\in E_\a}\exp(\rm{i}Q_3(n))\cdot F\circ T_2^{\ell_1 n} = \frac{1}{N}\sum_{n = 1}^N 1_{E_\a}(n)\cdot\exp(\rm{i}Q_3(n))\cdot F\circ T_2^{\ell_1 n}.\]

Next, the indicator function $1_{E_\a}$ corresponds to a quadratic
Bohr set, and so among $1$-bounded functions on $\bbN$ it can be
approximated in density by linear combinations of gen-polynomial
maps taking values in $\Sone$.  Appealing once again to
multilinearity, it follows that we need only prove convergence of
the averages
\[\frac{1}{N}\sum_{n = 1}^N
\exp(\rm{i}Q_4(n))\cdot F\circ T_2^{\ell_1 n}\] for a
suitably-enlarged list of possible gen-polynomials $Q_4$.

The convergence of these now follows from the results of Bergelson
and Leibman in~\cite{BerLei07} (or could probably also be deduced
from the results of Host and Kra in their related
paper~\cite{HosKra07}). In particular, a simple appeal to the
spectral theorem and Corollary 0.26 in~\cite{BerLei07} shows that
whenever $(U^t_1)_{t\in\bbR}$ and $U_2$ are respectively a unitary
flow and a unitary operator acting on a Hilbert space $\frH$ and
$Q'_1:\bbZ\to\bbR$ and $Q'_2:\bbZ\to\bbR$ are generalized
polynomials, then the sequence of operator averages
\[\frac{1}{N}\sum_{n=1}^NU_1^{Q'_1(n)}U_2^{Q'_2(n)}\]
converges in the strong operator topology. (In fact this result lies
just between two further corollaries that Bergelson and Leibman
obtain explicitly in~\cite{BerLei07}, Corollary 0.27 concerning
tuples of flows and Corollary 0.28 concerning tuples of single
operators.) This implies the convergence we need in the case when
$\frH = L^2(\mu)$, $U^t_1$ is multiplication by $\exp(\rm{i}t)$,
$U_2$ is the Koopman operator of $T_2^{\ell_1}$, $Q'_1(n) := Q_4(n)$
and $Q'_2(n) := n$.

This completes the proof of Proposition~\ref{prop:reduced-polyconv},
and hence of Theorem~\ref{thm:polyconv}. \qed

\textbf{Remark}\quad In~\cite{HosKra01} Host and Kra augment their
proof of convergence with a description of the limit function that
emerges.  Although the last step in our proof of convergence above
is rather similar to their argument, the other stages in our
reduction leave it much less clear just how the limit function can
be described in our case, even after passing to a suitable extended
system. \fin

\appendix

\section{Moore cohomology}\label{app:cohom}

We collect here the definition of Moore's measurable cohomology
theory for locally compact groups and some of its basic properties
that are needed in Section~\ref{sec:cohom}.  Some of the result proved below can be improved using the continuity results of~\cite{Aus--cohomcty}, but we have left them in the form in which they were presented before the appearance of that paper in order to remain consistent with the main text above.

The most convenient definition of this cohomology theory for our
purposes is in terms of the measurable homogeneous bar resolution.
We recall this here for completeness, noting that it is shown by
Moore to be equivalent to various more abstract definitions, and to
support the usual functorial cohomological machinery of discrete
group cohomology (particularly the procedure of dimension-shifting and the Hochschild-Serre
spectral sequence).

\begin{dfn}[Measurable cohomology for locally compact groups]\label{dfn:meas-cohom}
If $A$ is a locally compact group, $R$ is a Polish Abelian group and
$\a:A\curvearrowright R$ is a continuous left-action by
automorphisms, then we define the \textbf{measurable cohomology of
$A$ with coefficients in $(R,\a)$} as the (discrete) cohomology of
the chain complex
\[0 \longrightarrow R\stackrel{d}{\longrightarrow} \C(A,R)\stackrel{d}{\longrightarrow} \C(A^2,R)\stackrel{d}{\longrightarrow}\ldots\]
with chain maps defined by
\begin{multline*}
d\phi(a_1,a_2,\ldots,a_{n+1}) := \a^{a_{n+1}}(\phi(a_1,a_2,\ldots,a_n))\\
+
\sum_{i=1}^n(-1)^{n+1-i}\phi(a_1,a_2,\ldots,a_i+a_{i+1},\ldots,a_{n+1})
+ (-1)^{n+1}\phi(a_2,a_3,\ldots,a_{n+1}).
\end{multline*}
We write $\Z^n(A,R) := \ker d|_{\C(A^n,R)}$ for the subgroup of
\textbf{cocycles} in $\cal{C}(A^n,R)$ and $\B^n(A,R):=
\rm{img}\,d|_{\C(A^{n-1},R)}$ for the subgroup of
\textbf{coboundaries}, and in these terms the cohomology groups are
the \emph{discrete} groups
\[\rmH^n(A,R) := \frac{\Z^n(A,R)}{\B^n(A,R)}.\]
\end{dfn}

We warn the reader that this definition of differential is
`back-to-front' compared with the usual conventions of discrete
group cohomology (see Section 6.5 of Weibel~\cite{Wei94}) so as to
be better adapted to our present setting; it is clear that this
makes only a cosmetic difference to the theory.

It is easy to find examples in which the measurability condition on
the above cochains makes a large difference to the cohomology groups
that result.  Perhaps most simply, it is easy to check that for any
Polish Abelian group $A$ with trivial $\bbR$-action we have that
$\rmH^1(\bbR,A)$ is isomorphic to the group of continuous
homomorphisms $\bbR\to R$, whereas $\rmH^1(\bbR_{\rm{discrete}},A)$
is a discrete Abelian group of uncountable rank in general.

Moore also gives some discussion in~\cite{Moo76(gr-cohomIV)} of
possible topologies on the cohomology groups themselves.  However,
the obvious candidate topologies are often badly behaved (for
example, by being non-Hausdorff, as in the well-known case when $A =
\bbZ$, $R = \C(X,\mu)$ and $\a^n(f) = f\circ T^n$ for some
nontrivial aperiodic action $T:\bbZ\curvearrowright (X,\mu)$), and
we will not need a topology on these groups here.

We now state three important calculational results from Moore's
papers that we will need later.  Their proofs employ the basic
functorial machinery of this cohomology theory that are set up
there, particularly the Hochschild-Serre spectral sequence and its
corollary, the restriction-inflation exact sequence; we omit them
here.

\begin{prop}[Second cohomology and the fundamental
group]\label{prop:second-cohom-of-tori} If $Z$ is a compact
connected Lie group with fundamental group $\pi_1(Z)$, and
$\pi_1(Z)_{\rm{tor}}$ is the torsion subgroup of $\pi_1(Z)$, then
there is a canonical isomorphism $\rmH^2(Z,\bbT)\cong
\widehat{\pi_1(Z)_{\rm{tor}}}$.  In particular, $\rmH^2(\bbT^d,\bbT)
= 0$ for all $d\geq 1$.
\end{prop}

\textbf{Proof}\quad This is Proposition 2.1 in part I
of~\cite{Moo64(gr-cohomI-II)}. \qed

\begin{prop}[Continuity of $\rmH^2$ under inverse and direct
limits]\label{prop:inv-lim-cohom} If $Z = \lim_{m\leftarrow}Z_{(m)}$
is an inverse limit of compact groups and $A =
\lim_{m\rightarrow}A_{(m)}$ is a direct limit of countable discrete
groups with trivial $Z$-action then
\begin{enumerate}
\item $\rmH^2(Z,A)$ is isomorphic to the direct limit of the
groups $\rmH^2(Z_{(m)},A_{(m)})$ under the compositions of the
inflation maps $\rm{inf}:\rmH^2(Z_{(m)},A_{(m)})\into
\rmH^2(Z,A_{(m)})$ with the embeddings $A_{(m)}\to A$, and
\item $\rmH^2(Z,\bbT)$ is similarly isomorphic to the direct limit of the
groups $\rmH^2(Z_{(m)},\bbT)$ under the inflation maps
$\rm{inf}:\rmH^2(Z_{(m)},\bbT)\to \rmH^2(Z,\bbT)$.
\end{enumerate}
\end{prop}

\textbf{Proof}\quad These are special cases of Theorems 2.1 and 2.2
of Part I of~\cite{Moo64(gr-cohomI-II)} (observing that any compact
Abelian group is almost connected). \qed

\begin{lem}[Real cohomology of compact Abelian
groups]\label{lem:vanishing-real-cohom} If $Z$ is a compact Abelian
group then $\rmH^1(Z,\bbR) = \rmH^2(Z,\bbR) = 0$. If $Z$ is a
finite-dimensional compact Abelian group then this extends to
$\rmH^n(Z,\bbR) = 0$ for all $n > 0$.
\end{lem}

\textbf{Proof}\quad The first conclusion is part of Theorem 2.3 in
Part I of Moore~\cite{Moo64(gr-cohomI-II)}, and the second follows
from the identification for compact Lie groups of Moore's measurable
cohomology with the cohomology theory for topological groups defined
using classifying spaces, as outlined by Moore at the end
of~\cite{Moo76(gr-cohomIII)} and described in detail by Wigner
in~\cite{Wig73}. \qed

\begin{lem}[Integral degree-$2$ cohomology]\label{lem:integer-cohom} If $Z$ is a compact Abelian group
then $\rmH^2(Z,\bbZ) \cong \widehat{Z}$, where the isomorphism is
given by assigning to $\g \in \widehat{Z}$ the $2$-cocycle
\[\k_\g(z,w) := \lfr\{\g(z)\} + \{\g(w)\}\rfr.\]
\end{lem}

\textbf{Proof}\quad Suppose that $\k:Z\times Z\to \bbZ \subset \bbR$
is a Borel $2$-cocycle.  By the previous lemma we know there is some
$a:Z\to \bbR$ such that $da = \k$, but of course this $a$ may not be
$\bbZ$-valued.   However, since $\k$ does take values in $\bbZ$, we
know that
\[a(z) + a(w) - a(z+w) + \bbZ = \k(z,w) + \bbZ = \bbZ\]
almost surely, so on composing with the quotient map $\bbR\to \bbT$
our $1$-cochain $a$ must descend to a measurable (and hence
continuous) character $\g \in \widehat{Z}$.  The map $a'(z) :=
\{\g(z)\} \in [0,1)$ clearly does give $\g$ upon composing with the
quotient, and on the other hand a direct computation gives
\[a'(z) + a'(w) - a'(z+w) = \k_\g(z,w)\]
(since $a + b - \{a + b\} \equiv \lfr a + b\rfr$ for $a,b\in
[0,1)$).  Therefore $\k - \k_\g = d(a - a')$ with $a - a'$ taking
values in $\bbZ$.

On the other hand any two $2$-cocycles of the form $\k_\g$ must give
rise to different homomorphisms above, and so they cannot be
cohomologous in $\Z^2(Z,\bbZ)$.  This completes the proof. \qed

\textbf{Remark}\quad In fact for $Z = \bbT^d$ the preceding lemma is
a special case of a rather more far-reaching description of the
integral cohomology. With the standard definition of cup product,
the cohomology ring $\rmH^\ast(\bbT^d,\bbZ)$ is isomorphic to the
polynomial ring $\bbZ[X_1,X_2,\ldots,X_d]$ graded so that each free
variable $X_i$ has degree two (so, in particular,
$\rmH^n(\bbT^d,\bbZ) = 0$ when $n$ is odd), and for even $n$ the
cochains
\begin{multline*}
c(\bf{t}_1,\bf{t}_2,\ldots,\bf{t}_n)\\ :=
\Big(\prod_{j=1}^d\prod_{i=1}^{\ell_j}\lfr\{t_{2i-1,2\ell_1 +
2\ell_2 + \cdots + 2\ell_{j-1} + j}\}+\{t_{2i,2\ell_1 + 2\ell_2 +
\cdots + 2\ell_{j-1} + j}\}\rfr\Big)
\end{multline*}
corresponding to the monomials $X_1^{\ell_1}X_2^{\ell_2}\cdots
X_d^{\ell_d}$ with $2\ell_1 + 2\ell_2 + \cdots + 2\ell_d = n$
comprise a free set of generators of $\rmH^n(\bbT^d,\bbZ)$, where we
write $\bf{t}_i = (t_{i,1},t_{i,2},\ldots,t_{i,d})\in\bbT^d$.  In
all cases these calculations can be performed directly using the
measurable versions of standard group cohomological machinery,
particularly the Hochschild-Serre spectral sequence, that are set up
in Moore's earlier papers~\cite{Moo64(gr-cohomI-II)}; or,
alternatively, they can be deduced from results of
Wigner~\cite{Wig73} showing that for $\bbT^d$ and these particular
target modules the Moore cohomology can be identified with various
other cohomology theories (such as that defined in terms of
classifying spaces, developed in detail for compact Abelian groups
by Hofmann and Mostert in~\cite{HofMos73}). \fin

The proof of Proposition~\ref{prop:combining-indiv-cobdry-eqns} in
Section~\ref{sec:cohom} will rest on the following rather more
detailed cohomological calculations.

\begin{lem}\label{lem:more-explicit-2-cocycs}
Suppose that $F$ is a finite Abelian group, $r\geq 0$, $G$ is
another locally compact Abelian group on which $\bbT^r\times F$ acts
trivially,
\[\k:(\bbT^r\times F)\times (\bbT^r\times F)\to G\]
is a $2$-cocycle and
\[\b:(\bbT^r\times F)^3\to G\]
is a $3$-cocycle.

Then
\begin{enumerate}
\item if $G = \bbT$ then $\k$ is cohomologous to a $2$-cocycle $\k'$
that depends only on the coordinates in $F$;
\item if $G = \bbZ$ then $\b$ is cohomologous to a $3$-cocycle $\b'$
that depends only on the coordinates in $F$;
\item if $G = \bbT$ and $\k$ is a $\bbT$-valued coboundary on $\bbT^r\times
F$ and depends only on coordinates in $F$, then $\k$ is is a
$\bbT$-valued coboundary on $F$;
\item if $G = \bbZ/n\bbZ$ then $\k$ is cohomologous to a $2$-cocycle $\k'$
of the form $\k'(z,w) := \k''(z,w) + \lfr\{\g(z)\} + \{\g(w)\}\rfr +
n\bbZ$ for some $\g \in \widehat{\bbT^r\times F}$ and some
$2$-cocycle $\k''$ that depends only on coordinates in $F$.
\end{enumerate}
\end{lem}

\textbf{Proof}\quad\textbf{1.}\quad The first conclusion follows
from the spectral sequence calculations of Section 3 in Part I of
Moore~\cite{Moo64(gr-cohomI-II)}.  In particular, the first two
layers of the Hochschild-Serre spectral sequence introduce a
filtering of groups
\[\rmH^2(\bbT^r\times F,\bbT) \geq K_1 \geq K_2 \geq \{0\}\]
where $K_1$ is identified with the subgroup of cohomology classes
containing a representative $2$-cocycle $\k$ such that
$\k|_{\bbT^r\times \bbT^r} = 0$ (that is, the kernel of the
restriction map to $\bbT^r$), $K_2$ with the further subgroup of
classes containing a representative that depends only coordinates in
$F$ (that is, the image of the inflation map), and such that
$K_2/K_1 \cong \rmH^1(F,\rmH^1(\bbT^r,\bbT))$ (where
$\rmH^1(\bbT^r,\bbT)$ is given the discrete topology).

However, Proposition~\ref{prop:second-cohom-of-tori} tells us that
$\rmH^2(\bbT^r,\bbT) = 0$, so for any $2$-cocycle $\k:(\bbT^r\times
F)\times (\bbT^r\times F)\to \bbT$ we can find some
$\a:\bbT^r\to\bbT$ such that $\k|_{\bbT^r\times \bbT^r} = d\a$. If
we lift $\a$ to $\bbT^r\times F$ under the coordinate projection
map, it follows that $\k - d\a$ is a cohomologous $2$-cocycle that
vanishes on $\bbT^r\times \bbT^r$, and so we have shown that in our
setting $\rmH^2(\bbT^r\times F,\bbT) = K_1$.

In addition, we know that $\rmH^1(\bbT^r,\bbT) = \widehat{\bbT^r}
\cong \bbZ^r$ is torsion-free, and so $\rmH^1(F,\rmH^1(\bbT^r,\bbT))
\cong \rm{Hom}(F,\bbZ^r) = 0$.  Thus in fact $\rmH^2(\bbT^r\times
F,\bbT) = K_2$, giving the first conclusion is proved.

\quad\textbf{2.}\quad This will follow from Part 1 and the
switchback maps of the long exact sequence
\begin{multline*}
\ldots \to \rmH^n(\bbT^r\times F,\bbZ)\to \rmH^n(\bbT^r\times
F,\bbR)\to \rmH^n(\bbT^r\times F,\bbT)\\
\stackrel{\rm{switchback}}{\longrightarrow} \rmH^{n+1}(\bbT^r\times
F,\bbZ)\to \rmH^{n+1}(\bbT^r\times F,\bbR)\to \ldots
\end{multline*}
corresponding to the presentation $\bbZ\into \bbR\onto \bbT$. By
Lemma~\ref{lem:vanishing-real-cohom} we have $\rmH^n(\bbT^r\times
F,\bbR)= 0$ for all $n\geq 1$, so this long exact sequence collapses
to a collection of isomorphisms
\[\rmH^n(\bbT^r\times F,\bbT) \cong \rmH^{n+1}(\bbT^r\times
F,\bbZ)\] which for $n=2$ directly enables us to appeal to Part 1.

More explicitly, given any $3$-cocycle $\b:(\bbT^r\times F)^3 \to
\bbZ$, we can express it as the coboundary of an $\bbR$-valued
$2$-cochain $\k:(\bbT^r\times F)\times (\bbT^r\times F)\to \bbR$,
and now since $\b$ takes values in $\bbZ$ it follows that $\k+\bbZ$
is a $\bbT$-valued $2$-cocycle. Therefore by Part 1 we can find some
$\a_0:\bbT^r\times F\to \bbT$ such that $\k_0' := (\k+\bbZ) - d\a_0$
depends only on coordinates in $F$.  Now let $\a:\bbT^r\times
F\to\bbR$ be a lift of $\a_0$ and $\k':(\bbT^r\times F)\times
(\bbT^r\times F)\to \bbR$ a lift of $\k_0'$ that depends only on
coordinates in $F$, so we must have that $\k'':= \k - d\a - \k'$ is
$\bbZ$-valued. It follows that $\b = d\k = d\k' + d\k''$, where
$\k'$ depends only on coordinates in $F$ and $\k''$ is
$\bbZ$-valued, as required.

\quad\textbf{3.}\quad We need to show that the inflation map
$\rm{inf}:\rmH^2(F,\bbT)\to\rmH^2(\bbT^r\times F,\bbT)$ is
injective. This follows from another consequence of Moore's spectral
sequence calculations: the measurable analog of Lyndon's
inflation-restriction exact sequence, derived in Section I.5 of Part
I of~\cite{Moo64(gr-cohomI-II)}. In our case this specializes to
\begin{multline*}
0 \to \rmH^1(F,\bbT)\stackrel{\rm{inf}}{\longrightarrow}
\rmH^1(\bbT^r\times F,\bbT) \stackrel{\rm{res}}{\longrightarrow}
\rmH^1(\bbT^r,\bbT)\\ \stackrel{\rm{tg}}{\longrightarrow}
\rmH^2(F,\bbT)\stackrel{\rm{inf}}{\longrightarrow}
\rm{inf}(\rmH^2(F,\bbT)) \leq \rmH^2(\bbT^r\times F,\bbT),
\end{multline*}
where $\rm{tg}$ is the so-called `transgression' map.  We do not
need the precise definition of $\rm{tg}$, but only the result of
Moore that it is zero for a split extension such as $\bbT^r\times
F\onto F$, so that the desired injectivity follows.

\quad\textbf{4.}\quad In view of the presentation
\[\bbZ\into n\bbZ\onto \bbZ/n\bbZ\]
any $2$-cocycle $\k:(\bbT^r\times F)\times (\bbT^r\times F)\to
\bbZ/n\bbZ$ lifts to a $2$-cochain $\k':(\bbT^r\times F)\times
(\bbT^r\times F)\to \bbZ$, whose coboundary now defines a
$3$-cocycle $d\k':(\bbT^r\times F)\times (\bbT^r\times F)\times
(\bbT^r\times F)\to n\bbZ$.  By Part 2 this is cohomologous as an
$n\bbZ$-valued $3$-cocycle to some cocycle depending only on the
coordinates in $F$: that is, there are a $2$-cochain
$\a:(\bbT^r\times F)\times (\bbT^r\times F)\to n\bbZ$ and a
$3$-cocycle $\b:F\times F\times F\to n\bbZ$ such that $d\k' = d\a +
\b$.

Therefore $\b = d(\k' - \a)$ is a $3$-cocycle depending only on
coordinates in $F$ that can be expressed as the coboundary of some
$\bbZ$-valued $2$-cochain on $\bbT^r\times F$, say $\xi_1 \in
\C((\bbT^r\times F)^2,\bbZ)$.  We will next show that $\xi_1$ can
also be taken to depend only on coordinates in $F$.

Using once again the presentation $\bbZ\into \bbR\onto \bbT$ and
Lemma~\ref{lem:vanishing-real-cohom} we see that $\b$ can
alternatively be expressed as the coboundary of some $\bbR$-valued
$2$-cochain on $F$, say $\xi_2 \in \C(F^2,\bbR)$. Now $d(\xi_2 -
\xi_1) = 0$, so $\xi_2 - \xi_1$ is an $\bbR$-valued $2$-cocycle on
$\bbT^r\times F$, so another appeal to the vanishing of real-valued
cohomology allows us to write it as $d\g_1$ for some Borel
$\g_1:\bbT^r\times F\to\bbR$. Recalling that $\xi_1$ is
$\bbZ$-valued, composing with the quotient map $\bbR\onto \bbR/\bbZ$
we deduce that $d(\g_1 + \bbZ) = \xi_2 + \bbZ$. Therefore the
$\bbT$-valued $2$-cocycle $\xi_2 + \bbZ$ on $F$ is a coboundary when
lifted to $\bbT^r\times F$, and so by Part 3 above it is actually a
coboundary among cochains that depend only on $F$. Letting $\g_2$ be
a cochain $F\to \bbR$ such that $d(\g_2 +\bbZ) = \xi_2 + \bbZ$, it
follows that we have $\b = d\xi_2 = d(\xi_2 - d\g_2)$ where $\xi_2 -
d\g_2$ takes values in $\bbZ$. Thus we have shown that $\b$ is
actually a $3$-coboundary for $\bbZ$-valued cochains depending only
on coordinates in $F$, and hence we can write $\b = d\k''$ for some
$\k'':F\times F\to \bbZ$.

Therefore $d(\k' - \a - \k'') = 0$, so now $\k' - \a - \k''$ is a
$\bbZ$-valued $2$-cocycle on $\bbT^r\times F$, and hence by
Lemma~\ref{lem:integer-cohom} there are some $\g \in
\widehat{\bbT^r\times F}$ and cochain $\rho:\bbT^r\times F\to \bbZ$
such that
\[(\k' - \a - \k'')(z,w) = d\rho(z,w) + \lfr\{\g(z)\} + \{\g(w)\}\rfr,\]
and so finally since $\a$ takes values in $n\bbZ$, passing back down
through the quotient map $\bbZ\onto \bbZ/n\bbZ$ we obtain
\[\k(z,w)= (\k'' + n\bbZ)(z,w) + d(\rho + n\bbZ)(z,w) + (\lfr\{\g(z)\} + \{\g(w)\}\rfr + n\bbZ).\]
Since $\k''$ depends only on coordinates in $F$ this is of the form
desired. \qed

\textbf{Remark}\quad For Part 2 above we made use of the injectivity
of certain inflation maps from $\rmH^\ast(F,\,\cdot\,)$ to
$\rmH^\ast(F\times H,\,\cdot\,)$ for a direct product group $F\times
H$. In the setting of finite groups $F$ and $H$ this simple result
can be proved by hand using the homogeneous bar resolution. However,
in the setting of measurable cohomology on non-finite groups this
approach runs into trouble because it relies on sampling cochains on
zero-measure subsets of the product group, and our cochains are only
defined up to negligible sets.  For this reason rigorous proofs
require some more careful machinery (particularly the
Hochschild-Serre spectral sequence), and take rather more work. \fin

\parskip 0pt

\bibliographystyle{abbrv}
\bibliography{bibfile}

\begin{thebibliography}{10}

\bibitem{Aus--cohomcty}
T.~Austin.
\newblock Continuity properties of {M}oore cohomology.
\newblock Preprint, available online at \verb|arXiv.org|: 0030818.

\bibitem{Aus--newmultiSzem}
T.~Austin.
\newblock Deducing the multidimensional {S}zemer\'edi {T}heorem from an
  infinitary removal lemma.
\newblock To appear, \emph{{J}. d'{A}nalyse {M}ath.}

\bibitem{Aus--ergdirint}
T.~Austin.
\newblock Extensions of probability-preserving systems by measurably-varying
  homogeneous spaces and applications.
\newblock Preprint, available online at \verb|arXiv.org|: 0905.0516.

\bibitem{Aus--lindeppleasant1}
T.~Austin.
\newblock Pleasant extensions retaining algebraic structure, {I}.
\newblock Preprint, available online at \verb|arXiv.org|: 0905.0518.

\bibitem{Aus--lindeppleasant2}
T.~Austin.
\newblock Pleasant extensions retaining algebraic structure, {II}.
\newblock Preprint, available online at \verb|arXiv.org|: 0910.0907.

\bibitem{Aus--nonconv}
T.~Austin.
\newblock On the norm convergence of nonconventional ergodic averages.
\newblock {\em Ergodic Theory Dynam. Systems}, 30(2):321--338, 2009.

\bibitem{BerLei02}
V.~Bergelson and A.~Leibman.
\newblock A nilpotent {R}oth theorem.
\newblock {\em Invent. Math.}, 147(2):429--470, 2002.

\bibitem{BerLei07}
V.~Bergelson and A.~Leibman.
\newblock Distribution of values of bounded generalized polynomials.
\newblock {\em Acta Math.}, 198(2):155--230, 2007.

\bibitem{ConLes84}
J.-P. Conze and E.~Lesigne.
\newblock Th\'eor\`emes ergodiques pour des mesures diagonales.
\newblock {\em Bull. Soc. Math. France}, 112(2):143--175, 1984.

\bibitem{ConLes88.1}
J.-P. Conze and E.~Lesigne.
\newblock Sur un th\'eor\`eme ergodique pour des mesures diagonales.
\newblock In {\em Probabilit\'es}, volume 1987 of {\em Publ. Inst. Rech. Math.
  Rennes}, pages 1--31. Univ. Rennes I, Rennes, 1988.

\bibitem{ConLes88.2}
J.-P. Conze and E.~Lesigne.
\newblock Sur un th\'eor\`eme ergodique pour des mesures diagonales.
\newblock {\em C. R. Acad. Sci. Paris S\'er. I Math.}, 306(12):491--493, 1988.

\bibitem{Fur77}
H.~Furstenberg.
\newblock Ergodic behaviour of diagonal measures and a theorem of {S}zemer\'edi
  on arithmetic progressions.
\newblock {\em J. d'Analyse Math.}, 31:204--256, 1977.

\bibitem{FurWei96}
H.~Furstenberg and B.~Weiss.
\newblock A mean ergodic theorem for
  $\frac{1}{N}\sum_{n=1}^{N}f({T}^nx)g({T}^{n^2}x)$.
\newblock In V.~Bergleson, A.~March, and J.~Rosenblatt, editors, {\em
  Convergence in Ergodic Theory and Probability}, pages 193--227. De Gruyter,
  Berlin, 1996.

\bibitem{HewRos79}
E.~Hewitt and K.~A. Ross.
\newblock {\em Abstract Harmonic Analysis, I (second ed.)}.
\newblock Springer, 1979.

\bibitem{HofMos73}
K.~H. Hofmann and P.~S. Mostert.
\newblock {\em Cohomology theories for compact abelian groups}.
\newblock Springer-Verlag, New York, 1973.
\newblock With an appendix by Eric C. Nummela.

\bibitem{HosKra01}
B.~Host and B.~Kra.
\newblock Convergence of {C}onze-{L}esigne averages.
\newblock {\em Ergodic Theory Dynam. Systems}, 21(2):493--509, 2001.

\bibitem{HosKra05}
B.~Host and B.~Kra.
\newblock Nonconventional ergodic averages and nilmanifolds.
\newblock {\em Ann. Math.}, 161(1):397--488, 2005.

\bibitem{HosKra07}
B.~Host and B.~Kra.
\newblock Uniformity seminorms on $\ell^\infty$ and applications.
\newblock Preprint, available online at \verb|arXiv.org|: 0711.3637, 2007.

\bibitem{KuiNie74}
L.~Kuipers and H.~Niederreiter.
\newblock {\em Uniform distribution of sequences}.
\newblock Wiley-Interscience [John Wiley \& Sons], New York, 1974.
\newblock Pure and Applied Mathematics.

\bibitem{Lei09}
A.~Leibman.
\newblock A canonical form and the distribution of values of generalized
  polynomials.
\newblock Preprint, available online at
  \verb|http://www.math.ohio-state.edu/~leibman/preprints/bas.pdf|, 2009.

\bibitem{Les93}
E.~Lesigne.
\newblock \'{E}quations fonctionnelles, couplages de produits gauches et
  th\'eor\`emes ergodiques pour mesures diagonales.
\newblock {\em Bull. Soc. Math. France}, 121(3):315--351, 1993.

\bibitem{Moo64(gr-cohomI-II)}
C.~C. Moore.
\newblock Extensions and low dimensional cohomology theory of locally compact
  groups. {I}, {II}.
\newblock {\em Trans. Amer. Math. Soc.}, 113:40--63, 1964.

\bibitem{Moo76(gr-cohomIII)}
C.~C. Moore.
\newblock Group extensions and cohomology for locally compact groups. {III}.
\newblock {\em Trans. Amer. Math. Soc.}, 221(1):1--33, 1976.

\bibitem{Moo76(gr-cohomIV)}
C.~C. Moore.
\newblock Group extensions and cohomology for locally compact groups. {IV}.
\newblock {\em Trans. Amer. Math. Soc.}, 221(1):35--58, 1976.

\bibitem{Rud93}
D.~J. Rudolph.
\newblock Eigenfunctions of {$T\times S$} and the {C}onze-{L}esigne algebra.
\newblock In {\em Ergodic theory and its connections with harmonic analysis
  ({A}lexandria, 1993)}, volume 205 of {\em London Math. Soc. Lecture Note
  Ser.}, pages 369--432. Cambridge Univ. Press, Cambridge, 1995.

\bibitem{Wei94}
C.~A. Weibel.
\newblock {\em An introduction to homological algebra}, volume~38 of {\em
  Cambridge Studies in Advanced Mathematics}.
\newblock Cambridge University Press, Cambridge, 1994.

\bibitem{Wig73}
D.~Wigner.
\newblock Algebraic cohomology of topological groups.
\newblock {\em Trans. Amer. Math. Soc.}, 178:83--93, 1973.

\bibitem{Zie07}
T.~Ziegler.
\newblock Universal characteristic factors and {F}urstenberg averages.
\newblock {\em J. Amer. Math. Soc.}, 20(1):53--97 (electronic), 2007.

\end{thebibliography}

\vspace{10pt}

\small{\textsc{Department of Mathematics, University of California,
Los Angeles CA 90095-1555, USA}}

\vspace{5pt}

\small{Email: \verb|timaustin@math.ucla.edu|}

\vspace{5pt}

\small{URL: \verb|http://www.math.ucla.edu/~timaustin|}

\end{document}